\documentclass[11pt, twoside]{article}
\usepackage{amsmath,amsfonts, amssymb, mathabx}
\usepackage{enumerate}

\usepackage{array}
\usepackage{color}

\usepackage{geometry}
\usepackage{hyperref}
\usepackage{fancyhdr}
\usepackage{tikz}

\usepackage{graphicx}

\begin{document}

\title{Nonstandard Universes}
\author{P. Ouwehand}\date{}
\maketitle
\tableofcontents

\abstract{These notes are concerned with the existence and the basic properties of the set-theoretic universes for nonstandard analysis, compiled by a beginner in the subject. It assumes a basic background in first-order logic, though the necessary material is revised in Appendix \ref{appendix_logic}. Needless to say, none of the material presented here is original, but has been adapted from the following sources: \cite{goldblatt1998lectures}, \cite{lindstrom1988invitation},\cite{chang1990model}\cite{albeverio2009nonstandard},\cite{loeb2015nonstandard},\cite{vath2007nonstandard}, \cite{lammert2015good}.}
\pagestyle{fancy}
\fancyhf{}
\pagenumbering{arabic}

\newtheorem{theorem}{Theorem}[section]
\newtheorem{lemma}[theorem]{Lemma}
\newtheorem{corollary}[theorem]{Corollary}
\newtheorem{proposition}[theorem]{Proposition}
\newtheorem{definition&proposition}[theorem]{Definition and Proposition}
\newtheorem{definition&theorem}[theorem]{Definition and Theorem}
\newtheorem{example}[theorem]{Example}
\newtheorem{examples}[theorem]{Examples}
\newtheorem{definition}[theorem]{Definition}
\newtheorem{remarks}[theorem]{Remarks}
\newtheorem{remark}[theorem]{Remark}
\newtheorem{convention}[theorem]{Convention}
\newtheorem{remarks*}[theorem]{Remarks$^*$}
\newtheorem{exercise}[theorem]{Exercise}
\newtheorem{exercise*}[theorem]{Exercise$^*$}
\newtheorem{exercises}[theorem]{Exercises}
\newtheorem{stelling}[theorem]{Stelling}
\newtheorem{teorema}[theorem]{Teorema}
\newtheorem{definisie}[theorem]{Definisie}
\newtheorem{voorbeelde}[theorem]{Voorbeelde}
\newtheorem{voorbeeld}[theorem]{Voorbeeld}
\newtheorem{opmerkinge}[theorem]{Opmerkinge}
\newtheorem{opmerking}[theorem]{Opmerking}
\newtheorem{oefeninge}[theorem]{Oefeninge}
\newtheorem{oefening}[theorem]{Oefening}

\def\<{\langle}
\def\>{\rangle}
\def\bproof{\noindent{\bf Proof: }}
\def\eproof{\begin{flushright}$\dashv$\end{flushright}}
\def\endbox{\begin{flushright}$\Box$\end{flushright}}
\def\ess{\text{\rm ess }}
\def\dom{\text{\rm dom}}
\def\ran{\text{\rm ran}}

\fancyhead[LE,RO]{\thepage}
\fancyhead[RE]{Nonstandard Frameworks}

\section{Universes}
\fancyhead[LO]{Universes}
\subsection{Definition of a Universe}
We work in  a set theory with atoms (individuals). In our universe $\mathbb U$, there will be two kinds of entities, namely {\em individuals} and {\em sets}.  Individuals are entities that contain no members --- at least, no members that belong to the universe. Thus from $\mathbb U$'s point of view, an element $a\in\mathbb U$ is an individual if and only if $a\cap\mathbb U=\varnothing$, yet $a\neq \varnothing$. Note that the empty set is not regarded as an individual, but as a set.  Sets are entities that are sets in the usual sense, but have the property that each of their elements also belong to the universe.

\begin{convention} 
\rm We will use lower case letters $a,b,c\dots$ to range over both sets and individuals, and reserve upper case letters $A,B,C,\dots$ to range over sets.\endbox
\end{convention}

\vskip0.3cm
\begin{definition}\rm (Transitive Set) 
A set $A\in\mathbb U$  is said to be {\em transitive} if and only if elements of elements of $A$ are elements of $A$:
\[\forall a\in A\;\forall b\in a\;(b\in A).\]
(This is vacuously true if $a\in A$  is an individual, as the quantifiers range over members of $\mathbb U$.)
Equivalently, if $a\in A$ and $a$ is a set, then $a\subseteq A$.
\endbox
\end{definition}

\begin{definition}\rm (Universe)
\begin{enumerate}[1.]\item
A {\em universe} $\mathbb U$ is a set with the following properties:
\begin{enumerate}[(a)]\item $\mathbb U$ is strongly transitive, i.e. for every $A\in\mathbb U$ there is a set $B\in\mathbb U$ such that $B$ is transitive and $A\subseteq B\subseteq \mathbb U$.
\item If $a,b\in \mathbb U$, then $\{a,b\}\in\mathbb U$.
\item If $A,B\in\mathbb U$, then $A\cup B\in\mathbb U$.
\item If $A\in\mathbb U$, then $\mathcal P(A)\in\mathbb U$, where $\mathcal P(A)$ is the powerset of $A$.
\end{enumerate}
\end{enumerate}
\endbox
\end{definition}
It follows directly from strong transitivity that a universe is transitive with respect to its set-members:  If $A\in\mathbb U$ is a set, then $A\subseteq \mathbb U$. We say that $\mathbb U$ is {\em transitive over sets}.

A set $T$ is called {\em supertransitive} if and only if whenever $A\in T$, then $A\cup\mathcal P(A)\subseteq T$. Clearly, a supertransitive set is transitive. 

\begin{lemma} If \;$\mathbb U$ is a universe, then every $A\in\mathbb U$ is an element of some supertransitive set $T\in\mathbb U$. 
\end{lemma}
\bproof
Suppose that $S\in\mathbb U$ is transitive, and define $T:=S\cup\mathcal P(S)$. Then $T\in\mathbb U$, by the definition of a universe. We claim that $T$ is supertransitive. Indeed if $A\in T$, then either $A\in S$ or $A\in\mathcal P(S)$. Either way, we see that $A\subseteq S$, and hence $A\subseteq T$. Now if $B\in\mathcal P(A)$, then $B\subseteq A\subseteq S$, so $B\in\mathcal P(S)\subseteq T$. Hence $\mathcal P(A)\subseteq T$. It follows that if $A\in T$, then $A\cup\mathcal P(A)\subseteq T$, so that $T$ is supertransitive.

Now suppose $A\in\mathbb U$. As $\mathbb U$ is strongly transitive, there is a transitive set $S\in\mathbb U$ such that $A\subseteq S$. Let $T=S\cup\mathcal P(S)$. Then $T$ is supertransitive and $A\in T$.
\eproof

\begin{definition}\rm (Universe over $X$)
If $X$ is a set, then $\mathbb U$ is said to be a {\em universe over $X$} if and only if $X\in\mathbb U$, $\varnothing\not\in X$, and no element of a member of $X$ belongs to $\mathbb U$, i.e. $\bigcup X\cap \mathbb U=\varnothing$.
Thus, from the point of view of $\mathbb U$, no member of $X$ has elements, yet none are the empty set --- they are individuals.
\endbox
\end{definition}
Note that if $\mathbb U$ is a universe over $X$, and $Y\subseteq X$, then $\mathbb U$ is also a universe over $Y$. Similarly, if $\mathbb U$ is a universe over $X$ for each $X\in\mathcal X\in\mathbb U$, then $\mathbb U$ is a universe over $\bigcup\mathcal X$.

\begin{example}\label{example_superstructures}\rm (Superstructures)\newline
Superstructures are the most common universes in practice.

Suppose that $X$ is a set.  The superstructure over $X$, denoted $V(X)$ is defined inductively as follows:
\begin{align*} V_0(X)&:=X\\
V_{n+1}(X)&:=V_n(X)\cup\mathcal P(V_n(X))\\
V(X)&:=\bigcup_{n<\omega}V_n(X)
\end{align*}
Clearly $V_{n}(X)\subseteq V_{n+1}(X)$ for all $n<\omega$, and $a\in V_{n+1}(X)$ if and only if $a\in V_n(X)$ or $a\subseteq V_n(X)$. 

It is easy to show by induction that $V_{n+1}(X)=X\cup\mathcal P(V_n(X))$: This is obvious in the case $n=0$. Next, suppose that $V_n(X)=X\cup\mathcal P(V_{n-1}(X))$. If $a\in V_{n+1}(X)$, then either (i) $a\in V_{n}(X)$  or (ii) $a\subseteq V_{n}(X)$ (or both). Hence either (i) $a\in X$ or $a\subseteq V_{n-1}(X)$, or (ii) $a\subseteq V_{n}(X)$. Since $V_{n-1}(X)\subseteq V_n(X)$ it follows that either (i) $a\in X$, or (ii) $a\subseteq V_{n}(X)$ (or both). It follows that $V_{n+1}(X)\subseteq X\cup \mathcal P(V_n(X))$. The reverse inclusion is obvious. 

In order for $V(X)$ to be a universe over $X$, one requirement is that the members of $X$ act like individuals, i.e. that $x\cap V(X)=\varnothing$ for all $x\in X$. The set $X$ is said to be  a {\em base set} if $\varnothing \not\in X$, and $\forall x\in X(x\cap V(X)=\varnothing)$.  Note that it  is always possible to replace a set $X$ by a base set of the same size: For example, given an infinite ordinal $\alpha$, pick a set $Y$ with the same cardinality as $X$ such that every element of an element of $Y$ has rank $\alpha$. One can always choose $\alpha$ sufficiently large so that this is possible. It is then easy to see by induction that each element of $V_n(Y)$ has  a rank $\beta$ where either $\beta<n$ or $\alpha<\beta\leq \alpha+n+1$. Now if $z\in y\in Y$, then $\text{rank}(z)=\alpha$, so $z\not\in V(Y)$, i.e. $y\cap V(Y)=\varnothing$. 

Now assume that $X$ is a base set. It is easy to see that each by induction $V_{n}(X)$ is transitive over sets: Certainly if $A\in V_{n+1}(X)=X\cap\mathcal P(V_{n}(X))$, then $A\subseteq V_n(X)$, since $A$ is a set, i.e. $A\not\in X$. It follows that $A\subseteq V_{n+1}(X)$. In particular, it follows that $V(X)$ is strongly transitive.

It is also easy to see that if $a,b\in V_n(X)$, then $\{a,b\}\in V_{n+1}(X)$. Next, if $A,B\in V_{n}(X)$, then $A\cup B\subseteq V_{n-1}(X)$, so $A\cup B\in V_{n}(X)$. Further, if $A\in V_n(X)$, then $A\subseteq V_{n-1}(X)$, so $\mathcal P(A)\subseteq\mathcal P(V_{n-1}(X))\subseteq V_n(X)$, and hence $\mathcal P(A)\in V_{n+1}(X)$. 

Hence if $X$ is a base set, then $V(X)$ is a universe over $X$. 

Observe that the sets $V_n(X)$ that make up a superstructure $V(X)$ over $X$ are supertransitive.

\endbox
\end{example}

\begin{remarks}\rm Suppose that $\mathbb U$ is a universe over $X$. In that case it easy to prove by induction  on $n$ that each $V_n(X)\subseteq\mathbb U$, from which it follows that $V(X)\subseteq\mathbb U$. Thus $V(X)$ is the smallest universe over $X$, assuming that one exists. 

Not every universe is of the form $V(X)$, however.  For example, consider $\mathbb U:=V_{\omega+\omega}$ in the usual cumulative hierarchy of sets. This is a universe over $\varnothing$, and there are no individuals, since if $y\in x\in\mathbb U$, then $y\in\mathbb U$. Yet $\mathbb U\neq V(\varnothing)$, since $V(\varnothing)=V_{\omega}$.
\endbox\end{remarks}

\subsection{Closure Properties of a Universe}

Observe the following closure properties of a universe, which are easy consequences of the definition of universe.  (Recall again our convention concerning upper  and lower case letters.)
\begin{enumerate}[1.]\item $a\in \mathbb U$ implies $\{a\}\in\mathbb U$.
\item $A_1,\dots,A_m\in\mathbb U$ implies $A_1\cup\dots\cup A_m\in\mathbb U$.
\item If $A\subseteq\mathbb U$ is finite, then $A\in\mathbb U$ .\newline(By 1., 2.)
\item If $A\subseteq B$ and $B\in\mathbb U$, then $A\in\mathbb U$.\newline
(For $A\in\mathcal P(B)\in\mathbb U$, and $\mathbb U$ is transitive.)
\item If $\{A_i:i\in I\}\subseteq B\in\mathbb U$, then $\bigcup_{i\in I} A_i\in\mathbb U$.\newline
(For there is transitive $T\in\mathbb U$ such that $B\subseteq T$. Then each $A_i\in T$, so each $A_i\subseteq T$, and hence $\bigcup_{i\in I}A_i\subseteq T$. Now apply 4.)
\item If $B=\{A_i:i\in I\}\in\mathbb U$, then $\bigcup B= \bigcup_{i\in I}A_i \in\mathbb U$.\newline
(Follows directly from 5.)
\item If $\{A_i:i\in I\}\subseteq \mathbb U$ is a non-empty family  of sets, then $\bigcap_{i\in I} A_i\in\mathbb U$.\newline (The intersection is a subset of $A_{i_0}$, for any $i_0\in I$. Now apply 4.)
\item $a,b\in \mathbb U$ implies $(a,b)\in \mathbb U$. More generally, if $a_1,\dots, a_m\in \mathbb U$, then $(a_1,\dots,a_m)\in\mathbb U$.\newline (Because $(a,b):=\{\{a\},\{a,b\}\}$. Then $(a,b,c):=((a,b),c)$, etc.)
\item If $A,B\in\mathbb U$ and $R\subseteq A\times B$, then $R\in\mathbb U$. More generally, if $A_1,\dots, A_m\in\mathbb U$ and $R\subseteq A_1\times \dots \times A_m$, then $R\in \mathbb U$.
\newline(Since $R\subseteq \mathcal P\mathcal P(A\cup B)$, the result follows from  2. and  4.)
\item If $R\in\mathbb U$ is a binary relation, then $\dom(R), \ran(R), R^{-1}\in\mathbb U$. Furthermore, if $C\subseteq \dom(R)$, then $R[C]\in\mathbb U$.\newline(For $\dom(R),\ran(R)\subseteq \bigcup\bigcup R\in\mathbb U$, so by 4. we have that   $\dom(R),\ran(R)\in\mathbb U$. Also $R^{-1}\subseteq\ran(R)\times\dom(R)$, so that $R^{-1}\in\mathbb U$ follows by 9. Finally $R[C]\subseteq \ran(R)\in\mathbb U$.)
\item If $A,B\in \mathbb U$ and $f:A\to B$, then $f\in \mathbb U$. Furthermore, if $A'\subseteq A, B'\subseteq B$, then $f[A'], f^{-1}[B']\in\mathbb U$.
\newline(First, $f\subseteq A\times B\in\mathbb U$. Then $f[A']\subseteq B\in\mathbb U$ and $f^{-1}[B']\subseteq A\in\mathbb U$.)
\item If $A,B\in\mathbb U$, then $B^A\in\mathbb U$, where $B^A$ is the set of all functions from $A$ to $B$..\newline
(For $B^A\subseteq \mathcal P(A\times B)$.)
\item If $\{A_i:i\in\mathbb I\}\in\mathbb U$,and $I\in\mathbb U$, then $\prod_{i\in I}A_i\in\mathbb U$.\newline
(Because $\prod_{I}A_i\subseteq \;(\bigcup_{ I} A_i)^I$.)
\end{enumerate}

\subsection{Nonstandard Embeddings}
We assume some familiarity with basic first-order logic, including the basics of model theory. Refer to Appendix \ref{appendix_logic} for a quick reminder of the basic notions used below.

Let $\mathcal L_\in$ denote  a first-order language with equality $=$ and a single binary relation symbol $\in$.   We assume that there is a countable collection of variables, and take as basic propositional connectives the connectives $\lnot$ (not) and $\land$ (and), and as basic quantifier $\forall$ (for all). The other connectives $\lor$ (or), $\to$ (then), $\leftrightarrow$ (if and only if) are defined in terms of $\land,\lnot$ in the usual way, and the existential quantifier is defined in terms of $\lnot,\forall$ in the usual way.

In addition, $\exists!y\;\psi(y)$ abbreviates the formula $\exists y\;(\psi(y)\land\forall z\;(\psi(z)\to z=y))$, which states that $y$ is the unique element for which $\psi$ holds. 
\vskip0.3cm

In nonstandard analysis, the following types of $\mathcal L_\in$-formula play a central role:
\begin{definition}\rm (Bounded Formula)
A {\em bounded} $\mathcal L_\in$-formula is an $\mathcal L_\in$-formula all of whose quantifiers are bounded, i.e. of the form $\forall x\in y$ or $\exists x\in y$, where $\forall x\in y\;\varphi(x,y)$ is an abbreviation of $\forall x\;(x\in y\to\varphi(x,y))$, and $\exists x\in y\;\varphi(x,y)$ abbreviates $\exists x\;(x\in y\land\varphi(x,y))$.
\endbox\end{definition}

\begin{definition}\label{defn_transfer_map}\rm (Transfer Map)\\
A {\em transfer map} for a set $X$ is a function $*:\mathbb U\to\mathbb V$ between two universes $\mathbb U,\mathbb V$ with the properties that:
\begin{enumerate}[1.]\item $\mathbb U$ is a universe over $X$.
\item ${}^*x=x$ for every $x\in X$, and ${}^*\varnothing=\varnothing$
\item {\bf Transfer:} $*$ is a {\em bounded elementary embedding}, i.e. if $\varphi(x_1,\dots, x_n)$ is an $\mathcal L_\in$-formula, and $u_1,\dots, u_n\in\mathbb U$, then\[\mathbb U\vDash \varphi[u_1,\dots, u_n] \quad \text{iff}\quad \mathbb V\vDash \varphi[{}^*u_1,\dots,{}^*u_n].\]
\end{enumerate}
If $\mathbb U=V(X),\mathbb V=V(Y)$ are superstructures over $X,Y$ respectively, then it is usually also required that ${}^*X=Y$.
\endbox
\end{definition}
The transfer property in the preceding definition states that the $*$-map transfers properties  that are definable by bounded formulas from $\mathbb U$ to $\mathbb V$.

\vskip0.3cm\noindent {\bf Notation:} For $A\in\mathbb U$, define:\[{}^\sigma A:=*[A]:=\{{}^*a: a\in A\}\] to be the image of the set $A$ under the $*$-map. 

\begin{definition}\label{defn_nonstandard_embedding}\rm (Nonstandard Framework/Nonstandard Embedding)\\
A transfer map  $*:\mathbb U\to\mathbb V$ for $X$ is a {\em nonstandard framework}, or a {\em nonstandard embedding},  if there is a countable $C\in\mathbb U$ such that ${}^\sigma C$ is a proper subset of ${}^*C$, i.e. ${}^\sigma C\subsetneqq {}^*C$.
\endbox
\end{definition}

Note that if $*$ is a transfer map, then always ${}^\sigma A\subseteq {}^*A$, as $\mathbb U\vDash a\in A$ implies $\mathbb V\vDash {}^*a\in {}^*A$. If $A$ is finite, then we shall see that ${}^\sigma A={}^*A$. If, however, the $*$-map is a nonstandard framework --- so that $C$ is a proper subset of ${}^*C$ for some countable $C\in\mathbb U$--- then it will transpire that ${}^\sigma A\subsetneqq {}^*A$ whenever $A\in\mathbb U$ is infinite. In that case, we can think of ${}^*A$ as a version of $A\in\mathbb U$ that lives in $\mathbb V$ --- in that $A$ and ${}^*A$ satisfy the same bounded sentences --- but  where ${}^*A$ has additional  elements that do not correspond to members of $A$.

\begin{convention}\label{convention_expanded_language} \rm The language $\mathcal L_\in$ in which we work has no constant symbols. However, in the interests of brevity we will often write formulas as if there are constant symbols for every member $a\in\mathbb U$.
For example, when we write \[\mathbb U\vDash \exists x\in A\;\forall y\in B\;(c\in y\land \psi(x,y, d)),\qquad (\text{where }A,B,c, d\in\mathbb U)\] this should be taken to mean
\[\mathbb U\vDash \varphi[A,B,c, d],\quad\text{where }\varphi(u,v,w,t)\equiv\exists x\in u\;\exists y\in v\;(z\in y\land\psi(x,y, t)).\]
By transfer, we then have\[\mathbb V\vDash \varphi[{}^*A,{}^*B,{}^*c, {}^*d]\quad\text{i.e.}\quad \mathbb V\vDash\exists x\in {}^*A\;\forall y\in {}^*B\;({}^*c\in y\land\psi(x,y, {}^*d)).\]
Thus without loss of generality, we may assume that, when working with a transfer map $*:\mathbb U\to\mathbb V$, the language $\mathcal L_\in$ is expanded  to a language --- denoted $\mathcal L_\mathbb U$ --- which has a constant symbol $c_a$ for  every entity $a\in \mathbb  U$. Naturally, the constant symbol $c_a$ is to be interpreted as the  entity $a$ in the model $\mathbb U$. If $*:\mathbb U\to\mathbb V $ is a transfer map, then $c_a$ will be interpreted as ${}^*a$ in $\mathbb V$.  When we replace all occurrences of these constants in a formula $\varphi$ by their $*$-value, we obtain the $*$-transform ${}^*\varphi$ of the formula. Thus, for example\[{}^*\Big(\exists x\in A\;\forall y\in B\;(c\in y\land \psi(x,y, d))\Big)\equiv \exists x\in {}^*A\;\forall y\in {}^*B\;({}^*c\in y\land\psi(x,y, {}^*d)).\]
It is not hard to see how to define ${}^*\varphi$ for formulas $\varphi$ by induction on the  complexity of $\varphi$. 

The transfer property is then easily seen to be equivalent to the following: If $\varphi$ is a bounded {\em sentence} of $\mathcal L_\mathbb U$, then $\mathbb U\vDash\varphi$ if and only if $\mathbb V\vDash {}^*\varphi$. \endbox
\end{convention}

\begin{lemma} \label{lemma_V_over_*X} Suppose that $*:\mathbb U\to\mathbb V$ is a transfer map for $X$. Then $\mathbb V$ is a universe over ${}^*X$.
\end{lemma}
\bproof
We need only show that $\varnothing\not \in {}^*X$ and that $\bigcup{}^*X\cap\mathbb V=\varnothing$. Since ${}^*\varnothing=\varnothing$, we see that $\varnothing\not\in X$ transfers to $\varnothing\not\in{}^*X$. Next, if $\bigcup{}^*X\cap\mathbb V\neq \varnothing$, then $\mathbb V\vDash \exists x\in {}^*X\;\exists y\in x\;(y=y)$ (which simply says that there is an $x\in {}^*X$ which has an element $y\in\mathbb V$). Then transfer implies that  $\bigcup X\cap\mathbb U\neq \varnothing$, which contradicts the fact that $\mathbb U$ is a universe over $X$.
\eproof

In the next section, we will discuss some of the properties of nonstandard frameworks, assuming that they exist. The question of existence is dealt with in Section \ref{section_existence_nonstandard_frameworks}, but, assuming familiarity with the model-theoretic concepts in Appendix  \ref{appendix_logic}, this can be read now.
We'll have more to say about Definition \ref{defn_nonstandard_embedding} later on. For now, note that this property is essential for nonstandard analysis to have any real power via the introduction of nonstandard objects.

\section{Properties of the $*$-map}
\fancyhead[LO]{Properties of the $*$-map}
In this section, we assume that $*:\mathbb U\to\mathbb V$ is a transfer map between two universes $\mathbb U,\mathbb V$.

Note that the $*$-map is injective: For if $\mathbb V\vDash {}^*a={}^*b$, then $\mathbb U\vDash a=b$, by the transfer property.

Further observe that $*$ maps individuals in $\mathbb U$ to individuals in $\mathbb V$: For suppose that $a\in\mathbb U$ is an individual. Then $\mathbb U\vDash \lnot\exists x\in a\;(x=x)$, i.e. $a$ has no elements (in common with $\mathbb U$). By transfer, $^*a$ has no elements (in common with $\mathbb V$). Since $^*$ is injective, $^*a\neq {}^*\varnothing$. Now $^*\varnothing =\varnothing$, by definition of the $*$-map. Thus $^*a\in\mathbb V$ is an element which has no members (in common with $\mathbb V$), yet is not the empty set, i.e. $^*a$ is an individual in $\mathbb V$.

The following lemma shows that certain basic operations can be expressed by bounded $\mathcal L_\in$-formulas $\varphi_n$. In what follows, we will be able to improve readability by using these abbreviations instead of the $\varphi_n$ inside bounded $\mathcal L_\in$-formulas.
\begin{lemma} \label{lemma_define_simple_notions} Let $\mathbb U$ be a universe. There are bounded formulas $\varphi_0,\dots,\varphi_7$ such that for all elements $a_n\in\mathbb U$ and all {\bf sets} $A_n\in\mathbb U$, the following hold:
\begin{enumerate}[(a)]\item $A_1=\varnothing$ if and only if $\mathbb U\vDash\varphi_0[A_1]$.
\item $A_1=\{a_1,\dots,a_n\}$ if and only if $\mathbb U\vDash\varphi_{1,n}[A_1,a_1,\dots,a_n]$.
\item $A_1=(a_1,\dots, a_n)$ if and only if $\mathbb U\vDash\varphi_{2,n}[A_1,a_1,\dots,a_n]$.
\item $A_1\subseteq A_2$ if and only if $\mathbb U\vDash \varphi_3[A_1,A_2]$.
\item $A_1=A_2\times A_3$ if and only if  $\mathbb U\vDash \varphi_4[A_1,A_2, A_3]$.
\item $A_1:A_2\to A_3$ if and only if $\mathbb U\vDash \varphi_5[A_1,A_2, A_3]$.
\item If $\mathbb U$ is a universe over $X$, so that $V(X)\subseteq\mathbb U$, there is $\varphi_{6,n}$ such that $a_1\in V_n(X)$ if and only if $\mathbb U\vDash\varphi_{6,n}[X,a_1]$.
\item If $\mathbb U$ is a universe over $X$, so that $V(X)\subseteq\mathbb U$, there is $\varphi_{7,n}$ such that $A_1$ is a {\em set} in $V_n(X)$ if and only if $\mathbb U\vDash\varphi_{7,n}[X,A_1]$.
\end{enumerate}\end{lemma}
\bproof
(a) Take $\varphi_0(x)\equiv \forall y\in x\;(y\neq y)$. (Recall that $A_1$ is required to be a set.)

(b) Take $\varphi_{1,n}(x, y_1,\dots, y_n)$ to be the formula \[y_1\in x\land y_2\in x\land \dots\land y_n \in x \land \forall y\in x\;(y=y_1\lor y= y_2\lor \dots\lor y=y_n).\]

(c) For $n=2$, we see that $A_1=(a_1,a_2)=\{\{a_1\},\{a_1,a_2\}\}$ if and only if \[\mathbb U\vDash \exists x\in A_1\;\exists y\in A_1\;(A_1=\{x,y\} \land x=\{a_1\} \land y=\{a_1,a_2\}).\] The formulas inside the brackets are abbreviations of $\varphi_1$. This defines a bounded formula $\varphi_{2,2}$ for the case $n=2$. Then for $n=3$, we can proceed in a similar way to define $\varphi_{2,3}$, using the just defined $\varphi_{2,2}$, and the  fact that $(a_1,a_2,a_3):=((a_1,a_2),a_3)$, etc.

(d) Take $\forall u\in x\;(u\in y)$ for $\varphi_3(x,y)$.

(e) Note that $A_1=A_2\times A_3$ if and only if \[\forall u\in A_1\;\exists v\in A_2\;\exists w\in A_3\;(u=(v,w))\land\forall v\in A_2\;\forall w\in A_3\;\exists u\in A_1\;(u=(v,w)),\] where statements of the form $x=(y,z)$ are to be replaced by versions of $\varphi_2(x,y,z)$. From here, $\varphi_4(x,y,z)$ is apparent.

(f) Note that $A_1$ is a function from $A_2$ to $A_3$ if and only if \[A_1\subseteq A_2\times A_3 \land \forall v\in A_2\;\exists!w\in A_3\;\exists u\in A_1\;(u=(v,w)), \] where $\exists!y\in x\;\psi(x,y)$ abbreviates the bounded formula $\exists y\in x\;(\psi(x,y)\land\forall z\in x\;(\psi(x,z)\to z=y))$, which states that $y$ is the unique member of $x$ for which $\psi(x,y)$ holds. The required formula $\varphi_5$ can now be constructed easily with the aid of $\varphi_3,\varphi_4$.


(g) This is proved by induction on $n$. For the case $n=0$, note that $V_0(X)=X$, so take $\varphi_{6,0}(x,y)\equiv y\in x$. Then if  $a_1\in V_1(X)=X\cup\mathcal P(X)$, necessarily $a_1\in X$ or $a_1\subseteq X$, so take $\varphi_{6,1}(x,y)\equiv y\in x\lor \forall z\in y\;(z\in x)$. Similarly, note that if $a_1\in V_{n+1}(X)=X\cup \mathcal P(V_n(X))$, then $a_1\in X$ or $a_1\subseteq V_n(X)$, so define $\varphi_{6,n+1}(x,y)\equiv y\in x\lor\forall z\in y\;\varphi_{6,n}(x,z)$.

(h) Take $\varphi_{7,n}(x,y)$ to be $\varphi_{6,n}(x,y)\land y\not\in x$. Then $\varphi_{7,n}(X,A_1)$ holds if and only if $A_1$ is a {\em set} in $V_n(X)$.
\eproof


\begin{theorem}\label{thm_*_comprehension} {\rm ($*$-Comprehension)} Let $\varphi(y,x_1,\dots, x_n)$ be a bounded $\mathcal L_\in$-formula. For all $u_1,\dots, u_n, a\in\mathbb U$,
\[{}^*\{y\in a: \mathbb U\vDash \varphi(y,u_1,\dots,u_n)\}=\{y\in {}^*a: \mathbb V\vDash \varphi(y,{}^*u_1,\dots,{}^*u_n)\}.\]
\end{theorem}

\bproof By definition of the $*$-map, $^*\varnothing=\varnothing$, so the result is true if $a\in\mathbb U$ is an individual.  Assume therefore that $a$ is a set in $\mathbb U$. Define\[B:=\{y\in a: \mathbb U\vDash \varphi(y,u_1,\dots,u_n)\}.\] Then $B\in\mathcal P(a)\subseteq\mathbb U$, and \[ \mathbb U\vDash\forall y\in a\;(y\in B \leftrightarrow \varphi(y,u_1,\dots, u_n)).\] By transfer,  \[ \mathbb V\vDash\forall y\in {}^*a\;(y\in {}^*B \leftrightarrow \varphi(y, {}^*u_1,\dots, {}^*u_n)),\] from which it follows that ${}^*B\cap{}^*a=\{y\in{}^*a: \varphi(y, {}^*u_1,\dots, {}^*u_n)\}$. But $B\subseteq a$, so ${}^*B\subseteq{}^*a$, again by transfer, since the formula $B\subseteq a$ is clearly equivalent to a bounded formula of $\mathcal L_\in$. Thus the result follows.\eproof

Observe that if $*:\mathbb U\to\mathbb V$  is a transfer map for a set $X$, then:
\begin{itemize}\item $a\in B$ if and only if ${}^*a\in{}^*B$ and $a=b$ if and only if ${}^*a={}^*b$.
\item $A\subseteq B$ if and only if ${}^*A\subseteq{}^*B$.
\item If $A\subseteq X$, then $A\subseteq{}^*A\subseteq{}^*X$. In particular, $X\subseteq{}^*X$: This is because ${}^*x=x$ for $x\in X$.
\item ${}^*(A\cup B)={}^*A\cup{}^*B$, ${}^*(A\cap B)={}^*A\cap{}^*B$, ${}^*(A - B)={}^*A -{}^*B$: For example, by the preceding Lemma, ${}^*(A\cup B)={}^*\{x\in A\cup B: x\in A\lor x\in B\}=\{x\in {}^*(A\cup B): x\in {}^*A\lor x\in {}^*B\}={}^*(A\cup B)\cap ({}^*A\cup{}^*B)$. Since ${}^*A,{}^*B\subseteq {}^*(A\cup B)$, it follows that ${}^*(A\cup B)={}^*A\cup{}^*B$.
\item If $A=\{a_1,\dots,a_n\}$ is a finite set, then ${}^*A={}^\sigma A=\{{}^*a_1,\dots, {}^*a_n\}$: This is because $A=\{x\in A: x=a_1\lor\dots\lor x=a_n\}$.
\item If $A$ is a transitive set, then ${}^*A$ is a transitive set: For $A$ is transitive if and only if the bounded formula $\forall x\in A\;\forall y\in x\;(y\in A)$ holds in $\mathbb  U$.
\item ${}^*\mathcal P(A)\subseteq\mathcal P({}^*A)$: This follows by transfer of the formula $\forall X\in\mathcal P(A)\;(X\subseteq A)$. Hence every element $X$ of ${}^*\mathcal P(A)$ is a subset of ${}^*A$ and thus in $\mathcal P({}^*A)$. 
\item $^*(a_1,\dots, a_n)=({}^*a_1,\dots, {}^*a_n)$:  This follows from Lemma \ref{lemma_define_simple_notions}(c).
\item ${}^*(A_1\times \dots\times A_n)={}^*A_1\times\dots\times {}^*A_n$: This follows from Lemma \ref{lemma_define_simple_notions}(e).
\item If $R$ is an $n$-ary relation, so is ${}^*R$: For then $R\subseteq A_1\times \dots\times A_n$ for some sets $A_1,\dots, A_n\in\mathbb U$.
\item If $R$ is a binary relation, then
\begin{itemize}\item ${}^*\dom(R)=\dom ({}^*R)$: For if $x\in\dom(R)$, then $(x,y)=\{\{x\},\{x,y\}\}$ belongs to $R$ for some $y$. Hence, with $a:=\{\{x\}, \{x,y\}\}$, $b:=\{x\}, c:=\{x,y\}$, we have $x\in \dom(R)$ iff
\[\exists a\in R\;\exists b\in a\;\exists c\in a\;(a=\{b,c\}\land b=\{x\}\land\exists y\in c\;(c=\{x,y\})).\] The result now follows by Lemma \ref{lemma_define_simple_notions}(b). (Or, write it all out:
\[\aligned \exists a\in R\;\exists b\in a\;\exists c\in a\;&(\forall d\in a\;(d=b\lor d=c)\land x\in b\land\forall w\in b\;(w=x)\\&\land x\in c\land \exists y\in c\;\forall w\in c(w=x\lor w=y)),\endaligned\] a bounded formula.)
\item ${}^*\ran(R)=\ran(^*R)$: The proof is very similar to the preceding.
\item ${}^*(R^{-1})=({}^*R)^{-1}$: Choose $T\in\mathbb U$ transitive so that $R\cup R^{-1}\subseteq T$. Then if $(x,y)\in R$, it follows that $x,y\in T$. Then \[\mathbb U\vDash \forall x\in T\;\forall y\in T\;((x,y)\in R^{-1} \leftrightarrow (y,x)\in R),\] and hence \[\mathbb V\vDash \forall x\in {}^*T\;\forall y\in {}^*T\;((x,y)\in {}^*(R^{-1}) \leftrightarrow (y,x)\in {}^*R),\] using Lemma \ref{lemma_define_simple_notions}(c). As ${}^*R\cup {}^*(R^{-1})\subseteq {}^*T$, it follows that ${}^*(R^{-1})= ({}^*R)^{-1}$.
\item If $C\subseteq\dom(R)$, then ${}^*(R[C])= {}^*R[{}^*C]$: Consider the bounded formula \\$\forall y\in\ran(R)\;(y\in R[C]\leftrightarrow\exists x\in C\;((x,y)\in R))$.
\item If $D\subseteq \ran(R)$, then ${}^*(R^{-1}[D])= ({}^*(R^{-1}[{}^*D])$: For $\ran(R)=\dom(R^{-1})$.
\item If $R,S$ are binary relations, then ${}^*(R\circ S)={}^*R\circ{}^*S$:  Choose a transitive $T\in\mathbb U$ such that $R\cup S\subseteq T$. Then use transfer on the bounded formula\[\forall x\in T\;\forall z\in T\;((x,z)\in R\circ S\leftrightarrow \exists y\in T\;((x,y)\in S\land (y,z)\in R)).\]
\end{itemize}
\item If $f:A\to B$ is a function in $\mathbb U$, then ${}^*f:{}^*A\to{}^*B$ is a function in $\mathbb V$. Moreover, ${}^*(f(a))={}^*f({}^*a)$ for $a\in A$, and ${}^*f$ is injective/surjective if and only if $f$ is injective/surjective: Since $f$ is a binary relation with $\dom(f)=A, \ran(f)\subseteq B$, we immediately see that ${}^*f$ is a binary relation with $\dom({}^*f)={}^*A, \ran({}^*f)\subseteq {}^*B$. Since $\forall a\in A\exists!b\in B\;((a, b)\in f)$, it follows by transfer that ${}^*f$ is a function. Furthermore, transfer of  $\forall a\in A\forall b\in B\;((a,b)\in f \leftrightarrow b=f(a))$ now leads to $({}^*f)({}^*a)={}^*f(a)$. Next, $f$ is injective if and only if the bounded formula $\forall a_1\in A\;\forall a_2\in A\;(\exists b\in B\;((a_1,b)\in f\land (a_2,b)\in f)\rightarrow a_1=a_2)$ is satisfied, from which it follows by transfer ${}^*f$ is injective if and only if $f$ is. Finally $f$ is surjective if and only if $\forall b\in B\;\exists a\in A\;((a,b)\in f)$.
\end{itemize}

\begin{lemma} If $*:\mathbb U\to\mathbb V$  is a transfer map for a set $X$ and $V_{n}(X)\subseteq B\in\mathbb U$, then ${}^*V_n(X)=V_n({}^*X)\cap {}^*B$.
\end{lemma}
\bproof By Lemma \ref{lemma_define_simple_notions}(g), we have \[V_n(X)=V_n(X)\cap B=\{x\in B:\mathbb U\vDash \varphi_{6,n}(X,x)\},\] and thus \[{}^*V_n(X)=\{x\in {}^*B:\mathbb V\vDash \varphi_{6,n}({}^*X,x)\}=V_n({}^*X)\cap {}^*B.\] \eproof

The following result is crucial, as we wish to be able to talk about structures in our universes. For example, suppose that $*:V(\mathbb R)\to V({}^*\mathbb R)$ is a transfer map over the set $\mathbb  R$ between two superstructures. When talking about the reals, we also wish to take into account the operations and relations $+,\cdot,-, {}^{-1}, \leq$, so that we want to be able to talk about the model $(\mathbb R,+,\cdot, -, {}^{-1},\leq)$. This will transfer to a model $({}^*\mathbb R,{}^*+,{}^*\cdot, {}^*-, {}^*{}^{-1},{}^*\leq)$. The following result shows that these two models are elementarily equivalent, i.e. that they satisfy the same first-order sentences. Indeed, the $*$-map induces an elementary embedding from the one into the other:

\begin{theorem} Suppose that $*:\mathbb U\to\mathbb V$ is a transfer map for $X$, and that $A\in\mathbb U$. Let $\mathfrak{A}=(A, \mathcal L^\mathfrak A)$ be a model of a first-order language $\mathcal L$ (which need not be the language $\mathcal L_\in$ of $\mathbb U$). Then the restriction $*\restriction_A$ of $*$ to $A$ is an elementary embedding from $\mathfrak{A}$ into ${}^*\mathfrak{A}$.
\end{theorem}

\bproof To simplify notation, suppose that $\mathfrak{A}=(A,R^\mathfrak A)$ where $R\in\mathcal L$ is a single binary relation symbol. Then $R^\mathfrak A\subseteq A\times A$, and thus $R^\mathfrak A\in\mathbb U$. Hence $\mathfrak{A}\in\mathbb U$, and ${}^*\mathfrak{A}=({}^*A, {}^*R^\mathfrak A)$.  Choose $T\in\mathbb U$ transitive so that $\mathfrak{A}\in T$. We will show by induction on complexity that for every $\mathcal L$-formula $\varphi(x_1,\dots, x_n)$ there is a bounded $\mathcal L_\in$-formula $\bar{\varphi}$ such that\[\mathfrak{A}\vDash \varphi[a_1,\dots, a_n]\qquad\text{if and only if}\qquad \mathbb U\vDash\bar{\varphi}[a_1,\dots,a_n,\mathfrak A,T] \qquad\text{(for $a_1,\dots,a_n\in A$)}.\]
Suppose first that $\varphi(x,y)\equiv R(x,y)$ is atomic. Define the bounded $\mathcal L_\in$ formula $\bar{\varphi}$ by  \[\bar{\varphi}(x,y,z,t)\equiv\exists u\in t\;\exists v\in t\;\exists w\in t\;\Big(w=(x,y)\land z= (u,v)\land w\in v\Big).\]
Then $\bar{\varphi}[a_1,a_2,\mathfrak A,T]$ asserts that there are $u,v,w\in T$ such that $w=(a_1,a_2), z=(A,R^\mathfrak A)$ and $w\in R^\mathfrak A$.
Now clearly $\mathfrak{A}\vDash R[a_1,a_2]$ if and only if  $\mathbb U\vDash\exists u\in T\;\exists v\in T\;\exists w\in T\;\Big(w=(a_1,a_2)\land (u,v)=\mathfrak A\land w\in v\Big)$, i.e. if and only if $\mathbb U \vDash \bar{\varphi}[a_1,a_2,\mathfrak A,T]$. This deals with the case where $\varphi$ is an atomic $\mathcal L$ formula. The propositional connectives are easily handled. If $\varphi(x_1,\dots, x_n)\equiv\exists z\;\psi(z,x_1,\dots, x_n)$, then 
 \[\aligned\mathfrak{A}\vDash \varphi[a_1,\dots, a_n] &\Leftrightarrow \mathfrak{A}\vDash\psi[b,a_1,\dots,a_n]\quad\text{some $b\in A$},\\
  &\Leftrightarrow\mathbb U\vDash \bar{\psi}[b,a_1,\dots, a_n,\mathfrak{A},T]\quad\text{some $b\in A$, by induction hypothesis},\\
 &\Leftrightarrow\mathbb U\vDash\underbrace{\exists u\in T\;\exists v\in T\;\exists z\in T\;((u,v)=\mathfrak{A}\land z\in u\land \bar{\psi}(z,a_1,\dots, a_n, \mathfrak{A},T))}_{\bar{\varphi}[a_1.\dots,a_n,\mathfrak{A},T]}
 \endaligned\]
 This completes the induction.
 
 Since in $\mathbb V$ we have that ${}^*\mathfrak{A}\in{}^*T$ and ${}^*T$ is transitive, we obtain similarly that \[{}^*\mathfrak{A}\vDash \varphi[a_1,\dots, a_n]\qquad\text{if and only if}\qquad \mathbb V\vDash\bar{\varphi}[a_1,\dots,a_n,{}^*\mathfrak A,{}^*T] \qquad\text{(for $a_1,\dots,a_n\in {}^*A$)}\] for every $\mathcal L$ formula $\varphi$. But since $\bar{\varphi}$ is then a bounded $\mathcal L_\in$-formula and $*$ is a transfer map, we have 
\[\mathbb U\vDash\bar{\varphi}[a_1,\dots,a_n,\mathfrak A,T] \qquad\text{if and only if}\qquad \mathbb V\vDash\bar{\varphi}[{}^*a_1,\dots,{}^*a_n,{}^*\mathfrak A,{}^*T],\] from which it follows that
\[\mathfrak{A}\vDash \varphi[a_1,\dots, a_n]\qquad\text{if and only if}\qquad {}^*\mathfrak{A}\vDash \varphi[{}^*a_1,\dots, {}^*a_n].\] Thus $*:\mathfrak{A}\to{}^*\mathfrak{A}$ is an elementary embedding.
 \eproof

\section{Standard, Internal and External Objects} 
\fancyhead[LO]{Standard, Internal and External Objects}

In this section, we suppose  that we're working with a transfer map $*:\mathbb U\to\mathbb V$ for a base set $X\in\mathbb U$.
\subsection{Standard and Internal Objects}
\begin{definition}\rm \begin{enumerate}[(a)]\item An element $v\in\mathbb V$ is said to be {\em standard} if and only if there is $u\in\mathbb U$ such that ${}^*u=v$. Thus the standard objects in $\mathbb V$ are those in the range of the $*$-map.\\ We denote the set of standard objects in $\mathbb V$ by ${}^\sigma\mathbb U$.
\item An element $v\in \mathbb V$ is said to be {\em internal} if and only if there is $A\in\mathbb U$ such that $v\in {}^*A$. Thus the internal objects in $\mathbb V$ are those that belong to a standard set.\\
We denote the set of internal objects in $\mathbb V$ by ${}^*\mathbb U$.
\item Sets in $\mathbb V$ which are not internal are said to be {\em external}.
\end{enumerate}\endbox\end{definition}

\begin{remarks}\rm
\begin{enumerate}[(a)]
\item ${}^*\mathbb U\subseteq \mathbb V$ is transitive, i.e. if $A\in\mathbb V$ is internal, and $a\in A$, then $a$ is internal: Since $A$ is internal, there is $U\in\mathbb U$ such that $A\in {}^*U$. By definition of universe, there is a transitive $T\in\mathbb U$ such that $U\subseteq T$, and hence $A\in {}^*T$. Since ${}^*T$ is also transitive, we have $a\in {}^*T$, which shows that $A$ is internal.
\item Note that since ${}^*x=x$ for every member $x\in X$ of the base set, each $x\in X$ is standard.
\item Note that every standard object is internal. Indeed, ${}^*u\in \{{}^*u\}={}^*\{u\}$ for every standard object ${}^*u\in\mathbb V$. Hence ${}^\sigma\mathbb U\subseteq{}^*\mathbb U$.
\item Note that $v\in\mathbb V$ is internal if and only if there is transitive $T\in\mathbb U$ such that $v\in {}^*T$: Indeed, if $v$ is internal, there is $U\in\mathbb U$ such that $v\in {}^*U$. But since $\mathbb U$ is a universe, there is a transitive $T\in\mathbb U$ such that $U\subseteq T$, and so $v\in {}^*U\subseteq{}^*T$.
\item Better yet, given any internal objects $v_1,\dots, v_n$ in $\mathbb V$, there is a transitive $T\in\mathbb U$ such that $v_1,\dots,v_n\in {}^*T$: For there are $A_1,\dots, A_n\in\mathbb U$ such that $v_i\in {}^*A_i$. But then there is transitive $T\in\mathbb U$ such that $A_1\cup\dots\cup A_n\subseteq T$, from which it follows that $v_i\in {}^*A_i\subseteq{}^*T$.
\end{enumerate}\endbox\end{remarks}

 Here is the reason that internal objects play an important role: A transfer map $*:\mathbb U\to\mathbb V$ transfers the truth of bounded $\mathcal L_\in$-formulas from $\mathbb U$ to $\mathbb V$. Thus, for example if $a,b\in\mathbb U$, and  $\varphi\equiv\forall x\in a\;\exists y\in b\;\exists z\in y\;\psi(x,y,z, a,b)$, then with ${}^*\varphi\equiv \forall x\in {}^*a\;\exists y\in {}^*b\;\exists z\in y\;\psi(x,y,z,{}^*a,{}^*b)$ we have $\mathbb U\vDash \varphi$ if and only if $\mathbb V\vDash {}^*\varphi$.  Now to check if $\mathbb V\vDash {}^*\varphi$, we need merely check ${}^*\varphi$ over elements $x\in {}^*a, y\in {}^*b,z\in y$, i.e. we need only consider internal $x,y,z$. Then if ${}^*\varphi$ is true in $\mathbb V$, $\varphi$ is true in $\mathbb U$. 

\begin{theorem} {\rm (Internal/Standard Definition Principle)} Let $\varphi(y,x_1,\dots, x_n)$ be a  bounded $\mathcal L_\in$-formula, and let $B, A_1,\dots, A_n\in\mathbb V$ be internal (resp. standard). Then the set\[\{y\in B: \mathbb V\vDash \varphi[y, A_1,\dots, A_n]\}\] is internal (resp. standard).
\end{theorem}
\bproof The case where $B, A_1,\dots, A_n\in\mathbb V$ are standard follows directly from Proposition \ref{thm_*_comprehension}.

Now suppose that $B, A_1,\dots, A_n\in\mathbb V$ are internal. Let $T\in\mathbb U$ be a transitive set such that $B,A_1,\dots, A_n\in {}^*T$. 
Consider the bounded $\mathcal L_\in$-formula
\[\psi(t, p)\equiv \forall b\in t\;\forall  x_1\in t\;\dots\forall x_n\in t\;\exists u\in p\;\forall y\in t\;(y\in u \leftrightarrow (y\in b\land\varphi(y,x_1,\dots, x_n)).\]
Then $\mathbb U\vDash \psi[T, \mathcal P(T)]$, since if $b, x_1,\dots, x_n\in T$, then the set $\{y\in b:\varphi[y, x_1,\dots, x_n]\}$ is a subset of $b$, thus of $T$, and hence a member of $\mathcal P(T)$. By transfer $\mathbb V\vDash \psi[{}^*T, {}^*\mathcal P(T)]$. In particular, it follows that
\[\mathbb V\vDash \exists u\in {}^*\mathcal P(T)\;\forall y\in {}^*T\;(y\in u\leftrightarrow(y\in B\land\varphi(y,A_1,\dots, A_n)).\]
Thus there is  $u\in{}^*\mathcal P(T)$ such that \[u\cap {}^*T=\{y\in B:\varphi(y,A_1,\dots, A_n)\}\cap {}^*T.\] But since $u\in {}^*\mathcal P(T)$, we have $u\subseteq{}^*T$, so that $u\cap {}^*T=u$. On the other hand, since $B\subseteq {}^*T$ (by transitivity of ${}^*T$), we have that $\{y\in B:\varphi(y,A_1,\dots, A_n)\}\cap {}^*T= \{y\in B:\varphi(y,A_1,\dots, A_n)\}$, and hence $u= \{y\in B:\varphi(y,A_1,\dots, A_n)\}$. Finally, since $u\in {}^*\mathcal P(T)$, $u$ is internal.
\eproof

\begin{lemma} Suppose that $U,V$ are sets, that $U$ is transitive and that $U\subseteq V$. Then $(U,\in)$ is a bounded elementary submodel of $(V,\in)$.
\end{lemma}
\bproof
The proof  that $U\vDash \varphi$ if and only if $V\vDash\varphi$, for any bounded $\mathcal L_\in$-formula $\varphi$ is by induction on formula complexity. The only troublesome case in the induction step for formulas of the form $\varphi(x,y)\equiv\exists z\in x\;\psi(x,y,z)$, where $\psi$ is a bounded $\mathcal L_\in$-formula. 
Suppose that $a,b\in U$ and that $V\vDash\varphi(a,b)$. Then there is $c\in V$ such that  $c\in a$  and $V\vDash\psi(a,b,c)$. As $U$ is transitive and $c\in a$, we also have $c\in U$. By induction hypothesis, therefore, we have $U\vDash \psi(a,b,c)$, from which it follows that $U\vDash \varphi(a,b)$.

\eproof

\begin{proposition} The map $*:(\mathbb U,\in)\hookrightarrow({}^*\mathbb U,\in)$ and the inclusion $({}^*\mathbb U,\in)\hookrightarrow(\mathbb V,\in)$ are bounded elementary embeddings.
\end{proposition}
\bproof  We already know that $({}^*\mathbb U,\in)$ is a transitive submodel of $(\mathbb V,\in)$. Hence  by the previous Lemma, it follows that $({}^*\mathbb U,\in)\hookrightarrow(\mathbb V,\in)$ is a bounded elementary embedding.

Now if $\varphi(x,y)$ is a bounded $\mathcal L_\in$-formula and $a,b\in\mathbb U$, then we have $\mathbb U\vDash\varphi(a,b)$ if and only if $\mathbb V \vDash \varphi({}^*a,{}^*b)$ if and only if ${}^*\mathbb U \vDash \varphi({}^*a,{}^*b)$, which shows that $*:(\mathbb U,\in)\hookrightarrow({}^*\mathbb U,\in)$ is a bounded elementary embedding.
\eproof

\subsection{Examples of Standard and Internal Objects}
In this section, assume that $*:\mathbb U\hookrightarrow \mathbb V$ is a transfer map over $X$.

\begin{lemma} \label{lemma_internal_preservation}
\begin{enumerate}[(a)]
\item If $a_1,\dots, a_n$ are internal, so are $\{a_1,\dots, a_n\}$ and $(a_1,\dots, a_n)$.
\item If $A,B\in\mathbb V$ are internal, so are $A\cup B, A\cap B, A-B, A\times B$.
\item If $\mathcal A$ is internal, then $\bigcup\mathcal A$ and $\bigcap\mathcal A$ are internal.
\item If a binary relation $R$ is internal, then so are $\dom(R),\ran(R)$, and $R^{-1}$. If $C\subseteq\dom(R)$ is internal, then $R[C]\subseteq \ran(R)$ is internal.
\item If $R,S$ are internal binary relations, their composition $R\circ S$ is internal.
\item If a function $f$ is internal, and $a\in \dom(f)$, then $f(a)$ is internal.
\end{enumerate}\end{lemma}
\bproof
(a) For example, there is transitive $T\in\mathbb U$ such that $a_1,\dots,a_n\in {}^*T$.  Then \[\{a_1,\dots, a_n\}=\{x\in {}^*T:x=a_1\lor\dots\lor x=a_n\}.\] This set is internal, by the internal definition principle. Now if $a,b$ are internal, then so are $\{a\}, \{a,b\}$ and hence so is $\{\{a\},\{a,b\}\}=(a,b)$, etc.

(b) For example, if $A,B$ are internal, then there is transitive $T\in\mathbb U$ such that $A, B\subseteq {}^*T$. Since ${}^*T\times {}^*T={}^*(T\times T)$, we see that  \[A\times B=\{x\in {}^*(T\times T):\exists y\in A\;\exists z\in B\;(x=(y,z))\},\] which is internal, by the internal definition principle and the fact that the formula $x=(y,z)$ is bounded, according to Lemma \ref{lemma_define_simple_notions}(c).

(c) Choose transitive $T\in \mathbb U$ such that $\mathcal A\in {}^*T$. Then $\bigcup\mathcal A\subseteq{}^*T$, and hence \[\bigcup\mathcal A=\{x\in{}^*T:\exists y\in \mathcal A\;(x\in y)\},\] which is internal by the internal definition principle. For $\bigcap\mathcal A\subseteq{}^*T$, just replace $\exists y\in\mathcal A\;(x\in y)$ by $\forall y\in\mathcal A\;(x\in y)$.

(d) For example, if $T\in\mathbb U$ is transitive so that $R\in\mathbb {}^*T$, then $\dom(R),\ran(R)\subseteq {}^*T$, and  apply the internal definition principle to \begin{align*}\dom(R)&=\{x\in T:\exists y\in T\;((x,y)\in R)\},\\
\ran(R)&=\{y\in {}^*T:\exists x\in T\;((x,y)\in R)\},\\
R^{-1}&=\{w\in {}^*(T\times T):\exists x\in {}^*T\;\exists y\in {}^*T\;\exists v\in R\;(v=(x,y)\land w=(y,x))\},\\
R[C]&=\{y\in {}^*T: \exists x\in C\;((c,y)\in R)\}.
\end{align*}
Here, $(x,y)\in R$ is short for $\exists z\in R\;(z=(x,y))$, and $z=(x,y)$ is a bounded formula, by Lemma \ref{lemma_define_simple_notions}(c).

(e) Choose $T\in\mathbb U$ transitive so that $R,S\in {}^*T$. Then
\[R\circ S=\{w\in {}^*(T\times T): \exists x\in {}^*T\;\exists y\in {}^*T\;\exists z\in {}^*T\;((x,y)\in S\land (y,z)\in R\land w=(x,z))\}.\]

(f) $f$ is a binary relation, so $\ran(f)$ is internal. Since $f(a)\in\ran(f)$, $f(a)$ is internal, by transitivity of ${}^*\mathbb U$.
\eproof

\begin{lemma} Let $*:\mathbb U\hookrightarrow \mathbb V$ be a transfer map for $X$. Then
\[{}^*V_n(X)=V_n({}^*X)\cap {}^*\mathbb U\quad\text{ i.e.}\quad {}^*V_n(X)=\{x\in V_n({}^*X): x\text{ is internal}\}.\]
\end{lemma}
\bproof Recall that if $\mathbb U$ is a universe over $X$, then $V(X)\subseteq\mathbb U$. Since $\mathbb V$ is a universe over ${}^*X$ (by Lemma \ref{lemma_V_over_*X}),  we have $V({}^*X)\subseteq \mathbb V$. 
Now by Lemma \ref{lemma_define_simple_notions},
\[\mathbb U\vDash \forall x\in V_n(X)\;\varphi_{6,n}(X,x), \quad\text{so}\quad \mathbb V\vDash \forall x\in {}^*V_n(X)\; \varphi_{6,n}({}^*X,x),\] from which it follows that ${}^*V_n(X)\subseteq\{x\in V_n({}^*X): x\text{ is internal}\}$. Now suppose that
${}^*V_n(X)\subsetneqq\{x\in V_n({}^*X): x\text{ is internal}\}$.  Then there must be an internal $a\in V_n({}^*X)$ which is not  in ${}^*V_n(X)$. Since $a$ is internal, there is $B\in \mathbb U$ such that $a\in {}^*B$. Hence
\[\mathbb V\vDash\exists a\in {}^*B[\varphi_{6,n}({}^*X,a)\land a\not\in {}^*V_n(X)],\] from which we obtain
\[\mathbb U\vDash\exists a\in B[\varphi_{6,n}(X,a)\land a\not\in V_n(X)],\] which is impossible.\eproof

 \begin{lemma} \label{lemma_internal_subsets} Let $A\in \mathbb U$ be a set. Then ${}^*\mathcal P(A)=\mathcal P({}^*A)\cap{}^*\mathbb U$ is the set of all internal subsets of ${}^*A$.
\end{lemma}
\bproof Transfer of the true bounded sentence $\forall B\in\mathcal P(A)\;(B\subseteq A)$ shows that every $B\in{}^*\mathcal P(A)$ is a subset of ${}^*A$, and hence in $\mathcal P({}^*A)$. Since every $B\in{}^*\mathcal P(A)$ is obviously internal, we have  ${}^*\mathcal P(A)\subseteq\mathcal P({}^*A)\cap{}^*\mathbb U$.

Now suppose that $B$ is an arbitrary internal subset of ${}^*A$. Then there is $C\in\mathbb U$ such that $B\in {}^*C$. Now transfer of the true sentence $\forall B\in C\;(B\subseteq A\to B\in\mathcal P(A))$ shows that $B\in{}^*\mathcal P(A)$, from which  ${}^*\mathcal P(A)\supseteq\mathcal P({}^*A)\cap{}^*\mathbb U$.
\eproof 

 \begin{lemma}\label{lemma_set of subsets_internal} If $A\in\mathbb V$ is internal, then so is the set $\mathcal P(A)\cap{}^*\mathbb U$ of internal subsets of $A$.
 \end{lemma}
 \bproof
 Choose $T\in\mathbb U$ transitive such that $A\subseteq {}^*T$. If $B\subseteq A$ is internal, then $B$ is an internal subset of ${}^*T$, and hence $B\in {}^*\mathcal P(T)$. So
  \[\mathcal P(A)\cap{}^*\mathbb U=\{B\in {}^*\mathcal P(T): B\subseteq A\}.\]
  Since $A,B, {}^*\mathcal P(T)$ are internal, it follows by the internal definition principle that $\mathcal P(A)\cap{}^*\mathbb U$ is internal.
 \eproof

\begin{lemma}\label{lemma_internal_subsets_internal} Let $\mathcal A\in\mathbb U$ be a family of sets. Then\[{}^*\{\mathcal P(A):A\in \mathcal A\}=\{P_A:A\in{}^*\mathcal A\},\] where $P_A:=\mathcal P(A)\cap {}^*\mathbb U$ is the set of internal subsets of $A$.\end{lemma}
\bproof
Let $\mathcal X:=\{\mathcal P(A):A\in\mathcal A\}$, and observe that $\mathcal X$ is the range of the map $\mathcal A\to\mathcal P\mathcal P(\bigcup\mathcal A):A\mapsto\mathcal P(A)$, so that $\mathcal X\in\mathbb  U$. Let $T\in\mathbb U$ be transitive such that $\mathcal X\subseteq T$. Observe that \[\mathcal X=\{P\in T:\exists A\in\mathcal A\;\forall B\in T\;(B\in P \leftrightarrow B\subseteq A)\}, \] so that  \[{}^*\mathcal X=\{P\in {}^*T:\exists A\in{}^*\mathcal A\;\forall B\in {}^*T\;(B\in P \leftrightarrow B\subseteq A)\}.\] 
Thus $P\in{}^*\mathcal X$ if and only if $P\in {}^*T$ and there is $A\in{}^*\mathcal A$ such that $P\cap {}^*T=\mathcal P(A)\cap {}^*T$. Since $P\in{}^*T$ and ${}^*T$ is transitive, we have $P=\mathcal P(A)\cap {}^*T$ for some $A\in{}^*\mathcal A$. It follows that $P\subseteq P_A$. \\Now if $A\in{}^*\mathcal A$ and $B\in P_A$, then $B\subseteq A$ and there is $S\in \mathbb U$ such that $B\in {}^*S$. Since\[\mathbb U\vDash \forall B\in S\;\forall A\in \mathcal A\;(B\subseteq A\to B\in T),\]
we have that \[\mathbb V\vDash \forall B\in {}^*S\;\forall A\in {}^*\mathcal A\;(B\subseteq A\to B\in {}^*T),\] from which it follows that $B\in {}^*T$, and thus that $B\in P$. Hence $P=\mathcal P(A)\cap {}^*T$ is the set $P_A=\mathcal P(A)\cap {}^*\mathbb U$  of all internal subsets of $A$.  It therefore follows that $P\in{}^*\mathcal X$ if and only if  $P=P_A$ for some $A\in{}^*\mathcal A$.

\eproof

\begin{lemma} Let $\mathcal A\in \mathbb U$ be a family of sets. Then ${}^*(\bigcup\mathcal A)=\bigcup{}^*\mathcal A$.
\end{lemma}
 
 \bproof Just transfer the true bounded formula \[[\forall a\in \bigcup\mathcal A\;\exists A\in\mathcal A\;(a\in A)]\land[ \forall A\in\mathcal A\;\forall a\in A\;(a\in\bigcup\mathcal A)].\]
 \eproof
 (Note that every element of an element of ${}^*\mathcal A$ is internal, since ${}^*\mathbb U$ is transitive. Thus automatically $\bigcup {}^*\mathcal A\subseteq{}^*\mathbb U$, and we do not need to write  ${}^*(\bigcup\mathcal A)=(\bigcup{}^*\mathcal A )\cap {}^*\mathbb U$, as for some of the other operations in this section.)

 We often deal with unions, intersections and products of indexed families of sets.
 Suppose, for example,  that $\mathcal A:=\{A_i:i\in I\}\in\mathbb U$ is an indexed family of sets in $\mathbb U$, and let $f:I\to\mathcal A:i\mapsto A_i$ be the indexing function, where $f\in\mathbb U$. Then ${}^*f:{}^*I\to{}^*\mathcal A$ is a function, with ${}^*f({}^*i)={}^*(f(i))={}^*A_i$ when $i\in I$. For $i\in {}^*I- I$, define ${}^*A_i:={}^*f(i)$ --- {\bf but note that if $i\in {}^*I- I$, then ${}^*A_i$ is not necessarily the $*$-value of a member of $\mathbb U$.} Now since $f$ is surjective, so is ${}^*f$, and hence ${}^*\mathcal A=\{{}^*A_i: i\in {}^*I\}$. Hence by the previous result,  \[{}^*\left(\bigcup_{i\in I}A_i\right)={}^*\bigcup\mathcal A = \bigcup{}^*\mathcal A=\bigcup_{i\in{}^*I} {}^*A_i.\]

 \begin{lemma} Let $A, B\in \mathbb U$ be  sets. Then ${}^*(B^A)={}^*B^{{}^*\!A}\cap{}^*\mathbb U$ is the set of all internal functions ${}^*A\to{}^*B$.
\end{lemma}
\bproof Consider the bounded formula $\varphi_5(f,A,B)$ of Lemma \ref{lemma_define_simple_notions}, which asserts that $f:A\to B$, i.e. that $f\in B^A$. Transfer of the bounded formula
$\forall f\in B^A\;\varphi_5(f,A,B)$ shows that every member of ${}^*(B^A)$ is an internal function ${}^*A\to {}^*B$. Conversely, suppose that $f\in {}^*B^{{}^*A}$ is internal. Let $C\in\mathbb U$ so that $f
\in {}^*C$. Transfer of $\forall f\in C\;(\varphi_5(f,A,B)\to f\in B^A)$ shows that $f\in {}^*(B^A)$.
\eproof

\begin{lemma} Suppose that $A,B\in\mathbb V$ are internal. Then so is the set $B^A\cap {}^*\mathbb U$ of internal functions from $A$ to $B$.
\end{lemma}
\bproof We have seen that $A\times B$ is internal, i.e. there is a transitive $T\in\mathbb U$ such that $A\times B\in {}^*T$. Let $P:=\{X\subseteq A\times B: X\text{ is } internal\}$, which is internal by Lemma \ref{lemma_set of subsets_internal}.  Clearly, \[B^A\cap {}^*\mathbb U = \{f\in P: \varphi_5(f,A,B)\},\] where $\varphi_5(f,A,B)$ is the bounded formula of Lemma \ref{lemma_define_simple_notions} which asserts that $f:A\to B$. By the internal definition principle, $B^A\cap {}^*\mathbb U$ is internal.
\eproof

\begin{lemma}Let $\mathcal A\in\mathbb U$ be a family of sets. Then ${}^*(\prod\mathcal A)=(\prod{}^*\mathcal A)\cap{}^*\mathbb U$.\end{lemma}
\bproof
Recall that $\prod\mathcal A$ is the set of all choice functions, i.e. that  $f\in\prod\mathcal A$ if and only if $f:\mathcal A\to\bigcup\mathcal A$ is such that $f(A)\in A$ for all $A\in \mathcal A$.
Thus transfer of \[\forall f\in \prod\mathcal A\;\Big(\varphi_5(f,\mathcal A,\bigcup\mathcal A)\land\forall A\in\mathcal A\;(f(A)\in A)\Big)\] shows that every member of ${}^*(\prod\mathcal A)$ is an internal choice function  ${}^*\mathcal A\to {}^*(\bigcup\mathcal A)=\bigcup{}^*\mathcal A$, and thus a member of $\prod{}^*\mathcal A\cap {}^*\mathbb U$.

Conversely if $f\in \prod{}^*\mathcal A\cap {}^*\mathbb U$, then there is $C\in\mathbb U$ such that $f\in {}^*C$. Then transfer of
\[\forall f\in C\;\Big(\big[\varphi_5(f,\mathcal A,\bigcup\mathcal A)\land\forall A\in \mathcal A\;(f(A)\in A)\big]\to f\in\prod\mathcal A\Big)\] shows that $f\in{}^*(\prod\mathcal A)$.
\eproof

We can deal with indexed products in a manner very similar to the way we handled indexed unions: Let $\mathcal A=\{A_i:i\in I\}\in \mathbb U$ be an indexed family of sets. Then \[{}^*\Big(\prod_{i\in I}A_i\Big)={}^*(\prod\mathcal A)=\Big(\prod{}^*\mathcal A\Big)\cap{}^*\mathbb U=\Big(\prod_{i\in {}^*I}{}^*A_i\Big)\cap {}^*\mathbb U.\]

\begin{lemma}Let $\mathcal A, \mathcal B\in\mathbb U$ be  families of sets. Then\[{}^*\{A\times B:A\in\mathcal A,B\in\mathcal B\}=\{A\times B:A\in{}^*\mathcal A,B\in{}^*\mathcal B\}.\]\end{lemma}
\bproof Let $\mathcal X:=\{A\times B:A\in\mathcal A,B\in\mathcal B\}$, and choose a transitive $T\in\mathbb U$ such that $\mathcal X\subseteq T$. The bounded formula $\varphi_4$ of Lemma \ref{lemma_define_simple_notions} asserts that $\varphi_4(P,A,B)$ holds if and only if $P=A\times B$. Thus
\[\mathcal X=\{P\in T:\exists A\in\mathcal A\;\exists B\in\mathcal B\;\varphi_4(P,A,B)\}.\] Transfer then yields that ${}^*\mathcal X=\{A\times B:A\in{}^*\mathcal A,B\in{}^*\mathcal B\}\cap {}^*T$.
But if $A\in{}^*\mathcal A, B\in{}^*\mathcal B$, then $A\times B$ is internal, and so there is $S\in\mathbb U$ such that $A\times B\in {}^*S$. But\[\mathbb U\vDash\forall P\in S\;\Big(\exists A\in \mathcal A\;\exists B\in\mathcal B\;\varphi_4(P,A,B)\to P\in T\Big),\] so \[\mathbb V\vDash\forall P\in {}^*S\;\Big(\exists A\in {}^*\mathcal A\;\exists B\in{}^*\mathcal B\;\varphi_4(P,A,B)\to P\in {}^*T\Big).\] Hence $A\times B\in {}^*T$ for all $A\in{}^*\mathcal A, B\in {}^*\mathcal B$. Thus ${}^*\mathcal X=\{A\times B:A\in{}^*\mathcal A,B\in{}^*\mathcal B\}$.\eproof 

\subsection{${}^\sigma A$ is External if $A$ is Infinite}
In this subsection, we assume that $*:\mathbb U\to\mathbb V$ is a nonstandard framework for a set $X$, i.e. a transfer map with the property that  there is a countable set $C\in\mathbb U$ such that ${}^\sigma C\subsetneqq {}^*C$. 

\begin{theorem}\label{thm_nonstandard_framework_nonstandard_elements} If $*:\mathbb U\to\mathbb V$ is a nonstandard framework, then ${}^\sigma A$ is external whenever $ A\in\mathbb U$ is infinite. Hence ${}^\sigma A\subsetneqq {}^*A$ whenever $ A\in\mathbb U$ is infinite. 
\end{theorem} 

\bproof By definition of nonstandard framework, there is  a countable set $C\in\mathbb U$ such that ${}^\sigma C\subsetneqq {}^*C$.  We will begin by showing that the difference $D:={}^*C-{}^\sigma C$ is external. Suppose that $\{c_m:n\in\mathbb N\}$ enumerates $C$. This induces a well-ordering $\preceq$ on $C$ by $c_n\preceq c_m$ if and only if $n\leq m$. Observe that $\preceq\;\in\mathcal P(C\times C)\in\mathbb U$, so that $\preceq\;\in\mathbb U$.
Now the assertion that $\preceq$ is a well-ordering is a bounded sentence, where the fact that every non-empty subset of $C$ has a $\preceq$-least element is given by the bounded sentence.\[\mathbb U\vDash\forall X\in\mathcal P(C)\;[X\neq \varnothing \to \exists x_0\in X\;\forall x\in X\;(x_0\preceq x)].\]
Thus ${}^*\!\!\preceq\;\in\mathbb V$ is a linear ordering on ${}^*C$, and .\[\mathbb V\vDash\forall X\in{}^*\mathcal P(C)\;[X\neq \varnothing \to \exists x_0\in X\;\forall x\in X\;(x_0{}^*\!\!\preceq x)].\]
This does not assert that ${}^*\!\!\preceq$ is a well-ordering on ${}^*C$, however, because it may not be the case that ${}^*\mathcal P(C)=\mathcal P({}^*C)$.  What it does assert is that every  non-empty {\em internal} subset of ${}^*C$ has a ${}^*\!\!\preceq$-least element, since we know that ${}^*\mathcal P(C)=\mathcal P({}^*C)\cap{}^*\mathbb U$ is the set of internal subsets of ${}^*C$.

Now suppose that $D$ is internal, and let $d_0$ be the ${}^*\!\!\preceq$-least element of $D$. We shall obtain a contradiction.

Observe that, for every $n\in\mathbb N$, \[\mathbb U\vDash \forall x\in C\;(x=c_0\lor x=c_1\lor\dots\lor x=c_n\lor x\succ c_n), \] and thus
\[\mathbb V\vDash \forall x\in {}^*C\;(x={}^*c_0\lor x={}^*c_1\lor\dots\lor x={}^*c_n\lor x\:{}^*\!\!\succ {}^*c_n). \] 
Now since $d_0\not\in{}^\sigma C$, we see that $d_0\neq {}^*c_n$ for any $n\in\mathbb N$, and thus that $d_0\;{}^*\!\!\succ {}^*c_n$ for all $n$.
Furthermore, in $\mathbb U$ every element of $C$ has an immediate $\preceq$-predecessor, excepting of course the least element  $c_0$: \[\mathbb U\vDash \forall x\in C\;(x\neq c_0\to\exists y\in C\;(y\prec x\land \forall z\in C\;(z\prec x\to z\preceq y))).\] Hence  \[\mathbb V\vDash \forall x\in {}^*C\;(x\neq {}^*c_0\to\exists y\in {}^*C\;(y \;{}^*\!\!\prec x\land \forall z\in {}^*C\;(z\;{}^*\!\!\prec x\to z\;{}^*\!\!\preceq y))).\]
So if $D$ is internal, its least element $d_0$ has an immediate predecessor $d_{-1}\in {}^*C$. Now clearly we cannot have $d_{-1}= {}^*c_n$ for any $n\in\mathbb N$, for else necessarily $d_0={}^*c_{n+1}\in {}^\sigma C$. So $d_{-1}\in D$ --- contradicting the fact that $d_0$ is the least element of $D$. Thus the assumption that $D$ is internal leads to contradiction, i.e. $D$ is external. 

But as ${}^*C$ is internal, it follows that ${}^\sigma C$ is external, because the difference of two internal sets is internal.

Suppose now that $A\in\mathbb U$ is infinite, and that $f\in\mathbb U$ is a surjection $f:A\twoheadrightarrow C$. It follows easily that ${}^*f:{}^*A\twoheadrightarrow {}^*C$. Since ${}^*(f(a))={}^*f({}^*a)$ for all $a\in A$, we have ${}^*f[{}^\sigma A]= {}^\sigma C$. If ${}^\sigma A$ is internal, then by Lemma \ref{lemma_internal_preservation} it would follows that ${}^\sigma C$ is internal as well. As this is false, ${}^\sigma A$ is external.

Finally, it is always the case that ${}^\sigma A\subseteq {}^*A$. If ${}^\sigma A= {}^*A$, then ${}^\sigma A$ would be internal. If $A$ is infinite, this is not the case. Hence ${}^\sigma A\subsetneqq {}^*A$ if $A\in\mathbb U$ is  infinite.
\eproof

\section{Hyperfinite Sets}
\fancyhead[LO]{Hyperfinite Sets}
\subsection{The Set ${}^*\mathbb N$ of Hypernatural Numbers}
Throughout this section, assume that $*:\mathbb U\to\mathbb V$ is a nonstandard embedding for an infinite set $X$. Without loss of generality (e.g. by renaming elements) we may assume that $\mathbb N\subseteq X$. Then as $\mathcal P(X)\in\mathbb U$ and $\mathbb U$ is transitive, also $\mathbb N\in\mathbb U$. Since ${}^*x=x$ for all $x\in X$, it follows in particular that ${}^\sigma \mathbb N=\mathbb N$.  Since $*$ is nonstandard, $\mathbb N\subsetneqq{}^*\mathbb N$.

Given a fixed but arbitrary $N\in\mathbb N$, transfer of the true bounded sentence \[\forall n\in \mathbb N\;(n=1\lor n=2\lor \dots\lor n=N\lor n>N)\] shows that any member of ${}^*\mathbb N- \mathbb N$ is $>N$. As this is true for all $N\in\mathbb N$, it follows that every member of ${}^*\mathbb N-\mathbb N$ is greater than any natural number, and thus said to be {\em infinite}. We thus define $\mathbb N_\infty:={}^*\mathbb N-\mathbb N$ to be the set of infinite natural numbers, so that we have the disjoint union ${}^*\mathbb N=\mathbb N\cup \mathbb N_\infty$.  As $\mathbb N={}^\sigma\mathbb N$ is external, so is $\mathbb N_\infty$.

Now  the structure $(\mathbb N, +,\dotdiv,\cdot, \leq, 0,1)\in\mathbb U$ is elementarily equivalent to the corresponding structure  $({}^*\mathbb N, +,\dotdiv,\cdot,\leq, 0,1)\in{}^*\mathbb U$ (where  $n\dotdiv m:=n-m$ if $n\geq m$ and $:=0$ else). We should really have written $({}^*\mathbb N, {}^*+,{}^*\dotdiv,{}^*\cdot, {}^*\leq,0,1)$,  but we drop the stars on the arithmetic operations and order relation for easier reading. It is clear that if  $n\in\mathbb N_\infty$ then $n-1\in\mathbb N_\infty$ also, and thus
$\mathbb N_\infty$ has no least element, i.e. ${}^*\mathbb N$ is not well-ordered. However:

\begin{theorem} Every non-empty {\em internal} subset of ${}^*\mathbb N$ has a least element.
\end{theorem}
\bproof  Let $A\subseteq{}^*\mathbb N$ be internal and non-empty. Then $A\in\mathcal P({}^*\mathbb N)\cap{}^*\mathbb U={}^*\mathcal P(\mathbb N)$. Transfer of the true bounded sentence
\[\forall A\in\mathcal P(\mathbb N)\;[A\neq \varnothing\to\exists a_0\in A\;\forall a\in A\;(a_0\leq a)]\] now shows that $A$ has a least element.\eproof

\begin{theorem}{\rm (Overflow and Underflow)}
\begin{enumerate}[(a)]\item {\rm (Overflow)}  Let $N\in\mathbb N$, and suppose that $X\in{}^*\mathcal P(\mathbb N)$ is an {\em internal} subset of ${}^*\mathbb N$ with the property that whenever $n\in\mathbb N$ satisfies $n\geq N$, then $n\in X$. Then there is $M\in \mathbb N_\infty$ such that whenever $n\in {}^*\mathbb N$ satisfies $N\leq n\leq M$, then $n\in X$ .
\item {\rm (Underflow)} Let $M\in\mathbb N_\infty$, and suppose that $X\in{}^*\mathcal P(\mathbb N)$ is an {\em internal} subset of ${}^*\mathbb N$ with the property that whenever $n\in{}^*\mathbb N_\infty$ satisfies $n\leq M$, then $n\in X$. Then there is $N\in \mathbb N$ such that whenever $n\in {}^*\mathbb N$ satisfies $N\leq n\leq M$, then $n\in X$ .
\end{enumerate}
\end{theorem}
\bproof
(a) By the internal definition principle, the set $Y:=\{n\in {}^*\mathbb N: n>N\land n\in {}^*\mathbb N-X\}$ is internal.  If $Y=\varnothing$, then every $M\in\mathbb N_\infty$ has the desired property. If $Y$ is non-empty, then it must have a least element $K$, i.e. $K$ is least such that $K>N$ and $K\not\in X$. By the assumption on $X$ we must have $K\in \mathbb N_\infty$.  It follows that every infinite natural number which is $<K$ belongs to $X$, so let $M:=K-1$.

(b) By the internal definition principle, the set $Z:=\{k\in {}^*\mathbb N: \forall n\in {}^*\mathbb N\;(k\leq n\leq M\to n\in X)\}$ is internal.  Then $Z\neq \varnothing$, since $M\in\mathbb Z$, and so $Z$ has a least member $N$. As $N\leq M$,  the assumption on $X$ implies that we must have $N\in\mathbb N$. 
\eproof

\subsection{Hyperfinite sets}
In this section, we suppose that $*:\mathbb U\to\mathbb V$ is a transfer map over  a set $X$, where $\mathbb N\subseteq X$. Recall that the $<$-relation on $\mathbb N$ is just a subset of $<\;\subseteq \mathbb N\times\mathbb N$, so to say that $m<n$ is equivalent to saying $(m,n)\in\; <$.
\begin{definition}\rm  A set $A\in\mathbb V$ is said to be hyperfinite (or $*$-finite) if and only if there is an {\em  internal} bijection $f:\{0,1,\dots, n-1\}\to A$ for some $n\in {}^*\mathbb N$.\\ In that case, we denote $\#A=n$.\endbox\end{definition}

Observe that  if $n\in{}^*\mathbb N$, then $\{0,1,\dots, n-1\}=\{x\in {}^*\mathbb N:x<n\}$ is internal, by the internal definition principle. Hence if $f:\{0,1,\dots, n-1\}\to A$ is an internal bijection, then $A=\ran(f)$ is internal, i.e. every hyperfinite set is internal.

Observe also that there is a bounded formula $\psi(A,f,n,U)$ which asserts that $f:U\to A$ is a bijection, and that $U=\{0,\dots, n-1\}$: Indeed, $\psi$ is the conjunction of 
\begin{itemize}\item $U\subseteq\mathbb N\land\forall m\in\mathbb N\;(m<n \leftrightarrow m\in U)$  --- i.e. $U=\{0,1,\dots n-1\}$.
\item $\forall x\in f\;\exists u\in U\;\exists a\in A\; (x=(u,a))$ --- i.e. $f$ is a binary relation, with  $\dom(f)\subseteq U$ and $\ran(f)\subseteq A$. \\(Recall that the bounded formula $\varphi_{2,2}(c,a,b)$ of Lemma \ref{lemma_define_simple_notions} asserts that $c=(a,b)$.)
\item $\forall u\in U\;\exists a\in A\;((u,a)\in f)$ --- i.e. $\dom(f)\supseteq U$. 
\item $\forall a\in A\;\exists u\in U\;((u,a)\in f)$ --- i.e. $A\subseteq\ran(f)$.
\item $\forall u\in U\;\forall a\in A\;\forall b\in A\;(u,a)\in f\land (u,b)\in f\to a=b)$ --- i.e $f$ is a function.
\item $\forall u\in U\;\forall v\in U\;(\exists a\in A\;((u,a)\in f\land (v,a)\in f)\to u=v)$ --- i.e. $f$ is 1-1.
\end{itemize}
Now define the formula \qquad$\Psi(A,f,n)\equiv \exists U\in\mathcal P(\mathbb N)\;\psi(A,f,n,U)$. \\ Then $\Psi(A,f,n)$ is a bounded formula which asserts that $f$ is a bijection $f:\{0,\dots, n-1\}\to A$. Transfer guarantees that ${}^*\Psi(A,f,n)$ asserts the same for $A, f\in {}^*\mathbb U$ and $n\in {}^*\mathbb N$. We thus see that, by definition,  $A\in \mathbb V$ is hyperfinite if and only if there are an internal $f$ and an $n\in {}^*\mathbb N$ such that ${}^*\Psi(A,f,n)$.
\vskip0.3cm
Suppose that $B\in\mathbb U$, and recall that ${}^*\mathcal P(B)=\mathcal P({}^*B)\cap {}^* \mathbb U$ is the set of all internal subsets of $B$. Let  \[\mathcal P^{<\omega}(B):=\{C\subseteq B: C\text{ is finite}\}.\] 
Then ${}^*\mathcal P^{<\omega}(B)\subseteq {}^*\mathcal P(B)$. The next theorem characterizes the hyperfinite sets as members of some ${}^*\mathcal P^{<\omega}(B)$.

\begin{theorem} A set $A\in\mathbb V$ is hyperfinite if and only if there is $B\in\mathbb U$ such that $A\in{}^*\mathcal P^{<\omega}(B)$.
\end{theorem}

\bproof $(\Rightarrow)$: 
Suppose that $A\in\mathbb V$ is hyperfinite, and let $f:\{0,1,\dots, n-1\}\to A$ be an internal bijection, for some $n\in{}^*\mathbb N$. Let $U:=\{0,1,\dots, n-1\}=\dom(f)$ and  $A:=\ran(f)$. Then $U,A$ are internal so $U \in{}^*\mathcal P(\mathbb N)$,  and there is a transitive set $B\in\mathbb U$ such that $A\in {}^*B$. As ${}^*B$ is transitive also, we  have $A\subseteq {}^*B$, and since $f$ is internal, we must have  $f\in\mathcal P({}^*\mathbb N\times{}^*B)\cap{}^*\mathbb U={}^*\mathcal P(\mathbb N\times B)$. It follows that ${}^*\Psi(A,f,n)$ holds. Now observe that\[\mathbb U\vDash \forall X\in B\;\Big(\exists n\in\mathbb N\;\exists f\in \mathcal P(\mathbb N\times B)\;\Psi(X,f,n)\to X\in\mathcal P^{<\omega}(B)\Big)\]
Applying transfer, with $X=A$,   it follows that $A\in{}^*\mathcal P^{<\omega}(B)$.

$(\Leftarrow)$: Conversely, suppose that $A\in{}^*\mathcal P^{<\omega}(B)$ for some $B\in\mathbb U$. Observe that\[\mathbb U\vDash \forall X\in\mathcal P^{<\omega}(B)\;\exists n\in\mathbb N\;\exists f\in \mathcal P(\mathbb N\times B)\;\Psi(X,f,n),\] and hence by transfer that there exist $n\in {}^*\mathbb N$ and $f\in {}^*\mathcal P(\mathbb N\times B)$ such that ${}^*\Psi(A, f,n)$. Then $f$ is an internal bijection from $\{0,1,\dots, n-1\}$ onto $A$.
\eproof

\section{Enlargements and Saturation}
\fancyhead[LO]{Enlargements and Saturation}
\subsection{Definitions and Basic Properties}

Throughout this section, assume that $*:\mathbb U\to\mathbb V$ is a transfer map for a set $X$.

Recall Convention \ref{convention_expanded_language}.
Let $\mathcal L_{\mathbb U}$ denote the expansion of the langue $\mathcal L_\in$ (or some expansion thereof) with additional constant symbols for element of $\mathbb U$. Similarly, let $\mathcal L_{{}^*\mathbb U}$ denote the expansion with constant symbols for every internal set.

Recall also that a family $\mathcal A$ of sets has the {\em finite intersection property} (f.i.p.) if and only if the intersection of any finitely many members of $\mathcal A$ is non-empty.

\begin{definition}\label{defn_enlargement_saturation}\rm 
Suppose that $\kappa$ is an infinite cardinal.
\begin{enumerate}[(i)]\item $*:\mathbb U\to\mathbb V$ is a $\kappa$-{\em enlargement} if and only if for every set $\Sigma(x)$ of $<\kappa$-many bounded formulas of $\mathcal L_{\mathbb U}$, if $\Sigma(x)$ is finitely satisfiable in $\mathbb U$ by elements of some set  $T\in\mathbb U$, then $\Sigma(x)$ is satisfiable in $\mathbb V$ by an element of $^*\mathbb U$.\\
$*$ is an {\em enlargement} if it is a $|\mathbb U|^+$-enlargement, or, what is equivalent, if it is a  $\kappa$-enlargement for any cardinal $\kappa$.
\item $*:\mathbb U\to\mathbb V$ is $\kappa$-{\em saturated}  if and only if for every set $\Sigma(x)$ of $<\kappa$-many bounded formulas of $\mathcal L_{{}^*\mathbb U}$, if $\Sigma(x)$ is finitely satisfiable in $\mathbb V$ by elements of some set  $T\in{}^*\mathbb U$, then $\Sigma(x)$ is satisfiable in $\mathbb V$ by an element of $^*\mathbb U$.\\
$*$ is {\em polysaturated} if it is $|\mathbb U|^+$-saturated.
\end{enumerate}\endbox
\end{definition}

\begin{theorem}\label{thm_enlargement_saturation_intersection}
\begin{enumerate}[(a)]\item $*:\mathbb U\to\mathbb V$ is a $\kappa$-{\em enlargement} if and only if whenever $\mathcal A\subseteq\mathbb U$ is a family of sets of cardinality $<\kappa$ such that $\mathcal A$ has the f.i.p., then $\bigcap {}^\sigma\mathcal A=\bigcap\{{}^*A:A\in\mathcal A\}$ is non-empty.
\item $*:\mathbb U\to\mathbb V$ is $\kappa$-{saturated} if and only if whenever $\mathcal A\subseteq {}^*\mathbb U$ is a family of internal sets of cardinality $<\kappa$ such that $\mathcal A$ has the f.i.p., then $\bigcap\mathcal A$ is non-empty.
\end{enumerate}
\end{theorem}

\bproof 
(a) Suppose first that $*:\mathbb U\to\mathbb V$ is a $\kappa$-enlargement, and that $\mathcal A:=\{A_\beta:\beta<\alpha\}\subseteq\mathbb U$ is a family of sets of cardinality $<\kappa$ such that $\mathcal A$ has the f.i.p. By replacing each $A_\beta$ by $A_\beta\cap A_0$ we may without loss of generality assume that each $A_\beta\subseteq A_0$ --- This does not affect the f.i.p. nor the intersection of all the $A_\beta$.  Let $\Sigma(x):=\{\sigma_\beta(x):\beta<\alpha\}$, where $\sigma_\beta(x)\equiv x\in A_\beta$. Then $\Sigma(x)$ is a set of $<\kappa$-many bounded formulas of $\mathcal L_{\mathbb U}$ which is finitely satisfiable by elements of $A_0$. Hence $\Sigma(x)$ is satisfiable in $\mathbb V$ by an element of $^*\mathbb U$. Thus there is $a\in {}^*\mathbb U$ such that $a\in {}^*A_\beta$ for all $\beta<\alpha$, so that $\bigcap{}^\sigma\mathcal A\neq \varnothing$.

For the reverse direction,  suppose we have a set $\Sigma(x)=\{\sigma_\beta(x):\beta<\alpha\}$ of $<\kappa$-many bounded formulas of $\mathcal L_{\mathbb U}$, and that $\Sigma(x)$ is finitely satisfiable in $\mathbb U$ by elements of some set  $T\in\mathbb U$. Let $A_\beta:=\{t\in T: \mathbb U\vDash \sigma_\beta(t)\}$. As $\Sigma(x)$ is finitely satisfiable, the family $\mathcal A:=\{A_\beta:\beta<\alpha\}$ has the f.i.p. If $s\in\bigcap{}^\sigma\mathcal A$, then $s\in {}^*\mathbb U$. As $\mathbb U\vDash \forall x\in A_\beta(\sigma_\beta(x))$, we have $\mathbb V\vDash \forall x\in {}^*A_\beta(\sigma_\beta(x))$ so that $\mathbb V\vDash \Sigma(s)$.

(b) Suppose first that $*:\mathbb U\to\mathbb V$ is $\kappa$-saturated, and that $\mathcal A:=\{A_\beta:\beta<\alpha\}\subseteq{}^*\mathbb U$ is a family of internal sets of cardinality $<\kappa$ such that $\mathcal A$ has the f.i.p. By replacing each $A_\beta$ by $A_\beta\cap A_0$ we may without loss of generality assume that each $A_\beta\subseteq A_0$ --- This does not affect the f.i.p. nor the intersection of all the $A_\beta$.  Let $\Sigma(x):=\{\sigma_\beta(x):\beta<\alpha\}$, where $\sigma_\beta(x)\equiv x\in A_\beta$. Then $\Sigma(x)$ is a set of $<\kappa$-many bounded formulas of $\mathcal L_{{}^*\mathbb U}$ which is finitely satisfiable by elements of $A_0$. Hence $\Sigma(x)$ is satisfiable in $\mathbb V$ by an element of $^*\mathbb U$. Thus there is $a\in {}^*\mathbb U$ such that $a\in A_\beta$ for all $\beta<\alpha$, so that $\bigcap\mathcal A\neq \varnothing$.

For the reverse direction,  suppose we have a set $\Sigma(x)=\{\sigma_\beta(x):\beta<\alpha\}$ of $<\kappa$-many bounded formulas of $\mathcal L_{{}^*\mathbb U}$, and that $\Sigma(x)$ is finitely satisfiable in $\mathbb V$ by elements of some set  $T\in{}^*\mathbb U$. Let $A_\beta:=\{t\in T:  \sigma_\beta(t)\}$. By the internal definition principle, each $A_\beta$ is internal. As $\Sigma(x)$ is finitely satisfiable by members of $T$, the family $\mathcal A:=\{A_\beta:\beta<\alpha\}$ has the f.i.p. If $s\in\bigcap\mathcal A$, then $s\in {}^*\mathbb U$, and  $\mathbb V\vDash \Sigma(s)$.
\eproof 

Note that we do  not demand that $\mathcal A\in\mathbb U$ in (a), or that $\mathcal A\in\mathbb V$ in (b) of Theorem \ref{thm_enlargement_saturation_intersection}. However:

\begin{theorem} \label{thm_fip_element} To verify that a transfer map $*:\mathbb U\to\mathbb V$ for a set $X$ is a $\kappa$-enlargement, it suffices to consider sets $\mathcal A\in\mathbb U$ in (a)  of of Theorem \ref{thm_enlargement_saturation_intersection}.
\\To verify that a transfer map $*:\mathbb U\to\mathbb V$ for a set $X$ is $\kappa$-saturated, it suffices to consider sets  $\mathcal A$ which are subsets of standard sets --- and thus in $\mathbb V$ --- in (b) of of Theorem \ref{thm_enlargement_saturation_intersection}.
\end{theorem}

\bproof Suppose that $*$ satisfies (i) of Definition \ref{defn_enlargement_saturation} for sets $\mathcal A\in\mathbb U$. Let $\mathcal A'\subseteq\mathbb U$ be a set of cardinality $<\kappa$ which satisfies the f.i.p. Choose $A_0\in\mathcal A'$, and let $\mathcal A:=\{A\cap A_0:A\in\mathcal A'\}$. Observe that $|\mathcal A|<\kappa$, that $\mathcal A$ satisfies the f.i.p. as well, and that since $\mathcal A\subseteq\mathcal P(A_0)$ we have $\mathcal A\in\mathbb U$. By assumption, $\bigcap{}^\sigma\mathcal A\neq \varnothing$. Now as ${}^*(A\cap A_0)\subseteq {}^*A$ for all $A\in\mathcal A'$, we have $\varnothing\neq \bigcap {}^\sigma\mathcal A=\bigcap{}^\sigma\mathcal A'$.

Next, suppose that $*$ satisfies (ii) of Definition \ref{defn_enlargement_saturation} for $\mathcal B\subseteq {}^*\mathcal A$, where $\mathcal A\in\mathbb U$.  Let $\mathcal B'\subseteq{}^*\mathbb U$ be a family of internal sets of cardinality $<\kappa$ which satisfies the f.i.p. Choose $B_0\in\mathcal B'$, and let $\mathcal B:=\{B\cap B_0:B\in\mathcal B'\}$. Observe that $|\mathcal B|<\kappa$, that $\mathcal B$ has the f.i.p. also, and that  $\mathcal B\subseteq {}^*\mathbb U\cap \mathcal P(B_0)$. By Lemma \ref{lemma_set of subsets_internal},  $P:={}^*\mathbb U\cap \mathcal P(B_0)$ is internal, and hence there is some transitive $\mathcal A\in\mathbb U$ such that $P\in{}^*\mathcal A$. Then also $\mathcal B\subseteq {}^*\mathcal A$, and hence by assumption we have that $\bigcap\mathcal B\neq 0$. Hence $\varnothing\neq \bigcap\mathcal B=\bigcap\mathcal B'$.
\eproof

\begin{remark}\rm It follows from Theorem \ref{thm_enlargement_equiv} that $*:\mathbb U\to\mathbb V$ is an enlargement if and only if it is a $|\mathbb U|$-enlargement, as $\mathcal A\in\mathbb U$ implies $|\mathcal A|<|\mathbb U|$.\endbox
\end{remark}

Suppose that $B\in{}^*\mathbb U$ is an infinite internal set. Then for $b\in B$, the set $B_b:=\{c\in B:c\neq b\}$ $= B-\{b\}$ is internal, by the internal definition principle and the transitivity of ${}^*\mathbb U$. Clearly $\mathcal B:=\{B_b:b\in B\}$ has the f.i.p., yet $\bigcap\mathcal B=\varnothing$. Hence $*$ cannot be a $|B|^+$-saturated extension. In particular, if $*:\mathbb U\to\mathbb V$ is a  $\kappa$-saturated transfer map for an infinite set $X$ of atoms, then necessarily $|{}^*X|\geq\kappa$.
It is therefore impossible to find an extension which is $\kappa$-saturated for every cardinal $\kappa$.

The above argument also shows that:
\begin{proposition} If  $*:\mathbb U\to\mathbb V$ is $\kappa$-saturated and $B\in\mathbb V$ has cardinality $|B|<\kappa$, then $B$ is external.
\endbox
\end{proposition}

\begin{lemma}\label{lemma_polysaturated_implies_enlargement} If $*:\mathbb U\to\mathbb V$ is $\kappa$-saturated, it is a $\kappa$-enlargement. Hence a polysaturated extension is an enlargement.
\end{lemma}
\bproof If $\mathcal A\subseteq\mathbb U$ is a family of sets with the f.i.p., then $\mathcal B:={}^\sigma\mathcal A$ is a family of internal  (indeed, standard) sets with the f.i.p. For, given $A_1,\dots, A_n\in\mathcal A$, transfer of the  bounded sentence $\exists x\in A_1(x\in A_2\land\dots\land x\in A_n)$ shows that ${}^*A_1\cap\dots\cap{}^*A_n\neq \varnothing$. Clearly $|\mathcal B|=|\mathcal A|$.\eproof

\subsection{Enlargements, Concurrency and Hyperfinite Approximation}

\begin{theorem} A  transfer map $*:\mathbb U\to\mathbb V$ for an infinite set $X$ is $\omega_1$-enlargement if and only if it is a nonstandard embedding.

\end{theorem}
\bproof $(\Rightarrow)$: Suppose that $*$ is an $\omega_1$-enlargement. To show that $*$ is a nonstandard embedding, it suffices to show that there is a countable $C\in\mathbb U$ such that ${}^\sigma C\subsetneqq {}^*C$. So let $C\subseteq X$ be countable. Since $\mathcal P(X)\in\mathbb U$ and $\mathbb U$ is transitive, we have $C\in\mathbb U$. Let $\mathcal C:=\{C-\{c\}:c\in C\}$. Then $\mathcal C$ has the f.i.p. and $|\mathcal C|<\omega_1$, so $\bigcap {}^\sigma\mathcal C\neq \varnothing$. Pick $c_0\in \bigcap {}^\sigma\mathcal C$. Then $c_0\in {}^*(C-\{c\})={}^*C-\{{}^*c\}$ for every $c\in C$. Hence $c_0\in {}^*C-{}^\sigma C$.

$(\Leftarrow)$: Now assume that $*$ is a nonstandard embedding. Let $\mathcal A\subseteq \mathbb U$ be a family of sets with the f.i.p. such that $|\mathcal A|<\omega_1$, i.e. $|\mathcal A|$ is countable. Let $\mathcal A=\{A_1,A_2,\dots\}$ be an enumeration of $\mathcal A$. Without loss of generality (by renaming elements if necessary), we may assume that $\mathbb N\subseteq X$. Let $f:\mathbb N\twoheadrightarrow\mathcal P(A_1):n\mapsto\bigcap_{k\leq n}A_k$. Then $f\in\mathbb U$, and $\varnothing\not\in \text{ran} f$. By transfer, $\varnothing\not\in \text{ran}\,{}^*f$, since ${}^*\varnothing =\varnothing$ and ${}^*(\text{ran} f)= \text{ran}\,{}^*f$. In particular, for all $n\in {}^*\mathbb N$, ${}^*f(n)\neq \varnothing$.  
Now fix an arbitrary $n\in\mathbb N$, and note that \[\mathbb U\vDash\forall m\in \mathbb N\;(m>n\to f(m)\subseteq A_n).\] By transfer
\[\mathbb V\vDash \forall m\in {}^*\mathbb N\;(m>n\to {}^*f(m)\subseteq {}^*A_n).\] (Recall ${}^*x=x$ for all $x\in X$, and thus ${}^*n=n$.)
Since $*$ is a nonstandard embedding, it follows from Theorem \ref{thm_nonstandard_framework_nonstandard_elements} that there exists $m_0\in {}^*\mathbb N-{}^\sigma\mathbb N$, i.e. an infinite $m_0$. Then we have both ${}^*f(m_0)\neq \varnothing$ and $m_0\geq n$. It follows that ${}^*f(m_0)\subseteq {}^*A_n$. Since $n\in\mathbb N$ was arbitrary, we have that $\varnothing\neq  f(m_0)\subseteq\bigcap_{n\in\mathbb N} {}^*A_n=\bigcap{}^\sigma \mathcal A$. This shows that $*$ is an $\omega_1$-enlargement.
\eproof

As a corollary we immediately see that:
\begin{corollary}
Every enlargement, and hence every polysaturated extension, is a nonstandard embedding.\endbox
\end{corollary}

\begin{definition}\rm A binary relation $R$ is {\em concurrent} if and only if for every finite $\{x_1,\dots, x_n\}\subseteq\dom(R)$ there is $y\in \ran(R)$ such that $x_iRy$ for all $i=1,\dots, n$.\endbox
\end{definition}

For a set $A$, an  important example of a concurrent relation is the set \[R:=\{(a,F): a\in F\in\mathcal P^{<\omega}(A)\}\subseteq A\times\mathcal P^{<\omega}(A).\] Indeed, given $a_1,\dots, a_n\in A=\dom(R)$, we can define $F:=\{a_1,\dots, a_n\}$, and then observe that $a_iRF$ for all $i=1,\dots, n$.

\begin{theorem}\label{thm_enlargement_equiv} Suppose that $*:\mathbb U\to \mathbb V$ is a transfer map for a set $X$. Then the following are equivalent: 
\begin{enumerate}[(i)]\item $*:\mathbb U\to \mathbb V$ is a $\kappa$-enlargement.
\item For every concurrent binary relation $R\in\mathbb U$ with $|\dom(R)|<\kappa$ there is $y_0\in
\mathbb V$ such that ${}^*x({}^*R)y_0$ for all $x\in \dom(R)$.
\item (Hyperfinite Approximation) For each set $A\in \mathbb  U$ with $|A|<\kappa$ there is a hyperfinite subset $B$ of ${}^*A$ which contains all the standard members of ${}^*A$, i.e.\[{}^\sigma A\subseteq B\in {}^*\mathcal P^{<\omega}(A).\]
\end{enumerate}
\end{theorem}

\bproof
(i) $\Rightarrow$ (ii): Suppose $|\dom(R)|<\kappa$. For $x\in \dom(R)$, let $A_x:=\{y\in \ran(R): xRy\}$, and let $\mathcal A:=\{A_x:x\in\dom(R)\}$. Then $\mathcal A\subseteq\mathcal P(\ran(R))$, so $\mathcal A\in\mathbb  U$. Since $R$ is concurrent, $\mathcal A$ has the f.i.p. In addition $|\mathcal A|<\kappa$. Since $*$ is a $\kappa$-enlargement, there is $y_0\in\bigcap {}^\sigma\mathcal A$. Then $y_0\in {}^*A_x$ for all $x\in\dom(R)$. Now $\mathbb U\vDash \forall y\in A_x\;(xRy)$, so by transfer we have $\mathbb V\vDash \forall y\in{}^*A_x\;({}^*x({}^*R)y)$. It follows that ${}^*x({}^*R)y_0$ for all $x\in \dom(R)$.\\
(ii)$\Rightarrow$ (iii): Given a set $A\in\mathbb U$ with $|A|<\kappa$, define the concurrent relation $R$ by  \[R:=\{(a,F): a\in F\in\mathcal P^{<\omega}(A)\}.\] Then $|\dom(R)|<\kappa$, so by assumption, there is $B\in\mathbb V$ such that  ${}^*a({}^*R)B$ --- i.e. such that ${}^*a\in B$  --- for all $a\in A$. Thus ${}^\sigma A\subseteq B$. Now as $B\in \ran({}^*R)$, and since $\ran(R)\subseteq\mathcal P^{<\omega}(A)$, we have that $\ran({}^*R)\subseteq {}^*\mathcal P^{<\omega}(A)$, so that $B\in {}^*\mathcal P^{<\omega}(A)$.\\
(iii) $\Rightarrow$ (i):  We use Theorem \ref{thm_fip_element}. Let $\mathcal A\in \mathbb U$ have the f.i.p. with $|\mathcal A|<\kappa$, and let $T\in\mathbb U$ be transitive so that $\mathcal A\subseteq T$. Then $A\subseteq T$ for all $A\in\mathcal A$, and hence \[\mathbb U\vDash\forall \mathcal F\in\mathcal P^{<\omega}(\mathcal A)\;\exists x\in T\;\forall A\in\mathcal F\;(x\in A).\] By assumption, there is a hyperfinite $\mathcal B\subseteq {}^*\mathcal A$ such that ${}^\sigma\mathcal A\subseteq \mathcal B$. Then by transfer, with $\mathcal F=\mathcal B\in {}^*\mathcal P^{<\omega}(\mathcal A)$, there is $x\in {}^*T$ such that $x\in A$ for all $A\in\mathcal B$. In particular, $x\in {}^*A$ for all $A\in\mathcal A$.
\eproof

\subsection{Saturation and Concurrency}\label{section_saturation}

\begin{theorem}\label{thm_saturation_concurrency} Suppose that $*:\mathbb U\to \mathbb V$ is a nonstandard embedding.
Then the following are equivalent:
\begin{enumerate}[(i)]\item $*:\mathbb U\to \mathbb V$ is $\kappa$-saturated.
\item For every {\em internal} concurrent binary relation $R\in\mathbb V$ with and every (internal or external) $A\subseteq \dom(R)$ with  $|A|<\kappa$,  there is $y_0\in
\mathbb V$ such that $xRy_0$ for all $x\in A$.
\end{enumerate}
\end{theorem}

\bproof
(i) $\Rightarrow$ (ii): Suppose that $R\in\mathbb V$ is an internal concurrent relation, and that $A\subseteq\dom(R)$ is such that $|A|<\kappa$. For each $x\in A$, let $A_x:=\{y\in \ran(R): xRy\}$. By the internal definition principle, each $A_x$ is internal. Put $\mathcal A:=\{A_x:x\in A\}$, so that $\mathcal A$ is a family of internal sets with $|\mathcal A|<\kappa$. Since $R$ is concurrent, $\mathcal A$ has the f.i.p. By $\kappa$-saturation, there is $y_0\in\bigcap \mathcal A$. Then $y_0\in A_x$ for all $x\in A$, i.e. $xRy_0$ holds for all $x\in A$.

(ii) $\Rightarrow$ (i): Suppose that $\mathcal A$ is a family of internal sets with the f.i.p., where $|\mathcal A|<\kappa$. We must show that $\bigcap\mathcal A\neq\varnothing$. Fix $A_0\in \mathcal A$. By replacing each $A\in\mathcal A$ by $A\cap A_0$, we may assume that $A\subseteq A_0$ for each $A\in\mathcal A$. Thus $\mathcal A\subseteq \mathcal P(A_0)\cap{}^*\mathbb U$.  By Lemma \ref{lemma_set of subsets_internal}, $\mathcal P(A_0)\cap{}^*\mathbb U$ is an internal set. Consider the relation
\[R:=\{(A,a)\in(\mathcal P(A_0)\cap{}^*\mathbb U)\times A_0: a\in A\}.\]
By the internal definition principle, $R$ is internal. Now $\mathcal A\subseteq \dom(R)$, as each $A\in\mathcal A$ is non-empty (by the f.i.p.). Hence there is $a_0\in\mathbb V$ such that $(A,a_0)\in R$ for all $A\in\mathcal A$. But then $a_0\in\bigcap\mathcal A$.
\eproof

\subsection{Comprehensiveness}\label{section_comprehensiveness}

\begin{definition}\label{defn_comprehensiveness}\rm Let $\kappa$ be an infinite cardinal. A transfer map $*:\mathbb U \to\mathbb V$ is said to be {\em $\kappa$-comprehensive} if and only if for any sets $A,B\in\mathbb U$ such that $|A|<\kappa$, and any map $f:A\to {}^*B$, there is an {\em internal} function ${}^+\!f:{}^*A\to {}^*B$ with the property that ${}^+\!f({}^*a)=f(a)$ for all $a\in A$.\\
$*:\mathbb U \to\mathbb V$ is said to be {\em comprehensive} if and only if it is $\kappa$-comprehensive for every cardinal $\kappa$, or equivalently, if it is $|\mathbb U|^+$--comprehensive.\\
$*:\mathbb U \to\mathbb V$ is said to be {\em countably  comprehensive} if and only if it is $\omega_1$-comprehensive. \endbox\end{definition}

Here is a small improvement:
\begin{proposition}\label{propn_comprehensiveness} If a transfer map $*:\mathbb U \to\mathbb V$ is $\kappa$-comprehensive, then for any set $A\in\mathbb U$ of cardinality $|A|<\kappa$, and any {\em internal} set $B\in \mathbb V$, if $f:A\to B$, then there  is an {\em internal} function ${}^+\!f:{}^*A\to B$ with the property that ${}^+\!f({}^*a)=f(a)$ for all $a\in A$.
\end{proposition}

\bproof
Suppose that  $f:A\to B$, where $A\in\mathbb U$ and that $B\in \mathbb V$ is internal. Let $T\in\mathbb U$ be transitive such that $B\in {}^*T$. Since ${}^*T$ is transitive, we have that $B\subseteq {}^*T$, so that $f:A\to {}^*T$. By definition of comprehensiveness, there is an internal function ${}^+\!g:{}^*A\to {}^*T$ such that ${}^+\!g({}^*a)=f(a)$ for all $a\in A$. Choose an arbitrary $b_0\in B$, and define 
\[{}^+\!f:=\Big\{(a,b)\in {}^*A\times B: \exists c\in{}^*T\;\Big((a,c)\in{}^+\!g\land(c\in B\to b=c)\land (c\in {}^*T-B\to b=b_0)\Big)\Big\}.\]
By Lemmas \ref{lemma_define_simple_notions}, \ref{lemma_internal_preservation} and the internal definition principle, ${}^+\!f$ is internal.  Moreover, if $a\in {}^*A$, then if ${}^+\!g(a)\in B$, we have $(a,{}^+\!g(a))\in {}^+\!f$, whereas if ${}^+\!g(a)\in {}^*T-B$, then $(a, b_0)\in {}^+\!f$. Thus ${}^+\!f$ is an internal function ${}^+\!f:{}^*A\to B$. Finally, if $a\in A$, then ${}^+\!g({}^*a)=f(a)\in B$, and hence $({}^*a, {}^+g({}^*a))\in {}^+\!f$, i.e. ${}^+\!f({}^*a)=f(a)$ when $a\in A$.
\eproof

The most useful application of comprehensiveness is in the countable case:
\begin{corollary} Suppose that $*:\mathbb U \to\mathbb V$ is a countably comprehensive transfer map, where $\mathbb U$ is a universe over a set $X$ that contains $\mathbb N$. Suppose also that  $B$ is internal and that $(b_n:n\in\mathbb N)$ is a sequence of members of $B$. Then there is an internal sequence $(b_n:n\in{}^*\mathbb N)$ which extends the sequence $(b_n:n\in\mathbb N)$.
\end{corollary}
\bproof Define $f:\mathbb N\to B:n\mapsto b_n$, and let ${}^+\!f:{}^*\mathbb N\to B$ be the map provided by Proposition \ref{propn_comprehensiveness}.
Then since ${}^*n=n$ for all $n\in\mathbb N$, we have  ${}^+\!f(n)={}^+\!f({}^*n)=f(n)=b_n$ for $n\in\mathbb N$. For $n\in {}^*\mathbb N-\mathbb N$, define $b_n:={}^+\!f(n)$.\eproof

\begin{theorem}\label{thm_saturation_implies_comprehensive} Suppose that $*:\mathbb U\to\mathbb V$ is a $\kappa$-saturated nonstandard embedding. Then $*$ is $\kappa$-comprehensive.
\end{theorem}
\bproof
Suppose that $*$ is $\kappa$-saturated. Let $A,B\in\mathbb U$, where $|A|<\kappa$, and suppose that $f:A\to {}^*B$. 
\[B_a:=\{g\in{}^*\mathcal P(A\times B):  \text{$g$ is a function}\land \dom(g)={}^*A\land ({}^*a, f(a))\in g\}.\]
By the internal definition principle, each $B_a$ is internal.

 Let $\mathcal B:=\{B_a:a\in A\}$, so that $|\mathcal B|<\kappa$. To apply $\kappa$-saturation, we must show that $\mathcal B$ has the f.i.p.
 Now for each $n\in \mathbb N$, we have that \[\mathbb U\vDash\forall x_1\dots x_n\in B\;\exists g\in \mathcal P(A\times B)\;\Big(\text{$g$ is a function}\land\dom(g)=A\land (a_1, x_1)\in g\land\dots\land (a_n,x_n)\in g\Big).\]
 Thus by transfer, setting $x_i:=f(a_i)$ for $a_1,\dots, a_n\in A$, there is $g\in{}^*\mathcal P(A\times B)$ such that
 \[\mathbb V\vDash \text{$g$ is a function}\land\dom(g)={}^*A\land ({}^*a_1, f(a_1))\in g\land\dots\land ({}^*a_n,f(a_n))\in g,\]
 and then clearly $g\in\bigcap_{i\leq n}B_{a_i}$. As $n\in\mathbb N$ is arbitrary, it follows that $\mathcal B$ has the f.i.p.
 
By $\kappa$-saturation,  there is an element ${}^+\!f\in\bigcap\mathcal B$. Then as ${}^+\!f\in{}^*\mathcal P(A\times B)$, it is an internal function ${}^*A\to {}^*B$ with the property that ${}^+\!f({}^*a)=f(a)$ for all $a\in A$.
\eproof

In the case of $\omega_1$-saturation, we also have the converse:
\begin{theorem} Suppose that $*:\mathbb U\to\mathbb V$ is a nonstandard embedding, where $\mathbb U$ is a universe over a set $X$ that contains $\mathbb N$. The following are equivalent:
\begin{enumerate}[(i)]
\item $*$ is $\omega_1$-saturated, i.e. every countable family $\mathcal A$ of internal sets with the f.i.p. has non-empty intersection.
\item $*$ is countably comprehensive.

\end{enumerate}
\end{theorem}
\bproof 
(i) $\Rightarrow$ (ii): This follows directly from Theorem \ref{thm_saturation_implies_comprehensive}.

(ii) $\Rightarrow$ (i): Suppose that $*$ is countably comprehensive, and that $\mathcal A:=\{A_n:n\in\mathbb N\}$ is a countable family of internal sets with the f.i.p. By replacing $A_n$ with $\bigcap_{m\leq n}A_n$, we may assume that the $A_n$ form a decreasing sequence of non-empty internal sets. As $A_0$ is internal, there is a transitive $T$ such that $A_0\in {}^*T$. Then as ${}^*T$ is transitive, we see that each $A_n$ is an internal subset of ${}^*T$, and thus $A_n\in{}^*\mathcal P(T)$.
Thus we have a map $f:\mathbb N\to{}^*\mathcal P(T):n\to A_n$. By countable comprehensiveness, there is an internal map ${}^+\!f:{}^*\mathbb N\to {}^*\mathcal P(T)$ such that ${}^+\!f(n)=A_n$ for every $n\in\mathbb N$. For $n\in{}^*\mathbb N-\mathbb N$, define $A_n:={}^+\!f(n)$. Then $(A_n:n\in{}^*\mathbb N)$ is an internal hypersequence of internal subsets of ${}^*T$ that extends the sequence $(A_n:n\in\mathbb N)$. Let \[X:=\{n\in{}^*\mathbb N: \forall k\in{}^*\mathbb N\;(k\leq n\to {}^+\!f(k)\supseteq{}^+\!f(n))\land {}^+\!f(n)\neq\varnothing\}\]
Then $X$ is internal, by the internal definition principle. By assumption, $\mathbb N\subseteq X\subseteq{}^*\mathbb N$. As $\mathbb N$ is not internal, there is $N\in {}^*\mathbb N-\mathbb N$ such that $N\in X$. Then $\bigcap_{n\in\mathbb N}A_n\supseteq A_N\neq\varnothing$. (One could also apply overflow to the set $X$ to deduce that $X$ has an infinite member.)

\eproof

\section{Questions of Existence}\label{section_question_existence}
\fancyhead[LO]{Questions of Existence}
This section is concerned with the existence of nonstandard embeddings, enlargements, comprehensive extensions, and polysaturated extensions. Since every polysaturated extension is both an enlargement (Lemma \ref{lemma_polysaturated_implies_enlargement}) and comprehensive (Theorem \ref{thm_saturation_implies_comprehensive}), the reader --- if there is one --- may want to read only Section \ref{section_ultrapower_proof_existence NSU} and either Section \ref{section_polysatured_ultrapowers} (which in turn depends on the existence of {\em good ultrafilters} --- cf. Appendix \ref{appendix_good_ultrafilters}) or Section \ref{section_polysaturated_ultralimits}.

\subsection{Existence of Nonstandard Frameworks}\label{section_existence_nonstandard_frameworks}

\begin{theorem}\label{theorem_exist_NSU} Let $V(X)$ be a superstructure over an infinite base set $X$. Then there exists a transfer map (i.e. a bounded elementary embedding) $V(X)\buildrel{*}\over\to V(Y)$ of $V(X)$ into some superstructure $V(Y)$ with base set $Y$ such that \begin{enumerate}[(i)]\item $^*X=Y$. \item $^*\varnothing=\varnothing$.
\item There is a countable $A\subseteq X$ such that $^\sigma A:=\{^*a:a\in A\}$ is a proper subset of $^*A$.
\end{enumerate}\end{theorem}

We will provide two proofs of this result, the first via ultrapowers, and the second by use of the Compactness Theorem of first-order logic.
\subsubsection{Ultrapower Proof}\label{section_ultrapower_proof_existence NSU}
We start from a superstructure $V(X)$ over a base set $X$. Let $\mathcal U$ be a countably incomplete ultrafilter over some index set $I$. 
The construction of the nonstandard framework proceeds over 9 steps. The first 8 steps will define a transfer map $*:V(X)\to V(Y)$ from the superstructure $V(X)$ to a superstructure $V(Y)$ over base set $Y$. These steps do not require $\mathcal U$ to be countably incomplete. In the $9^\text{th}$ step, the countable incompleteness of $\mathcal U$ is used to show that $*$ induces a nonstandard framework.
\vskip0.5cm

\noindent{\bf Step 1:} First we construct the base set $Y$:  Define an equivalence relation $\sim_\mathcal U$ on $X^I$  by \[f\sim_\mathcal U g\qquad\text{ if and only if }\qquad \{i\in I: f(i)=g(i)\}\in\mathcal U.\]
Define $Y$ to be the family of all equivalence classes
\[Y:=X^I/\sim_\mathcal U,\] i.e. $Y$ is just the ultrapower $X^I/\mathcal U$.
We shall assume that $I$ is chosen so that $Y$ is a base set for $V(Y)$. This can be done as follows: Suppose that $X$ is a set of rank $\beta$, and let $I$ be a set of rank  $\gamma\geq\beta+\omega$. Let $\mathcal U$ be an ultrafilter over $I$, and let $f:I\to X$. If $\gamma$ is a successor ordinal, $\gamma=\delta+1$, then $I$ has an element  $i_0$ of rank $\delta$, but no elements of higher rank. It follows easily that the rank of $f$ is $\delta+3=\gamma+2$. On the other hand, if $\gamma$ is a limit ordinal, then $\sup_{i\in I}\text{rank}(i)=\gamma$, but the supremum is not attained. Hence $f$ has rank $\gamma$. It therefore follows that every $f:I\to X$ has the same rank, namely either $\gamma+2$ if $\gamma$ is a successor ordinal, or $\gamma$ if it is limit. Now an element of an element of $Y=X^{I}/\mathcal U$ is precisely a function $f:I\to X$. Thus all elements of elements of $Y$ have the same rank. The argument in Example \ref{example_superstructures} now shows that $Y$ is a base set, i.e. that $y\cap V(Y)=\varnothing$ for all $y\in Y$.
\vskip0.3cm
\noindent{\bf Step 2:} Now we define a structure $(W,=_\mathcal U,\in_\mathcal U)$ which forms part of an intermediate step in the definition of the transfer map $*:V(X)\to V(Y)$: Define binary relations $=_\mathcal U$ and $\in_\mathcal U$ on $V(X)^I$ as follows:
\begin{align*} f &=_\mathcal U g  \quad\Longleftrightarrow\quad \{i\in I: f(i)=g(i)\}\in\mathcal U,\\ f&\in_\mathcal U g \quad\Longleftrightarrow\quad \{i\in I: f(i)\in g(i)\}\in\mathcal U.\end{align*}
If $f=_\mathcal U g$, we say that $f=g$ almost everywhere, and if $f\in_\mathcal U g$, we say that $f\in g$ almost everywhere.

 We can associate with each $a\in V(X)$ the constant mapping $c_a:I\to V(X)$ which takes the value constant $a$.
 For each $n\in\mathbb N$, let
\[W_n:=\{f\in V(X)^I: f\in_\mathcal U c_{V_n(X)}\},\] i.e.
a function $f:I\to V(X)$ belongs to $W_n$ if and only if $f(i)\in V_n(X)$ for almost all $i$. Equivalently, $f\in W_n$ if and only if there is $g:I\to V_n(X)$ such that $f =_\mathcal U g$.
Observe that since $V_0(X)=X$, we have that $W_0$ is essentially just $X^I$, i.e. \[W_0=\{f\in V(X)^{I}:f(i)\in X\text{ for almost all $i\in I$}\}.\] Clearly\[W_0\subseteq W_1\subseteq\dots\subseteq W_n\subseteq\dots\]
Define\[W=\bigcup_nW_n.\]

\vskip0.3cm\noindent{\bf Step 3:} We now show that there is a {\em unique} map $\cdot/\mathcal U:W\to V(Y)$ such that 
\begin{enumerate}[(i)]
\item $f/\mathcal U=\{g\in X^I: f=_\mathcal U g\}$ if $f\in W_0$, and
\item
$f/\mathcal U =\{g/\mathcal U: g\in W \land g\in_\mathcal U f\}$ if $f\in W-W_0$.
\item $f/\mathcal U\in V_n(Y)$ whenever $f\in W_n$.
\end{enumerate}

To begin with, define $f/\mathcal U:=\{g\in X^I: f=_\mathcal U g\}$  when $f\in W_0$. Then $f/\mathcal U=g/\sim_\mathcal U$ for any $g\in X^I$ such that $f=_\mathcal U g$.  Since $g/\sim_\mathcal U\in X^I/\mathcal U=Y$, it follows that $f/\mathcal U\in V_0(Y)=Y$ when $f\in W_0$.

Now proceed by induction. Suppose we have shown that, for each $m\leq n$, there is a unique map $h_m:W_m\to V_m(Y)$ such that \[h_m(f)=f/\mathcal U\text{ if }f\in W_0,\qquad h_m(f)=\{h_m(g): g\in W \land g\in_\mathcal U f\} \text{ if } f\in W_m-W_0.\] (Note that this condition makes sense, since if $f\in W_m$ and $g\in W$ is such that $g\in_\mathcal U f$, then \[\{g\in_\mathcal U c_{V_m(X)}\}\supseteq \{g\in_\mathcal U f\}\cap\{f\in_\mathcal U c_{V_m(X)}\}\in\mathcal U,\] by transitivity of $V_m(X)$,  so that $g\in W_m$ also.)

Define a map $h_{n+1}:W_{n+1}\to V_{n+1}(Y)$ as follows:
First, for  $f\in W_n$, define $h_{n+1}(f)=h_n(f)$, so that $h_{n+1}\restriction W_n = h_n$.  Next,  suppose that $f\in W_{n+1}-W_n$, and that $g\in W$ is such that $g\in_\mathcal U f$.
Since also $f\in_\mathcal U c_{V_{n+1}(X)}$ (by definition of $W_{n+1}$), it is easy to see that $g\in_\mathcal U c_{V_n(X)}$, and hence that $g\in W_n$. Hence $h_{n+1}(g) = h_n(g)$ has already been defined, and so we may define $h_{n+1}(f):=\{h_{n+1}(g): g\in W, g\in_\mathcal U f\}$. Note that then $h_{n+1}(f)\subseteq V_n(Y)$, so indeed $h_{n+1}(f)\in V_{n+1}(Y)$.

Clearly, $h_0\subseteq h_1\subseteq h_2\subseteq \dots$. We therefore define the required map $\cdot/\mathcal U$ by $\cdot/\mathcal U:=\bigcup_n h_n$. Then clearly $\cdot/\mathcal U$ satisfies statements (i),(ii), (iii). Uniqueness is easily established as well, as any two maps satisfying (i), (ii) must agree on $W_0$, and then, by induction, on all $W_n$.

\vskip0.3cm\noindent{\bf Step 4:}  Suppose that $f,g\in W$. We show that
\begin{enumerate}[(i)]\setcounter{enumi}{3}\item  $g\in_\mathcal U f$ if and only if $g/\mathcal U \in f/\mathcal U$, and 
\item $g=_\mathcal U f$ if and only if $g/\mathcal U=f/\mathcal U$ . \end{enumerate}

The definition of the map $\cdot/\mathcal U$ ensures that (i) holds (even for  $W_0$, since $Y$ is a base set for $V(Y)$). Thus we need only prove (ii).

Since $f, g\in W$, there is $n\in\mathbb N$ such that $f,g\in W_n$.
If $f, g\in W_0$, the statement of the Lemma is obviously true. Next, suppose that the statement holds for members of $W_m$ whenever $m\leq n$, and that $f, g\in W_{n+1}$.

If $f/\mathcal U=g/\mathcal U$ are such that $f\neq _\mathcal U g$, then $\{i\in I :f(i)=g(i)\}\not\in\mathcal U$, and so one of the sets  $\{i\in I: f(i)-g(i)\neq \varnothing\}, \{i \in I: g(i)-f(i)\neq \varnothing\}$ belongs to $\mathcal U$. Suppose the former. Define $h:I\to V(X)$ by letting $h(i)\in f(i)-g(i)$ if this set is non-empty, and setting $h(i)=\varnothing$ otherwise. Then $h\in W_n$. Clearly then $h/\mathcal U\in f/\mathcal U-g/\mathcal U$, so that $f/\mathcal U\neq g/\mathcal U$ --- contradiction. 

Conversely, suppose that $f, g\in W_{n+1}$ are such that $f=_\mathcal U g$. Then if $h/\mathcal U\in f/\mathcal U$, it follows that $h\in_\mathcal U f$, i.e.  that $\{i\in I: h(i)\in f(i)\}\in\mathcal U$, and thus that $\{i\in I: h(i)\in g(i)\}\supseteq\{i\in I: h(i)\in f(i)\}\cap\{i\in I: f(i)=g(i)\}\in\mathcal U$.  Hence also $h\in_\mathcal U g$, so $h/\mathcal U\in g/\mathcal U$. It follows that $f/\mathcal U\subseteq g/\mathcal U$. By symmetry $f/\mathcal U=g/\mathcal U$.

\vskip0.3cm\noindent{\bf Step 5:} Now we define the embedding $^*:V(X)\to V(Y)$ by
\[^*a:= c_a/\mathcal U.\]  To decompose this definition, define $\iota:V(X)\to W:a\mapsto c_a$. Then if $a\in V_n(X)$, $\iota(a)=c_a\in W_n$. Then $*:V(X)\to V(Y)$ is just the composition
\[V(X)\overset{\iota}{\hookrightarrow} W\overset{\cdot/\mathcal U}{\rightarrow}V(Y).\]
It is clear that $*$ is an embedding, i.e that if $a\neq b$ belong to $V(X)$, then $^*a\neq {}^*b$.

\vskip0.3cm\noindent{\bf Step 6:} We show that $^*\varnothing=\varnothing$ and ${}^*X=Y$.
The first statement is obvious. Observe that \[^*X=c_X/\mathcal U=\{f/\mathcal U: f\in W\land f\in_\mathcal U c_{V_0(X)}\} =\{f/\mathcal U: f\in W_0\}=X^I/\mathcal U=Y.\]

\vskip0.3cm\noindent{\bf Step 7:} 
We show that if $a\in V_n(X)$, then $^*a\in V_n(Y)$

This is clear if $n=0$, and thus holds for all individuals. We must therefore prove it for sets. Suppose now that it holds for $n$. If $A\in V_{n+1}(X)$ is a set, then $^*A=\{f/\mathcal U: f\in W\land f\in_\mathcal U c_A\}$. Thus if $f/\mathcal U\in {}^*A$, then $\{i\in I: f(i)\in A\}\in\mathcal U$, and thus $\{i\in I: f(i)\in V_n(X)\}\in\mathcal U$. It follows that $f\in W_n$, and thus that $f/\mathcal U\in V_n(Y)$ (by (iii) of Step 3). We therefore see that $^*A\subseteq V_n(Y)$, and thus that $^*A\in V_{n+1}(Y)$.

\vskip0.3cm\noindent{\bf Step 8:} Next, we prove that $^*$ is a transfer map, i.e. a bounded elementary embedding. 
We must show that for every  bounded $\mathcal L_\in $--formula $\varphi(x_1,\dots, x_n)$ and every $a_1,\dots, a_n\in V(X)$, we have
\[V(X)\vDash\varphi[a_1,\dots, a_n] \quad \Longleftrightarrow\quad V(Y)\vDash\varphi[{}^*a_1,\dots, {}^*a_n].\tag{$\dagger$}\] We will first show that  $f_1,\dots, f_n\in W$, then
\[V(Y)\vDash \varphi[f_1/\mathcal U,\dots, f_n/\mathcal U] \quad \Longleftrightarrow\quad \{i\in I : V(X)\vDash \varphi[f_1(i),\dots ,f_n(i)]\}\in\mathcal U.\tag{$\ddagger$}\]

The proof is similar to that of \L o\' s' Theorem, and proceeds by induction on the complexity of $\varphi$. If $\varphi$ is an atomic formula, i.e. of the form $\varphi(x_1,x_2)\equiv x_1 \in x_2$ or $\varphi(x_1,x_2)\equiv x_1=x_2$,  the result is an easy consequence of Step 4. In the former case, for example, $V(Y)\vDash f/\mathcal U\in g/\mathcal U$ if and only if $f\in_\mathcal U g$ if and only if $\{i\in I: V(X)\vDash f(i)\in g(i)\}\in\mathcal U$. 

Now suppose that $\varphi\equiv\psi\land\chi$, and that the result has been proved for $\psi,\chi$. Then since $\mathcal U$ is closed under intersections and supersets, we have
\begin{align*}
{}&\phantom{\Longleftrightarrow} V(Y)\vDash \varphi[f_1/\mathcal U,\dots, f_n/\mathcal U]\\ &\Longleftrightarrow V(Y)\vDash \psi [f_1/\mathcal U,\dots, f_n/\mathcal U]\quad\text{and}\quad V(Y)\vDash \chi[f_1/\mathcal U,\dots, f_n/\mathcal U]\\
&\Longleftrightarrow\{i\in I: V(X)\vDash \psi[f_1(i),\dots, f_n(i)]\}\in\mathcal U \quad\text{and}\quad \{i\in I: V(X)\vDash \chi[f_1(i),\dots, f_n(i)]\}\in\mathcal U \\
&\Longleftrightarrow  \{i\in I: V(X)\vDash \psi[f_1(i),\dots, f_n(i)]\}\cap \{i\in I: V(X)\vDash \chi[f_1(i),\dots, f_n(i)]\}\in\mathcal U\\
&\Longleftrightarrow \{i\in I: V(X)\vDash \varphi[f_1(i),\dots, f_n(i)]\}\in\mathcal U
\end{align*}
Next, suppose that $\varphi\equiv\lnot\psi$, and that the result has been proved for $\psi$. Then since $\mathcal U$ is an ultrafilter, we have that, for every $A\subseteq I$, either $A\in \mathcal U$ or $A^c\in\mathcal U$. Hence
\begin{align*}
{}&\phantom{\Longleftrightarrow} V(Y)\vDash \varphi[f_1/\mathcal U,\dots, f_n/\mathcal U]\\
&\Longleftrightarrow V(Y)\not\vDash \psi[f_1/\mathcal U,\dots, f_n/\mathcal U]\\
&\Longleftrightarrow \{i\in I: V(X)\vDash \psi[f_1(i),\dots, f_n(i)]\}\not\in\mathcal U\\
&\Longleftrightarrow\{i\in I: V(X)\vDash \lnot\psi[f_1(i),\dots, f_n(i)]\}\in\mathcal U
\end{align*}
Finally, if $\varphi\equiv (\forall y\in x_1)\psi(y,x_1,\dots, x_n)$, where $y, x_1,\dots, x_n$ are variables, then
\begin{align*}
{}&\phantom{\Longleftrightarrow} V(Y)\vDash \varphi[f_1/\mathcal U,\dots, f_n/\mathcal U]\\
&\Longleftrightarrow V(Y)\vDash \psi[g/\mathcal U,f_1/\mathcal U,\dots, f_n/\mathcal U]\quad\text{ for some } g/\mathcal U\in f_1/\mathcal U\\
&\Longleftrightarrow \{i\in I: V(X)\vDash \psi[g(i),f_1(i),\dots, f_n(i)]\}\in\mathcal U \quad\text{ for some } g\in_\mathcal U f_1\\
&\Longleftrightarrow\{i\in I: V(X)\vDash \varphi[f_1(i),\dots, f_n(i)]\}\in\mathcal U
\end{align*}
We have now proved $(\ddagger)$. In particular, we have\[V(Y)\vDash \varphi[c_{a_1}/\mathcal U,\dots, c_{a_n}/\mathcal U] \quad \Longleftrightarrow\quad \{i\in I : V(X)\vDash \varphi[c_{a_1}(i),\dots ,c_{a_n}(i)]\}\in\mathcal U,\]
But since $c_a(i)=a$ for all $i\in I$, the set $\{i\in I : V(X)\vDash \varphi[c_{a_1}(i),\dots ,c_{a_n}(i)]\}$ is either all of $I$, in which case $V(X)\vDash \varphi[a_1,\dots,a_n]$, or it is $\varnothing$, in which case $V(X)\vDash \lnot\varphi[a_1,\dots,a_n]$.
This proves that $(\dagger)$, i.e. that $\cdot/\mathcal U$ is a transfer map for the language $\mathcal L_\in$.

\vskip0.3cm{\bf Step 9:}
In order to complete the proof of Theorem \ref{theorem_exist_NSU}, it  remains to show that $*:V(X)\to V(Y)$ is a nonstandard framework over $X$, i.e. that there is a countable subset $A$ of $X$ such  that $^\sigma A:=\{{}^*a: a\in A\}$ is a proper subset of $^*A$. This is the only place where we need the fact that the ultrafilter $\mathcal U$ is countably incomplete.
 In fact, we can directly prove: \begin{center}If $A\in V(X)$ is an infinite set, then $^\sigma A\subsetneqq {}^*A$.\end{center}

Since $\mathcal U$ is countably incomplete, we can partition $I$ into a countable sequence $I_n$ of sets, none of whom belong to $\mathcal U$. Suppose now that $A\in V(X)$ is an infinite set, with distinct elements $a_0, a_1,a_2,\dots$. There is $m\in \mathbb N$ such that $A\in V_m(X)$. Define $f\in W_m$ by $f(i)=a_n$ whenever $i\in I_n$. Since  ${}^*A=c_A/\mathcal U=\{g/\mathcal U: g\in W\land g\in_\mathcal U c_A\}$, we have that $f/\mathcal U\in {}^*A$. However, if $f/\mathcal U={}^*b$ for some $b\in A$, then  $f/\mathcal U=c_b/\mathcal U$, so $\{i\in I: f(i)=b\}\in\mathcal U$, by Step 4. Since $f$ only takes the values $a_n$, it follows that $b=a_{n_0}$ for some $n_0\in\mathbb N$, But then $I_{n_0}=\{i\in I: f(i)=a_{n_0}\}\in\mathcal U$ --- contradicting the fact that $I_n\not\in\mathcal U$ for all $n$. Hence $f/\mathcal U\in {}^*A-{}^\sigma A$.\eproof
\subsubsection{Compactness Theorem Proof}

For this section, first recall the Compactness Theorem of first order logic --- cf. Theorem \ref{thm_compactness}.

\begin{definition}\rm Let $\mathfrak B:=(B,E)$ be a model for $\mathcal L_\in$. A submodel $\mathfrak A$ of $\mathfrak B$ is said to be a transitive submodel if whenever $a\in A, b\in B$ and $bEa$, then $b\in A$.
\endbox
\end{definition}

\begin{lemma}\label{lemma_transitive_submodel} Suppose that $\mathfrak A=(A,E)$ is a transitive submodel of $\mathfrak B=(B,E)$. Then  $\mathfrak A$ is a bounded elementary submodel of $\mathfrak B$.
\end{lemma}

\bproof We show, by induction on the complexity of $\varphi$, that for all $a_1,\dots,a_n\in A$ we have \[\mathfrak B\vDash \varphi[a_1,\dots, a_n] \quad\Leftrightarrow\quad \mathfrak A\vDash \varphi[a_1,\dots, a_n] .\]
This is obvious for atomic formulas, and  then easy to verify for the propositional connectives. Suppose therefore that $\mathfrak B\vDash\exists x\in a_1\psi[x,a_1,\dots, a_n]$, where $a_1,\dots, a_n\in A$, and $\psi$ is a bounded formula. Then there is $b\in B$ such that $bEa_1$ and $\mathfrak B\vDash\psi[b,a_1,\dots, a_n]$. As $\mathfrak A$ is a transitive submodel, we have $b\in A$, so by induction we have $\mathfrak A\vDash \psi[b,a_1,\dots, a_n]$, and hence $\mathfrak A\vDash \exists x\in a_1\psi[x,a_1,\dots, a_n]$. As $\mathfrak A$ is a submodel of $\mathfrak B$, it is obvious that $\mathfrak A\vDash \exists x\in a_1\psi[x,a_1,\dots, a_n]$ implies $\mathfrak B\vDash \exists x\in a_1\psi[x,a_1,\dots, a_n]$.\eproof

Given a model $\mathfrak B=(B,E)$, and  an $X\in B$, we have that --- from $\mathfrak B$'s point of view --- the element $X$ is a base set if and only if \[\mathfrak B\vDash \text{BASE}[X],\quad\text{where}\quad \text{BASE}(x)\equiv\forall y\in x\;\forall z\in y\;(z\neq z),\] i.e. there are no $c,b\in B$ such that $cEbEX$: every $b\in B$ such that $bEX$ looks like an atom to $\mathfrak B$.

Given that $\mathfrak B\vDash\text{BASE}[X]$, we want to {\em truncate} $\mathfrak B$ by removing all sets which are not at a finite level over $X$. Recall the formulas $\nu_n(a,X)\equiv\varphi_{6,n}(X,a)$ which assert that  $a\in V_n(X)$:
\[\nu_0(y,x)\equiv y\in x,\qquad \nu_{n+1}(y,x)\equiv \nu_n(y,x)\lor\forall z\in y\;\nu_n(z,x).\]
Thus for a superstructure $V(X)$ over a base set $X$, we have \[V_n(X)=\{a\in V(X): (V(X),\in)\vDash \nu_n[a,X]\}, \qquad (V(X),\in)\vDash \text{BASE}[X].\] To define the truncation $\mathfrak A$ of  $\mathfrak B$ over $X\in B$, we imitate: Define\[A:=\{a\in B: \text{ there is }n<\omega\text{ such that }\mathfrak B\vDash \nu_n[a,X]\},\] and let $\mathfrak A=(A,E)$ be the resulting submodel of $\mathfrak B$.

\begin{lemma}\label{lemma_truncation_bd_elelmentary_submodel} Let $X$ be a base set, let $\mathfrak B=(B,E)$ be a bounded elementary extension of $(V(X),\in)$, and let $\mathfrak A$ be the truncation of $\mathfrak B$ over $X$. Then $\mathfrak A$ is a transitive submodel of $\mathfrak B$. Hence \[(V(X),\in)\preceq_b\mathfrak A\preceq_b\mathfrak B.\] In addition, the truncation of $\mathfrak A$ over $X$ is $\mathfrak A$ itself.
\end{lemma} 

\bproof
Observe that since $\text{BASE}(x)$ is a bounded formula and $V(X)\vDash \text{BASE}[X]$, we also have $\mathfrak B\vDash \text{BASE}[X]$. Now if $a\in A$, then there is a least $n<\omega$ such that $\mathfrak B\vDash \nu_n[a,X]$. If then $b\in B$ is such that $bEa$, then we cannot have $n=0$, because $\mathfrak B\vDash\text{BASE}[X]$, and hence necessarily $\mathfrak B\vDash \nu_{n-1}[b, X]$. It follows that $b\in A$, and thus that $\mathfrak A$ is a transitive submodel of $\mathfrak B$.

Furthermore, if $v\in V(X)$, then $V(X)\vDash \nu_n[v,X]$ for some $n$, As each $\nu_n(y,x)$ is a bounded formula, we have $\mathfrak B\vDash\nu_n[v,X]$, from which it follows that $v\in A$, and thus $(V(X),\in)$ is a transitive submodel of $\mathfrak A$.

By Lemma \ref{lemma_transitive_submodel}, we see that $(V(X),\in)\preceq_b(A,E)\preceq_b(B,E)$. 

Finally if $a\in A$, then $\mathfrak B\vDash \nu_n[a,X]$ for some $n<\omega$. Since $\mathfrak A\preceq_b\mathfrak B$, we have $\mathfrak A\vDash \nu_n[a,X]$, from which it follows that $\mathfrak A$ is its own truncation over $X$.\eproof

\begin{theorem} \label{thm_Mostowski_Collapse}{\rm (Mostowski Collapse)} Let $\mathfrak A=(A,E)$  be a model of $\mathcal L_\in$ with an element $X\in A$ such that
\begin{enumerate}[(i)]\item $\mathfrak A\vDash \text{\rm BASE}[X]$.
\item $\mathfrak A$ is its own truncation over $X$, i.e. for every $a\in A$ there is $n<\omega$ such that $\mathfrak A\vDash \nu_n[a, X]$.
\item $\mathfrak A$ is extensional over $X$: \[\mathfrak A\vDash\forall u\;\forall v\;(u\in x\lor v\in x\lor(u=v \leftrightarrow \forall z\;(z\in u\leftrightarrow z\in v))[X],\] i.e. two {\em sets} relative to $\mathfrak A$ are equal if and only if ($\mathfrak A$ thinks that) they have the same elements.
\item The set $Y:=\{a\in A:aEX\}$ is a base set.
\end{enumerate}
Then there is a unique bounded elementary embedding $h:\mathfrak A\hookrightarrow (V(Y),\in)$ with the properties that:
\begin{enumerate}[1)]\item $h(a)=a$ for all $a\in Y$,
\item $h(X)=Y$,
\item $\ran\;h$ is a transitive subset of $V(Y)$.
\end{enumerate}
\end{theorem}

\bproof For $n<\omega$, let $A_n:=\{a\in A:\mathfrak A\vDash\nu_n[a,X]\}$. Then $(A_n)_n$ is an increasing sequence of sets, and as $\mathfrak A$ is its own truncation over $X$, we have that $\bigcup_nA_n=A$. By definition of $\nu_0$, we have that $A_0=Y$. We now define $h\restriction A_n$ by induction, and then take $h:=\bigcup_nh\restriction A_n$.

For $a\in A_0=Y$, put $h(a)=a$.

Now suppose that $h\restriction A_n$ has been defined, and that $a\in A_{n+1}-A_n$. If $bEa$, then by definition of $\nu_{n+1}$ we must have $b\in A_n$, so that $h(b)$ is already defined. Thus put $h(a):=\{h(b):bEa\}$.

This completes the definition of $h$.

Now by definition of $h$ we have that $h(a)=a$ for all $a\in Y$. Furthermore, \[h(X)=\{h(a):aEX\}=\{h(a):a\in Y\}=\{a:a\in Y\}=Y.\]
Next, we show that $\ran\;h$ is transitive. Suppose that $b'\in a'\in \ran\;h$. We must show that $b'\in\ran\;h$. Now as $a'\in\ran\;h$, there exists $a\in A$ such that $h(a)=a'$, and so $b'\in h(a)$. By definition of $h$, therefore, there must be $b\in A$ such that $bEa$ and $h(b)=b'$. In particular, we see that $b'\in\ran(h)$. This demonstrates that $\ran\;h$ is a transitive set.

For $n<\omega$, consider the statement:\[
P_n\equiv\qquad \text{The restriction } h\restriction A_n\text{ is one-to one, and } h[A_n]\subseteq V_n(Y).\]
It is clear that $P_0$ holds. Now suppose that $P_n$ holds, and that $a\in A_{n+1}$. Then $b\in A_n$ whenever $bEa$, and so  $h(a)=\{h(b): bEa\}\subseteq V_n(Y)$, from which it follows that $h(a)\in V_{n+1}(Y)$, and thus that $h[A_{n+1}]\subseteq V_{n+1}(Y)$.  Moreover, if $a,a'\in A_{n+1}$ are such that $h(a)=h(a')$. Then for every $bEa$ there is $cEa'$ such that $h(b)=h(c)$. But necessarily $b,c\in A_n$, and as $P_n$ holds, we have  that $h(b)=h(c)$ implies $b=c$. Thus $bEa$ if and only if $bEa'$, and hence by the extensionality property we have $a=a'$. It follows that $h\restriction A_{n+1}$ is one-to one, and thus that $P_{n+1}$ holds. Thus, by induction, $P_n$ holds for all $n<\omega$.

It follows that $\ran\;\subseteq V(Y)$ and that $h$ is one-to-one.

Next note that, by definition of $h$, we see that $bEa$ implies $h(b)\in h(a)$.  Conversely, if $h(b)\in h(a)$, then $h(b)$ must be equal to $h(c)$ for some $cEa$. But as $h$ is one-to-one, we have that $b=c$. Hence $h(b)\in h(a)$ implies $bEa$. It follows that $h:(A,E)\hookrightarrow (V(Y),\in)$ is an embedding. As $(\ran\;h, \in)$ is a transitive submodel of $(V(Y),\in)$, it follows by Lemma \ref{lemma_transitive_submodel} that\[ (A,E)\cong (\ran\;h,\in)\preceq_n(V(Y),\in).\]
In particular, $h:\mathfrak A\to (V(Y),\in)$ is a bounded elementary embedding.

It remains to show that $h$ is the unique bounded elementary embedding with the properties 1)-3). Suppose that $h'$ is another such map. Again, we use induction to show that $h\restriction A_n = h'\restriction A_n$. This is clear if $n=0$. Now suppose that $h\restriction A_n=h'\restriction A_n$ and that $a\in A_{n+1}$. By definition of $h$ we see that $x\in h(a)$ implies $x=h(b)$ for some $bEa$. But then  $h'(b)\in h'(a)$. As necessarily $b\in A_n$, we have $h(b)=h'(b)$, and thus $x\in h'(a)$.It follows that $h(a)\subseteq h'(a)$. Conversely, if $x\in h'(a)$, then as $\ran\;h'$ is transitive, there is $b$ such that $h'(b)=x$, i.e. $h'(b)\in h'(a)$. But as $h'$ is a bounded elementary embedding, it follows that $bEa$, so that $b\in A_n$, and hence $h(b)=h'(b)$. Now $bEa$ implies $h(b)\in h(a)$, and thus $x\in h(a)$. It follows that $h'(a)\subseteq h(a)$, i.e. that $h'(a)=h(a)$, and hence that $h\restriction A_{n+1}=h'\restriction A_{n+1}$.

\eproof

We are now in a position to prove Theorem \ref{theorem_exist_NSU}:

\bproof
Suppose that $(V(X),\in)$ is a superstructure over an  infinite base  set $X$. Consider the language $\mathcal L_{V(X)}=\mathcal L\cup\{c_u:u\in V(X)\}$, and let $\Delta$ be the elementary diagram of $V(X)$ --- cf. Definition \ref{defn_theory_diagram}. For each infinite set $U\subseteq X$ in $V(X)$, let $d_U$ be a new constant symbol, and let $\Sigma=\Delta\cup\{\varphi_{U,u}: U\in V(X)\text{ an infinite subset of } X, u\in V(X)\}$, where
\[\varphi_{U,u}\equiv d_U\in c_U\land d_U\neq c_u.\]
If $\Sigma'\subseteq \Sigma$ is finite, then it refers to at most finitely many $U,u$, and hence $(V(X),\in)$ can be expanded to a model $\Sigma'$, where each $c_u$ is interpreted as the element $u\in V(X)$, and $d_U$ is interpreted to be a member of the set $U$. Hence $\Sigma$ is consistent, and therefore has a model $\mathfrak B=(B,E)_{c_u, d_U}$. As this is a model of the elementary diagram, of $(V(X),\in)$ we see that we have an elementary extension $(V(X),\in)\preceq(B,E)$ --- cf. Lemma \ref{lemma_diagram_embedding}. Moreover, if $b_U\in B$ is the interpretation of the constant $d_U$, then we have $b_U E U$, but $b_U\neq u$ for any $u\in V(X)$, i.e. for every infinite $U\subseteq  X$ there is $b\in B- V(X)$ such that $bEU$. By renaming, we may choose the set $B$ so that $Y:=\{b\in B: bEX\}$ is a base set. 

Let $\mathfrak A$ be the truncation of $\mathfrak B$ over $X$. By Lemma \ref{lemma_truncation_bd_elelmentary_submodel} , we have $(V(X),\in)\preceq_b\mathfrak A\preceq_b\mathfrak B$, and $\mathfrak A$ is its on truncation over $X$. Furthermore, as $\text{BASE}(x)$ is a bounded formula and $V(X)\vDash\text{BASE}[X]$, it follows that $\mathfrak A\vDash\text{BASE}[X]$. In addition, $Y=\{b\in B: \nu_0[b,X]\}\subseteq A$, so $Y=\{a\in A:aEX\}$ is a base set.

Next, we show that $\mathfrak A$ is extensional over $X$. Suppose that $a,b\in A$ are such that $a,b\not\in X$ and that $a\neq b$. The model $(V(X),\in)$ is certainly extensional over $X$, and as it is an elementary submodel of $\mathfrak B$ (i.e. not merely a bounded elementary submodel), it follows that $\mathfrak B$ is extensional over $X$. Since $a\neq b$ are members of $B$, we have  $\mathfrak B\vDash \exists z\:(z\in a\leftrightarrow z\not\in b)$. But this can also be written as a bounded formula: $\mathfrak B\vDash \exists z\in a(z\not\in b)\lor\exists z\in b(z\not\in a)$. Hence also $\mathfrak A\vDash \exists z\in a(z\not\in b)\lor\exists z\in b(z\not\in a)$, from which it follows that there is $c\in A$ such that $cEa$ if and only if $\lnot cEb$. This shows that $\mathfrak A$ is extensional over $X$.

We are now able to apply Theorem \ref{thm_Mostowski_Collapse} to deduce that there is a bounded elementary extension $h:\mathfrak A\hookrightarrow (V(Y),\in)$. Let $*:=h\restriction V(X)$. Then $*:(V(X),\in)\hookrightarrow (V(Y),\in)$ is a composition of two bounded elementary embeddings, and thus a bounded elementary embedding, with $^*x=x$ for all $x\in X$, and $^*X=Y$. Finally, if $U\subseteq X$ is infinite, there is $b\in B- V(X)$ such that $bEU$.  As $V(X)\vDash\forall u\in c_U\;(u\in c_X)[U,X]$, we also have that $bEU$ implies $bEX$ for all $b\in B$, and hence $bEU$ implies $b\in A$, by definition of truncation. Thus $h(b)$ is defined, and  $h(b)\in h(U)={}^*U$. The fact that $b\not \in X$ means that $h(b)\not\in \{h(x):x\in X\}=\{^*x:x\in X\}$, so $\{^*x:x\in U\}$ is a proper subset of $^*U$.\eproof

\subsection{Existence of Enlargements} 
\begin{lemma} Suppose that $*:\mathbb U\to\mathbb V$ is a $\kappa$-enlargement. Then $|\mathbb N^*|\geq \sup\{|A|:A\in\mathbb U, |A|<\kappa\}$.
\end{lemma}
\bproof Suppose that $A\in\mathbb U$ has cardinality $<\kappa$. By Theorem \ref{thm_enlargement_equiv}, there is a hyperfinite set $B$ such that ${}^\sigma A\subseteq B\subseteq {}^*A$. As $B$ is hyperfinite, there is a $n\in{}^*\mathbb N$ and a bijection $f:\{0,1,\dots, n-1\}\to B$, and thus an injection $h: B\hookrightarrow {}^*\mathbb N$.  Thus ${}^*\mathbb N$ has a subset of cardinality $|B|$, and thus one of cardinality $|A|=|{}^\sigma A|\leq |B|$. Hence $|{}^*\mathbb N|\geq |A|$ for any $A\in\mathbb U$ with $|A|<\kappa$.
\eproof

Using the above lemma, it can be seen that the  ultrapower construction does not automatically provide enlargements. For example, consider the ultrapower construction $*:V(X)\to V(Y)$ over an infinite base set $X$, with $Y=X^I/\sim_\mathcal U$ and $I=\mathbb N$.  we may assume that $\mathbb N\subseteq X$. Then $^*\mathbb N$ is of the form $c_\mathbb N/\mathcal U=\{g/\mathcal U:g\in_\mathcal U\mathbb N\}$, and so $|{}^*\mathbb N|\leq |\mathbb N^\mathbb N|=2^{\aleph_0}$. Now $\mathbb N,\mathcal P(\mathbb N),\mathcal P\mathcal P(\mathbb N)\dots$ are all members of $\mathbb U$, so if $\kappa>(2^{\aleph_0})^+$, then $\sup\{|A|:A\in\mathbb U, |A|<\kappa\}>2^{\aleph_0}$. Hence if  $\kappa>(2^{\aleph_0})^+$, then $*$ cannot be a $\kappa$-enlargement.

In order to use the ultrapower construction to obtain an enlargement, we have to be a bit more careful about the set $I$. 

\begin{theorem}\label{thm_enlargements_exist}
Given a base set $X$, let $I:=\mathcal P^{<\omega}(V(X))$ be the family of all finite subsets of $V(X)$. For each $a\in I$, define $I_a:=\{b\in I: a\subseteq b\}$. The family $\{I_a:a\in I\}$ has the f.i.p., so there is therefore an ultrafilter $\mathcal U$ over $I$ such that each $I_a\in\mathcal U$. Let $Y:=X^I/\mathcal U$. The ultrapower construction $*:V(X)\overset{\iota}{\hookrightarrow}W\overset{\cdot/\mathcal U}{\to}V(Y)$ given in the proof of Theorem \ref{theorem_exist_NSU} is an enlargement.\end{theorem}

\bproof Observe that if $a_1,\dots, a_n \in I$, then $a_1\cup\dots\cup a_n$, being finite, is a member of $I_{a_1}\cap\dots\cap I_{a_n}$, which shows that the family $\{I_a:a\in I\}$ has the f.i.p.
There is therefore an ultrafilter $\mathcal U$ over $I$ such that $I_a\in\mathcal U$ for every finite subset $a$ of $V(X)$. 

\underline{Method 1:} Suppose that $\mathcal B\in V(X)$ is a family of sets with the f.i.p. Define $f: I \to V(X)$ as follows: If $a\in I$, then $a\cap\mathcal B\in I$ also. If $a\cap\mathcal B\neq \varnothing$, choose $f(a)\in \bigcap (a\cap\mathcal B)$; else, put $f(a)=\varnothing$. Observe that if $\mathcal B\in V_n(X)$, then $f(a)\in V_n(X)$, so that $f\in W_n$. We claim that $f/\mathcal U\in \bigcap{}^\sigma\mathcal B$.

For if $B\in\mathcal B$, then $\{B\}\in I$, and hence $I_{\{B\}}=\{a\in I:B\in a\}\in\mathcal U$. Now \[a\in I_{\{B\}}\Rightarrow B\in a\Rightarrow a\cap \mathcal B\neq \varnothing\Rightarrow f(a)\in \bigcap (a\cap \mathcal B)\subseteq B,\] and hence $I_{\{B\}}\subseteq\{a\in I: f(a)\in B\}$. It follows that $f/\mathcal U\in c_B/\mathcal U={}^*B$. As $B\in\mathcal B$ is arbitrary, it follows that $f/\mathcal U\in \bigcap_{B\in\mathcal B}{}^*B=\bigcap {}^\sigma\mathcal B$. As $\mathcal B$ is arbitrary, it follows that $\bigcap{}^\sigma\mathcal B\neq \varnothing$ for all $\mathcal B\in V(X)$ with the f.i.p, and thus by Theorem \ref{thm_fip_element}, the extension $*:V(X)\to V(Y)$ is an enlargement.

\underline{Method 2:}  We use the hyperfinite approximation property: Let $B\in V(X)$, and define $g:I\to \mathcal P^{<\omega}(B)$ by $g(a):=a\cap B$. Then $\{a\in I: g(a)\in \mathcal P^{<\omega}(B)\}=I\in\mathcal U$, and thus $g/\mathcal U\in {}^*\mathcal P^{<\omega}(B)$, i.e. $g/\mathcal U$ is a hyperfinite subset of ${}^*B$. Now if $b\in B$ and $b\in a$, then $b\in a\cap B=g(a)$, and hence for $b\in B$ we have \[I_{\{b\}}:=\{a\in I: b\in a\}=\{a\in I:b\in g(a)\}.\] As $I_{\{b\}}\in\mathcal U$, we see that ${}^*b=c_b/\mathcal U\in g/\mathcal U$. Thus with $A:=g/\mathcal U\in {}^*\mathcal P^{<\omega}(B)$, we have ${}^\sigma B\subseteq A\subseteq {}^*B$. It follows that for every $B\in V(X)$ there is a hyperfinite set $A$ such that ${}^\sigma B\subseteq A\subseteq {}^*B$, and this is equivalent to $*$ being an enlargement, by Theorem \ref{thm_enlargement_equiv} .\eproof

\subsection{Existence of Comprehensive Transfer Maps}

Every nonstandard embedding obtained from an ultrapower construction is comprehensive:

 \begin{theorem} Let $V(X)$ be a superstructure over $X$, and let $*:V(X)\to V({}^*X)$ be the transfer map provided by an ultrapower construction, as in  Section \ref{section_ultrapower_proof_existence NSU}. Then $*$ is comprehensive.
 \end{theorem}
 
 \bproof Suppose that $V({}^*X)$ is obtained from $V(X)$ via an ultrafilter $\mathcal U$ over a set  $I$. Suppose further that $A, B\in V(X)$, and that $f:A\to {}^*B$. We 
must show that there is an internal map ${}^+\!f:{}*A\to {}^*B$ with the property that ${}^+\!f({}^a)=f(a)$ for all $a\in A$.

For $a\in V(X)$, let $c_a$ denote the constant map $c_a:I\to V(X):i\mapsto a$. Then ${}^*a:=c_a/\mathcal U$ --- see Step 5 of the ultrapower proof of Theorem \ref{theorem_exist_NSU}. Also, let $\rho_{f(a)}:I\to V(X)$ be such that $f(a)=\rho_{f(a)}/\mathcal U\in {}^*B$. As $V({}^*X)\vDash \rho_{f(a)}/\mathcal U\in c_B/\mathcal U$ we may, via \L os' Theorem, and without loss of generality, assume that $\rho_{f(a)}(i)\in B$ for all $i\in I$. For $i\in I$, define $f_i:A\to B:a\mapsto \rho_{f(a)}(i)$.
Now let $F:I\to V(X): i\mapsto f_i$. 

Observe first that if $A,B\in V_n(X)$, then each $f_i\in V_{n+2}(X)$, so that $F\in V_{n+2}(X)^I$, i.e. $F\in W_{n+2}$ has finite rank, and so $F/\mathcal U\in V({}^*X)$.
Since $V(X)\vDash F(i):c_A(i)\to c_B(i)$ holds for all $i\in I$, we see by \L os' Theorem that $V({}^*X)\vDash F/\mathcal U:c_A/\mathcal U\to c_B/\mathcal U$, where we use Lemma \ref{lemma_define_simple_notions}(f). Moreover, $F/\mathcal U\in c_{V_{n+2}(X)}/\mathcal U= {}^*V_{n+2}(X)$.

Thus if we define ${}^+\!f:=F/\mathcal U$, then we immediately see that ${}^+\!f:{}^*A\to {}^*B$ is internal.

Finally, ${}^+\!f({}^*a)= F/\mathcal U(c_a/\mathcal U)$. Now in $V(X)$, we have  that $F(i)(c_a(i))= f_i(a)=\rho_{f(a)}(i)$ for all $i\in I$. By \L os' Theorem, therefore,  we have that $F/\mathcal U(c_a/\mathcal U)=\rho_{f(a)}/\mathcal U$ holds in $V({}^*X)$, i.e. that ${}^+\!f({}^*a)=f(a)$.
\eproof

Since only special types of ultrapower constructions provide enlargements, not every comprehensive extension is an enlargement.

In Section \ref{section_saturation}, we showed that countably comprehensive transfer maps are $\omega_1$-saturated, and vice versa.  Example \ref{example_ultrafilter_omega1_good}  will show that every ultrafilter is $\omega_1$--good, which, combined with Theorem \ref{thm_transfer_map_good_ultrafilter_saturated}, shows that every transfer map induced by an ultrapower construction is $\omega_1$-saturated, and thus countably comprehensive.

\subsection{Existence of Polysaturated Extensions via Ultrapowers}\label{section_polysatured_ultrapowers}

Recall the following definitions from basic model theory:

\begin{definition}\rm
\begin{enumerate}[(a)]\item
Let $\mathcal L$ be a first-order language. Given a model $\mathfrak A$ of the language $\mathcal L$, we denote its universe $A$. If $X\subseteq A$, we denote by $\mathcal L_X$ the language $\mathcal L$ augmented with a set of new constant symbols $\{c_a:a\in X\}$. We expand the $\mathcal L$--structure $\mathfrak A$ to a $\mathcal L_X$-structure $(\mathfrak A,a)_{a\in X}$, where the constant symbol $c_a$ is interpreted as the element $a$  in $(A,a)_{a\in X}$.
\item A set of formulas $\Sigma(x)$ in one free variable $x$ is said to be {\em satisfiable} in a model $\mathfrak A$ if and only if there is $b\in A$ such that $\mathfrak A\vDash\Sigma[b]$.\\ $\Sigma(x)$ is said to be {\em finitely satisfiable} in $\mathfrak A$ if every finite subset of $\Sigma$ is satisfiable.
\item
Suppose that $\mathfrak A$ is a model of the language $\mathcal L$, and that $\kappa$ be an infinite cardinal. $\mathfrak A$ is said to be {\em $\kappa$-saturated} if and only if and only if the following condition holds: Given  a subset $X\subseteq A$ with $|X|<\kappa$ and a set of $\mathcal L_X$-formulas  $\Sigma(x)$ in one free variable $x$, then $\Sigma(x)$ is satisfiable in $(\mathfrak A, a)_{a\in X}$ whenever it is finitely satisfiable.
\end{enumerate}
\endbox
\end{definition}

We first deal with a simple case:
\begin{theorem}\label{thm_ultrapower_omega1_saturated}
	Suppose that $\mathcal L$ is a countable language, and that $\mathcal U$ is a countably incomplete ultrafilter over a set $I$. Then any ultraproduct $\prod_I\mathfrak A_i/\mathcal U$ is $\omega_1$-saturated
\end{theorem}
\bproof
As $\mathcal U$ is countably incomplete, there is a decreasing chain $(I_n)_{n\in\mathbb N}$ of elements of $\mathcal U$ such that $\bigcap_nI_n=\varnothing$. Without loss of generality, we may assume $I_0=I$.

We first show that the following claim holds:
\vskip0.2cm{\bf Claim:} If $\mathcal L$ is a countable language, and if  $\Sigma(x)$ is a set of $\mathcal L$-formulas which is finitely satisfiable in an ultraproduct $\prod_I\mathfrak A_i/\mathcal U$, then $\Sigma(x)$ is satisfiable.
\vskip0.2cm
So suppose that $\Sigma(x)$ is finitely satisfiable in $\prod_I\mathfrak A_i/\mathcal U$. Since $\mathcal L$ is countable, so is $\Sigma(x)$, and hence we can enumerate it:
\[\Sigma(x)=\{\sigma_n(x):n\geq 1\}.\]
Define \[U_0:=I,\qquad U_{n}=I_n\cap\{i\in I: \mathfrak A_i\vDash\exists x\bigwedge_{1\leq m\leq n}\sigma_m(x)\},\]so that each $U_n\in\mathcal U$, by \L os' Theorem. Now define $N(i):=\max\{n:i\in U_n\}$, and choose $a\in\prod_IA_i$ as follows:
If $N(i)=0$, let $a(i)\in A_i$ be arbitrary. Else, choose $a(i)$ so that
\[\mathfrak A_i\vDash \bigwedge_{1\leq m\leq N(i)}\sigma_m[a(i)].\]
Now note that if $n\geq 1$ and $i\in U_n$, then $N(i)\geq n$, and hence $\mathfrak A_i\vDash \sigma_n[a(i)]$. It follows that\[U_n\subseteq\{i:\mathfrak A_i\vDash \sigma_n[a(i)]\},\] As $U_n\in\mathcal U$, it follows that $\prod_I\mathfrak A_i/\mathcal U\vDash \sigma_n[a/\mathcal U]$. As $n\geq 1$ is arbitrary, the element $a/\mathcal U$ satisfies $\Sigma(x)$ in $\prod_I\mathfrak A_i/\mathcal U$. This proves the Claim.

Now to prove $\omega_1$-saturation: Suppose that $X=\{a_n/\mathcal U:n\in\mathbb N\}$ is a countable set of elements of $\prod_I\mathfrak A_i/\mathcal U$. Let $\Sigma(x)$ be a set of formulas of $\mathcal L(X)$ which is finitely satisfiable in $(\prod_I\mathfrak A_i/\mathcal U, (a_n/\mathcal U)_n)$. We must show that $\Sigma(x)$ is satisfiable in $(\prod_I\mathfrak A_i/\mathcal U, (a_n/\mathcal U)_n)$. Now the expanded language $\mathcal L(X)$ is still countable, and it is easy to verify that \[(\prod_I\mathfrak A_i/\mathcal U, (a_n/\mathcal U)_n)=\prod_I(\mathfrak A_i, (a_n(i))_n)/\mathcal U,\]i.e. $(\prod_I\mathfrak A_i/\mathcal U, (a_n/\mathcal U)_n)$ is an ultraproduct. By the Claim, in any ultraproduct modulo $\mathcal U$ of structures that interpret a countable language, every finitely satisfiable $\Sigma(x)$ is satisfiable. The result now follows by applying the claim to the ultraproduct $\prod_I(\mathfrak A_i, (a_n(i))_n)/\mathcal U$ of structures interpreting the countable language $\mathcal L(X)$.
\eproof

Suppose now that we have a family$\{\mathfrak A_i:i\in I\}$ of $\mathcal L$-structures, indexed by a set $I$, and an ultrafilter $\mathcal U$ on $I$. We seek conditions on $\mathcal U$ which will ensure that the ultraproduct $\prod_IA_i/\mathcal U$ is $\kappa$-saturated for $\kappa\geq\omega_1$. Moreover, as in Theorem \ref{thm_ultrapower_omega1_saturated}, we seek a condition on $\mathcal U$, i.e. one which is independent of the models $\mathfrak A_i$. 

\vskip0.5cm\noindent
{\bf Condition S:} Whenever $\Sigma(x)$ with $|\Sigma(x)|<\kappa$ is finitely satisfiable in an ultraproduct modulo $\mathcal U$, then it is  satisfiable. 
\vskip0.5cm

This condition is independent of the language $\mathcal L$ or the models that make up the ultraproduct.

Observe that if $\mathcal U$ satisfies (S) then any ultraproduct modulo $\mathcal U$ interpreting a language $\mathcal L$ of cardinality $<\kappa$ is $\kappa$-saturated. To see this, suppose that $X\subseteq\prod_IA_i/\mathcal U$ has $|X|<\kappa$. Then $\mathcal L_X$ is also a language of cardinality $<\kappa$. Automatically, therefore, any set of $\mathcal L_X$-formulas $\Sigma(x)$ has $|\Sigma(x)|<\kappa$. Condition (S) then immediately yields that any set of $\mathcal L_X$-formulas $\Sigma(x)$ which is finitely satisfiable in $\prod_I\mathfrak A_i/\mathcal U$ is satisfiable in  $\prod_I\mathfrak A_i/\mathcal U$.

We seek a property of $\mathcal U$ which guarantees that (S) holds.

Thus let $\Sigma(x)$ be a set of formulas of cardinality $<\kappa$ which is finitely satisfiable in $\prod_I\mathfrak A_i/\mathcal U$. Then if $\Theta$ is a finite subset of $\Sigma$, the set $\{i\in I: \mathfrak A_i\vDash \exists x(\bigwedge\Theta)\}$ belongs to $\mathcal U$, by \L os' Theorem. We thus have a map\[p:\mathcal P^{<\omega}(\Sigma)\to\mathcal U:\Theta\mapsto \{i\in I: \mathfrak A_i\vDash \exists x(\bigwedge\Theta)\}.\tag{$\star$}\]
Thus   $i\,\in p(\Theta)$ if and only if  $\mathfrak A_i\vDash\exists x\bigwedge\Theta$, i.e. if and only if 
$\Theta$ is satisfiable in $\mathfrak A_i$.

It should be clear that $\Theta\subseteq\Theta'\Rightarrow p(\Theta)\supseteq p(\Theta')$. Equivalently
\[p(\Theta\cup\Theta')\subseteq p(\Theta)\cap(\Theta')\]
\begin{lemma}\label{lemma_support_saturation} If there exists a sequence $(\Phi_i)_{i\in I}$ in $\mathcal P^{<\omega}(\Sigma)$ such that \begin{enumerate}[(i)]\item $i\in p(\Phi_i)$ for all $i\in I$, and
\item for all $\theta\in \Sigma$, $\{i\in I:\theta\in \Phi_i\}\in\mathcal U$,
\end{enumerate}
then $\Sigma$ is satisfiable in $\prod_I\mathfrak A_i/\mathcal U$.
\end{lemma}

\bproof
$p(\Theta)$ is just the set of all $i\in I$ for which $\Theta$ is satisfiable in $\mathfrak A_i$. Hence by if $i\in p(\Phi_i)$, then there is $a_i\in A_i$ such that $\mathfrak A_i\vDash\bigwedge\Phi_i(a_i)$. By (i), there is such an $a_i\in A_i$ for every $i\in I$. Thus by (ii),
\[\{i\in I: \mathfrak A_i\vDash \theta(a_i)\}\supseteq\{i\in I:\theta\in \Phi_i\}\in\mathcal U\qquad\text{for every }\theta\in \Sigma,\]
 so that 
by \L os' Theorem, 
$\prod_I\mathfrak A_i/\mathcal U\vDash \theta[(a_i)_{i\in I}/\mathcal U]$ for every $\theta\in \Sigma$. Thus the element $(a_i)_{i\in I}/\mathcal U$ satisfies $\Sigma$ in $\prod_I\mathfrak A_i/\mathcal U$.
\eproof

We say that a sequence $(\Phi_i)_{i\in I}$ of members of $\mathcal P^{<\omega}(\Sigma)$ {\em supports} $p$ if it satisfies (i), (ii) of Lemma \ref{lemma_support_saturation}. Thus if there is a sequence supporting the relation $p$ defined in $(\star)$ above, then $\Sigma$ is satisfiable in $\prod_{i\in I}\mathfrak A_i/\mathcal U$. 

Now the above observations do not in any way depend on the linguistic aspects of $\mathcal L, \Sigma$. We can therefore fruitfully move to a slightly more abstract realm:
\begin{definition}\rm\begin{enumerate}[(a)]\item
Let $X$ be a set. An {\em order-reversal} is a map $p:\mathcal P^{<\omega}(X)\to\mathcal U$ such that $\Theta\subseteq\Theta'\Rightarrow p(\Theta)\supseteq p(\Theta')$, or equivalently \[p(\Theta\cup\Theta') \subseteq p(\Theta)\cap p(\Theta').\] (Recall that $\mathcal P^{<\omega}(X)$ denotes the family of finite subsets of $X$.)
\item A sequence $(\Phi_i)_{i\in I}$ of members of $\mathcal P^{<\omega}(X)$ {\em supports} the order-reversal $p:\mathcal P^{<\omega}(X)\to\mathcal U$ if and only if \begin{enumerate}[(i)]\item $i\in p(\Phi_i)$ for all $i\in I$.
\item For all $\theta\in X$, $\{i\in I:\theta\in \Phi_i\}\in\mathcal U$, or equivalently, for all $\Theta\in\mathcal P^{<\omega}(X)$, $\{i\in I:\Theta\subseteq \Phi_i\}\in\mathcal U$.
\end{enumerate}
\item An order-reversal  $p$ is {\em anti-additive} if also
\[p(\Theta\cup\Theta')=p(\Theta)\cap p(\Theta').\]
\item An order-reversal $p$ is {\em locally finite} if for all $i\in I$
\[\sup\{|\Theta|:\Theta\in\mathcal P^{<\omega}(X), i\in p(\Theta)\}<\infty.\] This means that for all $i\in I$ there is $N_i\in\mathbb N$ such that $i\not\in\Theta$ whenever $|\Theta|>N_i$.
\end{enumerate}
\endbox
\end{definition}
%

Observe that the order-reversal defined in $(\star)$ need not be either anti-additive or locally finite. However, if an abstract order-reversal $p$ has a support $\Phi:=(\Phi_i)_{i\in I}$, then we can define another order-reversal $p_\Phi$ which has those properties:

\begin{lemma} Suppose that $p:\mathcal P^{<\omega}(X)\to\mathcal U$ is an order-reversal possessing a support $\Phi:=(\Phi_i)_{i\in I}$. Define a map \[p_\Phi:\mathcal P^{<\omega}(X)\to \mathcal U:\Theta\mapsto\{i\in I: \Theta\subseteq \Phi_i\}.\]
Then $p_\Phi$ is a anti-additive locally finite order-reversal with the properties that (i)	$p_\Phi\leq p$, and (ii) $\Phi$ is a support for $p_\Phi$ also.
\end{lemma}

\bproof
Clearly if $\Theta\subseteq\Theta'$ belong to $\mathcal P^{<\omega}(X)$, then $p_\Phi(\Theta)=\{i:\Theta\subseteq\Phi_i\}\supseteq\{i: \Theta'\subseteq\Phi_i\}=p_\Phi(\Theta')$. Hence $p_\Phi$ is an order-reversal. Moreover, since $\Phi_i\subseteq\Phi_i$, we have $i\in p_\Phi(\Phi_i)$ for all $i\in I$, and hence $p_\Phi$ is supported by $(\Phi_i)_{i\in I}$.

Note also that if $i\in p_\Phi(\Theta)$, then $\Theta\subseteq\Phi_i$, and hence $p(\Theta)\supseteq p(\Phi_i)$. Since also $i\in p(\Phi_i)$, we see that $i\in p(\Theta)$ whenever $i\in p_\Phi(\Theta)$, i.e. $p_\Phi(\Theta)\subseteq p(\Theta)$ for all $\Theta\in\mathcal P^{<\omega}(X)$.

It remains to show that $p_\Phi$ is anti-additive and locally finite. To prove that it is anti-additive, observe that
\[p_\Phi(\Theta\cup\Theta') = \{i:\Theta\cup\Theta'\subseteq \Phi_i\}=\{i:\Theta\subseteq \Phi_i\}\cap\{i:\Theta'\subseteq \Phi_i\}=p_\Phi(\Theta)\cap p_\Phi(\Theta').\]
Next, to prove that $p_\Phi$ is locally finite, observe that if $i\in p_\Phi(\Theta)$, then $\Theta\subseteq\Phi_i$. Hence $\sup\{|\Theta|:i\in p_\Phi(\Theta)\}=|\Phi_i|<\infty$.
\eproof

\begin{lemma} \label{lemma_antiadditive_locallyfinite_support} Suppose that $p:\mathcal P^{<\omega}(X)\to\mathcal U$ is a anti-additive and locally finite order-reversal. Then $p$ has a support $\Phi$ such that $p=p_\Phi$.
\end{lemma}
\bproof
Since $p$ is anti-additive, the set $\{\Theta\in\mathcal P^{<\omega}(X): i\in p(\Theta)\}$ is closed under finite unions, for each $i\in I$. Thus, as $p$ is locally finite, the set $\{\Theta\in\mathcal P^{<\omega}(X): i\in p(\Theta)\}$ has a maximum element, namely $\Phi_i:=\bigcup\{\Theta: i\in p(\Theta)\}$. Then certainly $i\in p(\Phi_i)$. Now observe that
\[i\in p(\Theta)\Rightarrow \Theta\subseteq\Phi_i\Rightarrow p(\Theta)\supseteq p(\Phi_i)\Rightarrow i\in p(\Theta),\] i.e. all the above are equivalent. In particular, for any $\Theta\in\mathcal  P^{<\omega}(X)$ we have\[\{i:\Theta\subseteq\Phi_i\}=\{i:i\in p(\Theta)\}=p(\Theta)\in\mathcal U,\] and hence $(\Phi_i)_{i\in I}$ is a support of $p$.

Finally,
\[p_\Phi(\Theta):=\{i: \Theta\subseteq \Phi_i\}=p(\Theta),\] i.e. $p_\Phi=p$.
\eproof

We now seek conditions that will ensure that any order-reversal on any set $X$ of cardinality $<\kappa$ has a support. Observe that if $p, p'$ are order-reversals so that $p'\leq p$, then $p'$ is locally finite if $p$ is, since $\{\Theta: i\in p'(\Theta)\}\subseteq\{\Theta:i\in p(\Theta)\}$ for all $i\in I$. Furthermore, if $\Phi$ is a support of $p'$, then it is a support of $p$ also, since $I=\{i:i\in p'(\Phi_i)\}\subseteq \{i:i\in p(\Phi_i)\}$.  Hence {\em if } we can define for every order-reversal $p$ two smaller order -eversals $Lp\leq p, Cp\leq p$ so that $Lp$ is locally finite, and $Cp$ is anti-additive, {\em then} $CLp\leq Lp\leq p$ is anti-additive and locally finite. Hence $CLp$ has a support, and this will be a support of $p$ also. 

According to Lemma \ref{lemma_support_saturation}, if the ultrafilter $\mathcal U$ has the property that any order-reversal $p:\mathcal P^{<\omega}(X)\to\mathcal U$ (where $|X|<\kappa$) has a support, then any ultraproduct modulo $\mathcal U$ interpreting a language with $|\mathcal L|<\kappa$ is $\kappa$-saturated.

So our aim is to find conditions on an ultrafilter which guarantee the existence of operators $C, L$.

We can easily deal with the operator $L$:

\begin{lemma}\label{lemma_locally_finite}
If $\mathcal U$ is a countably incomplete ultrafilter on a set $I$, then $L$ exists, i.e. for every order-reversal $p:\mathcal P^{<\omega}(X)\to\mathcal U$, there is a locally finite order-reversal $q\leq p$.
\end{lemma}
\bproof
Since $\mathcal U$ is countably incomplete, there exists a decreasing sequence $(I_n)_{n\in\mathbb N}$ of members of $\mathcal U$ such that $\bigcap_nI_n=\varnothing$. Without loss of generality, we may take $I_0=I$.

 For each $i\in I$, define $N(i):=\min\{n: i\not\in I_n\}$. As the $I_n$ form a decreasing sequence, we have $I_n=\{i: N(i)>n\}$.

Now suppose that we have an order-reversal $p:\mathcal P^{<\omega}(X)\to\mathcal U$. Define $Lp$ by
\[Lp(\Theta):=p(\Theta)\cap I_{|\Theta|},\qquad \text{i.e.}\] \[i\in (Lp)(\Theta) \quad\text{iff}\quad i\in p(\Theta) \land |\Theta|< N(i).\]
It is then easy to see that $Lp$ is an order-reversal and that $Lp\leq p$. Furthermore, $\sup\{|\Theta|: i\in (Lp)(\Theta)\}< N(i)<\infty$ for each $i\in I$, so $Lp$ is locally finite.
\eproof

%

To tackle the existence of $C$, we introduce the following definition:
\begin{definition}\rm
Let $\kappa$ be a cardinal. An ultrafilter $\mathcal U$ is $\kappa$-good if for every cardinal $\alpha<\kappa$ and every  order-reversal $p: \mathcal P^{<\omega}(\alpha)\to\mathcal U$,  $Cp$ exists, i.e. there is a anti-additive order-reversal $q: \mathcal P^{<\omega}(\alpha)\to\mathcal U$ such that $q\leq p$.\endbox
\end{definition}

To recapitulate:
\begin{itemize}\item Suppose that $\mathcal U$ is a countably incomplete $\kappa$-good ultrafilter on a set $I$, and that $(\mathfrak A_i)_{i\in I}$ is a family of models interpreting a language $\mathcal L$ of cardinality $<\kappa$.
\item Let $X\subseteq\prod_I\mathfrak A_i/\mathcal U$ be such that $|X|<\kappa$, and suppose that $\Sigma(x)$ is a family of formulas of $\mathcal L_X$ that is finitely satisfiable in $\prod_I\mathfrak A_i/\mathcal U$.
\item Then $|\Sigma(x)|<\kappa$.
\item If $\Theta\in\mathcal P^{<\omega}(\Sigma)$, then $\{i\in I:\mathfrak A_i\vDash\exists x \bigwedge \Theta\}\in\mathcal U$, since $\Theta$ is satisfiable in $\prod_I\mathfrak A_i/\mathcal U$. Thus the map
\[p:\mathcal P^{<\omega}(\Sigma)\to\mathcal U: \Theta\mapsto \{i\in I:\mathfrak A_i\vDash\exists x \bigwedge \Theta\}\] is an order-reversal.
\item By Lemma  \ref{lemma_support_saturation}, the set $\Sigma$ is satisfiable in $\prod_I\mathfrak A_i/\mathcal U$ when $p$ has a support.
\item Lemma \ref{lemma_antiadditive_locallyfinite_support} shows that every locally finite consistent order-reversal has a support. 
\item Since $\mathcal U$ is countably incomplete, Lemma \ref{lemma_locally_finite} shows that there is a locally finite $p'\leq p$.
\item The fact that $\mathcal U$ is $\kappa$-good then yields the existence of a $p''\leq p'$ which is consistent. That $p''$ will therefore also be locally finite, because $p'$ is. Hence $p''$ has a support $\Phi$. Then as $p''\leq p$, $\Phi$ will also be a support for $p$.
\item Hence $p$ has a support, and thus $\Sigma$ is satisfiable in $\prod_I\mathfrak A_i/\mathcal U$.
\item Since $X$, $\Sigma(x)$ were arbitrary, every finitely satisfiable family of formulas $\Sigma(x)$ in the expanded language $\mathcal L_X$ is satisfiable in $\prod_I\mathfrak A_i/\mathcal U$ (for $|X|<\kappa$). Thus $\prod_I\mathfrak A_i/\mathcal U$ is $\kappa$-saturated.
\end{itemize}
\vskip0.5cm

Thus we have shown:
\begin{theorem} \label{thm_saturated_ultrapower} If $\mathcal U$ is a countably incomplete $\kappa$-good ultrafilter, then every ultraproduct modulo $\mathcal U$ interpreting a language of cardinality $<\kappa$ is $\kappa$-saturated.\endbox	
\end{theorem}

It remains to address the existence of good ultrafilters. The following theorem requires some heavy-duty combinatorics, so its proof has been relegated to the appendix.
\begin{theorem} \label{thm_exist_good_ultrafilter} Suppose that $\kappa$ is an infinite cardinal. Then there exists a countably incomplete $\kappa^+$-good ultrafilter over $\kappa$.\endbox
\end{theorem}

\begin{example}\label{example_ultrafilter_omega1_good} \rm For the case $\kappa^+=\omega_1$,  is quite easy to show that there are $\omega_1$--good countably incomplete ultrafilters, for the simple reason that {\em every} ultrafilter is $\omega_1$-good. Indeed, suppose that $\mathcal U$ is an ultrafilter over a set $I$, and that $p:\mathcal P^{<\omega}(\omega)\to \mathcal U$ is an order-reversal. For $s\in\mathcal P^{<\omega}(\omega)$, define\[q(s):=p\{m\in\omega:m\leq \max s\}.\]
Then as $s\subseteq\{m\in\omega:m\leq \max s\}$, we have $q(s)\subseteq p(s)$.  Moreover, since for  $s, t\in\mathcal P^{<\omega}(\omega)$ we have \[\{m:m\leq \max (s\cup t)\}=\left\{\aligned \{m:m\leq \max s\}\quad&\text{if }\max s\geq \max t,\\\{m:m\leq\max t\}\quad&\text{else },\endaligned\right.\]
it follows that $q(s\cup t)=q(s)\cap q(t)$, i.e. that $q$ is anti-additive.

Hence Theorem \ref{thm_ultrapower_omega1_saturated} also follows from Theorem \ref{thm_exist_good_ultrafilter}.
\endbox
\end{example}

\subsubsection{Construction of Polysaturated Extensions via Good Ultrapowers}
Recall the ultrapower construction of a nonstandard extension $V(X)\overset{*}{\hookrightarrow}V(^*X)$: We start with a base set $X$ and an ultrafilter $\mathcal U$ over a set $I$, where $I$ is chosen so that the ultrapower $Y:=X^I/\mathcal U$ is another base set. We define relations $=_\mathcal U, \in_\mathcal U$ on $V(X)^I$ by\[f=_{\mathcal U}g\Leftrightarrow \{f=g\}\in\mathcal U,\qquad f\in_{\mathcal U}g\Leftrightarrow \{f\in g\}\in\mathcal U.\]
For $n\in\mathbb N$, we define \[W_n:=\{f\in V(X)^I: f\in_\mathcal U c_{V_n(X)}\},\quad\text{ and then } W:=\bigcup_{n\in\mathbb N} W_n,\]
where for $a\in V(X)$ the map $c_a:I\to V(X)$ is the constant map with value $a$. Observe that if $a\in V(X)$, then $a\in V_n(X)$ for some $n\in\mathbb N$, and hence $c_a\in W_n$. Thus there is a natural inclusion
\[\iota:V(X)\hookrightarrow W:a\mapsto c_a.\]
By induction, we construct a map $\cdot/\mathcal U:W\to V(Y)$, as follows: For $f\in W_0$, define \[f/\mathcal U:=\{g\in X^I: f=_\mathcal U g\},\]
and for $f\in W-W_0$, define\[f/\mathcal U:=\{g/\mathcal U: g\in W\land g\in_\mathcal U f\}.\]
The map $\cdot/\mathcal U$ has the property that $f/\mathcal U\in V_n(Y)$ whenever $f\in W_n$. 

 Then the $*$-map, defined as the composition $V(X)\overset{\iota}{\hookrightarrow}W\overset{\cdot/\mathcal U}{\to}V(Y)$, is a transfer map, with $Y={}^*X$. If the ultrafilter $\mathcal U$ is also countably incomplete, then $*:V(X)\to V({}^*X)$ is a nonstandard framework, in that $\{{}^*a:a\in A\}$ is a proper subset of ${}^*a$ whenever $A\in V(X)$ is an infinite set.
 \vskip0.3cm
Observe that the internal sets are precisely the sets of the form $f/\mathcal U$, for $f\in W$. Indeed, if $f\in W$, then $f\in W_n$ for some $n\in\mathbb N$, and hence $f\in_\mathcal U c_{V_n(X)}$, so that $f/\mathcal U\in{} ^*V_n(X)=c_{V_n(X)}/\mathcal U$, from which it follows that $f/\mathcal U$ is internal. Conversely, if $A\in V(Y)$ is internal, then $A\in{}^*B=c_B/\mathcal U$, where $c_B/\mathcal U:=\{f/\mathcal U: f\in W\land f\in_\mathcal U c_B\}$, from which it follows that $A=f/\mathcal U$ for some $f\in W$. Thus \[{}^*V(X)=\{f/\mathcal U:f\in W\}.\]

\begin{theorem}\label{thm_transfer_map_good_ultrafilter_saturated} Let $\kappa$ be a cardinal, and 
suppose that a nonstandard framework $V(X)\overset{*}{\hookrightarrow}V({}^*X)$ is obtained as an ultrapower construction via a $\kappa$-good countably incomplete ultrafilter $\mathcal U$ over a set $I$. Then $V({}^*X)$ is $\kappa$-saturated for then language $\mathcal L_{{}^*V(X)}$, i.e. if $\mathcal A$ is a family of internal sets with the f.i.p. such that $|\mathcal A|<\kappa$, then $\bigcap \mathcal A \neq \varnothing$.
\end{theorem}

\bproof
Let $\Gamma<\kappa$ be an ordinal, and suppose that $\mathcal A:=\{A_\gamma: \gamma<\Gamma\}$ is a family of internal sets with the f.i.p. We must show that $\bigcap \mathcal A\neq\varnothing$. 
Without loss of generality, by replacing $A_\gamma$ by $A_\gamma\cap A_0$, we may assume that $A_\gamma\subseteq A_0$ for all $\gamma<\Gamma$: This affects neither the f.i.p. nor the value of $\bigcap\mathcal A$.

As each $A_\gamma$ is internal, there is a function $a_\gamma\in W$ such that $A_\gamma=a_\gamma/\mathcal U$. 

As $\mathcal U$ is countably incomplete, there is a sequence $I=I_1\supseteq I_2\supseteq I_3\supseteq\dots$ of members of $\mathcal U$ such that $\bigcap_{n=1}^\infty I_n=\varnothing.$ Now define a map
\[f:\mathcal P^{<\omega}(\Gamma)\to\mathcal U:\Delta\mapsto I_n\cap\Big\{i\in I:\bigcap_{\gamma\in\Delta}a_{\gamma}(i)\neq \varnothing\Big\}.\] (Note that  $\bigcap_{\gamma\in\Delta}A_{\gamma}\neq \varnothing$ for any finite $\Delta\subseteq\Gamma$, by the f.i.p., so the set $\Big\{i\in I:\bigcap_{\gamma\in\Delta}a_{\gamma}(i)\neq \varnothing\Big\}$ belongs to $\mathcal U$.) Observe that $f$ is a reversal. Since $\Gamma<\kappa$ and $\mathcal U$ is $\kappa$-good, there is a strict reversal $g:\mathcal P^{<\omega}(\Gamma)\to\mathcal U$ such that $g\leq f$.

For $i\in I$, define \[\Gamma_i:=\{\gamma\in \Gamma: i\in g(\{\gamma\})\},\] so that $\gamma\in \Gamma_i$ if and only if $i\in g(\{\gamma\})$. We first show that each $\gamma_i$ is a finite set. Fix $i\in I$. Since $\bigcap_{n=1}^\infty I_n=\varnothing$, there is $n\in \mathbb N$ such that $i\not\in I_n$. We claim that $|\Gamma_i|<n$. Indeed, if $\gamma_1,\dots,\gamma_n$ are distinct elements of $\Gamma_i$, then \[i\in\bigcap_{m=1}^n g(\{\gamma_m
\})=g(\{\gamma_1,\dots,\gamma_n\})\subseteq f(\{\gamma_1,\dots,\gamma_n\})\subseteq I_n,\] which is impossible, as $i\not\in I_n$.
Thus each $\Gamma_i\in\mathcal  P^{<\omega}(\Gamma)$, so that $g(\Gamma_i)$ is defined.

Now observe that $g(\Gamma_i)=g(\bigcup_{\gamma\in \Gamma_i}\{\gamma\})=\bigcap_{\gamma\in \Gamma_i}g(\{\gamma\})$. As $i\in g(\{\gamma\})$ whenever $\gamma\in \Gamma_i$, we see that $i\in g(\Gamma_i)$. As $f(\Gamma_i)\supseteq g(\Gamma_i)$, we have $i\in f(\Gamma_i)$, so that $\bigcap_{\gamma\in \Gamma_i}a_\gamma(i)\neq \varnothing$, by definition of $f$.

Now choose a map $x\in W$ so that $x(i)\in\bigcap_{\gamma\in \Gamma_i}a_\gamma(i)$. (Recall that $A_\gamma\subseteq A_0$, and that $A_0=a_0/\mathcal U$ for some $a_0\in W$. Then $a_0\in W_n$ for some $n$, and hence we may take $x\in W_{n-1}\subseteq W$.) 

Now if $\gamma\in \Gamma_i$, then $x(i)\in a_{\gamma}(i)$, by definition of $x$. Hence if $\gamma\in \Gamma$, then
\[\{i\in I: x(i)\in a_{\gamma}(i)\}\supseteq \{i\in I: \gamma\in \Gamma_i\}=g(\{\gamma\})\in\mathcal U,\]
from which we see that $x/\mathcal U\in a_\gamma/\mathcal U=A_\gamma$. As $\gamma\in \Gamma$ was arbitrary, we have $x/\mathcal U\in\bigcap\mathcal A$.
\eproof

\subsection{Existence of Polysaturated Extensions via Ultralimits}\label{section_polysaturated_ultralimits}

\subsubsection{Limits of Chains of Superstructures}
Suppose that $\lambda $ is a limit ordinal, and that $\{X^\alpha: \alpha<\lambda\}$ is a collection of base sets. Suppose further that, for $\alpha\leq\beta<\lambda$, we have a chain of superstructures $(V(X^\alpha)))_{\alpha<\lambda}$ linked by bounded elementary embeddings $V(X^\alpha)\overset{\iota_{\alpha\beta}}{\to}V(X^\beta)$ with the following properties:
\begin{enumerate}[(i)]\item $\iota_{\alpha\alpha}=\text{id}_{V(X^\alpha)}$.
\item If $\alpha\leq\beta\leq \gamma<\lambda$, then $\iota_{\beta\gamma}\circ\iota_{\alpha\beta}=\iota_{\alpha\gamma}$.
\item If $\alpha\leq\beta<\lambda$, then $\iota_{\alpha\beta}(X^\alpha)=X^\beta$.
\end{enumerate}

We now want to construct a limit model $V(X^\lambda)$ and bounded elementary embeddings \\$V(X^\alpha)\overset{\iota_{\alpha\lambda}}{\to}V(X^\lambda)$ so that properties (i)-(iii) hold for $\alpha\leq\beta\leq\lambda$, i.e. for $\lambda$ as well.

\begin{lemma} The bounded elementary embeddings $\iota_{\alpha\beta}$ are rank-preserving, i.e. if $\alpha\leq \beta$ and $n\in\mathbb N$, then  $x\in V_n(X^\alpha)$ if and only if $\iota_{\alpha\beta}(x)\in V_n(X^\beta)$.
\end{lemma}

\bproof By Lemma \ref{lemma_define_simple_notions} there is a bounded formula $\varphi_{6,n}(X,x)$ such that $V(X)\vDash \varphi_{6,n}(X,x)$ if and only if $x\in V_n(X)$.
Thus $x\in V_n(X^\alpha)$ if and only if $V(X^\alpha)\vDash \varphi_{6,n}(X^\alpha,x)$ if and only if $V(X^\beta)\vDash \varphi_{6,n}(\iota_{\alpha\beta}(X^\alpha), \iota_{\alpha\beta}(x))$ if and only if $\iota_{\alpha\beta}(x)\in V_n(X^\beta)$, using the fact that $\iota_{\alpha\beta}(X^\alpha)=X^\beta$.\eproof

We now proceed to make a limit model out of the models $V(X^\alpha)$. For $n\in\mathbb N$, define
\[P_n:=\{(a,\alpha):\alpha<\lambda\land a\in V_n(X^\alpha)\},\qquad P:=\bigcup_{n\in\mathbb N}P_n.\]
Observe that $P_0\subseteq P_1\subseteq P_2\subseteq \dots$.

Define binary relations $\sim, E$ on $P$ as follows:
\[(a,\alpha)\sim (b,\beta) \Leftrightarrow\exists\gamma\geq\alpha,\beta\;\Big(\iota_{\alpha\gamma}(a)=\iota_{\beta\gamma}(b)\Big),\qquad (a,\alpha) E (b,\beta) \Leftrightarrow\exists\gamma\geq\alpha,\beta\;\Big(\iota_{\alpha\gamma}(a)\in\iota_{\beta\gamma}(b)\Big).\]

\begin{lemma}  \begin{enumerate}[(a)]\item
$\sim$ is an equivalence relation.\item have  that \[(a,\alpha)\sim (b,\beta) \Leftrightarrow\forall\delta\geq\alpha,\beta\;\Big(\iota_{\alpha\delta}(a)=\iota_{\beta\delta}(b)\Big),\] and that \[(a,\alpha) E (b,\beta) \Leftrightarrow\forall\delta\geq\alpha,\beta\;\Big(\iota_{\alpha\delta}(a)\in\iota_{\beta\delta}(b)\Big).\]
\item If $\alpha\leq \beta$, then $(a,\alpha)\sim(b,\beta)$ if and only if $\iota_{\alpha\beta}(a)=b$, and $(a,\alpha)E(b,\beta)$ if and only if $\iota_{\alpha\beta}(a)\in b$.
\item $(a,\alpha)\sim(b,\beta)$ if and only if $\forall (c,\gamma)\in P\;\Big((c,\gamma)E(a,\alpha)\leftrightarrow (c,\gamma)E(b,\beta)\Big)$.
\end{enumerate}
\end{lemma}

\bproof (a) It is clear that $\sim$ is reflexive and symmetric. If $(a,\alpha)\sim(b,\beta)\sim(c,\gamma)$, then there are $\eta\geq\alpha,\beta$ and $\xi\geq \beta,\gamma$ such that $\iota_{\alpha\eta}(a)=\iota_{\beta\eta}(b)$ and $\iota_{\beta\xi}(b)=\iota_{\gamma\xi}(c)$. Let $\delta\geq \eta,\xi$. Then \[\iota_{\alpha\delta}(a)=\iota_{\eta\delta}\circ\iota_{\alpha\eta}(a)=\iota_{\eta\delta}\circ\iota_{\beta\eta}(b)=\iota_{\beta\delta}(b)=\iota_{\xi\delta}\circ\iota_{\beta\xi}(b)=\iota_{\xi\delta}\circ\iota_{\gamma\xi}(c)=\iota_{\gamma\delta}(c),\] which shows that $(a,\alpha)\sim (c,\gamma)$, establishing transitivity.

(b) Suppose that $(a,\alpha)\sim (b,\beta)$ and that $\gamma\geq \alpha,\beta$ is such that $\iota_{\alpha\gamma}(a)=\iota_{\beta\gamma}(b)$. Let $\delta\geq\alpha,\beta$.
If $\delta\geq\gamma$, then \[\iota_{\alpha\delta}(a)=\iota_{\gamma\delta}\circ\iota_{\alpha\gamma}(a)=\iota_{\gamma\delta}\circ\iota_{\beta\gamma}(b)=\iota_{\beta\delta}(b).\]
If, on the other hand, $\alpha,\beta\leq\delta\leq \gamma$, then
\[i_{\delta\gamma}\circ\iota_{\alpha\delta}(a)=\iota_{\alpha\gamma}(a)=\iota_{\beta\gamma}(b)=\iota_{\delta\gamma}\circ\iota_{\beta\delta}(b).\] Then since $\iota_{\delta\gamma}$ is one-to-one, we conclude that $\iota_{\alpha\delta}(a) = \iota_{\beta\delta}(b)$ in this case also.\newline The proof for $E$ is similar.

(c) follows directly from (b): Since $\beta\geq\alpha,\beta$, it follows that $\iota_{\alpha\beta}(a)=\iota_{\beta\beta}(b)=b$.\newline The proof for $E$ is similar.

(d) Suppose that $(a,\alpha)\sim (b,\beta)$ and that $(c,\gamma)E(a,\alpha)$. Applying (b), we see that if $\delta\geq\alpha,\beta,\gamma$, then $\iota_{\gamma\delta}(c)\in\iota_{\alpha\delta}(a)=\iota_{\beta\delta}(b)$, so that $(c,\gamma)E(b,\beta)$. By symmetry, we see that if $(a,\alpha)\sim (b,\beta)$, then $(c,\gamma)E(a,\alpha)$ if and only if $(c,\gamma)E(b,\beta)$.\newline Conversely, suppose that for all $(c,\gamma)\in P$ we have  $(c,\gamma)E(a,\alpha)$ if and only if $(c,\gamma)E(b,\beta)$. Let $\gamma\geq\alpha,\beta$. If $(a,\alpha)\not\sim(b,\beta)$, then $\iota_{\alpha\gamma}(a)\neq \iota_{\beta\gamma}(b)$, so there is $c\in (\iota_{\alpha\gamma}(a)- \iota_{\beta\gamma}(b))\cup(\iota_{\beta\gamma}(b)-\iota_{\alpha\gamma}(a))$. Without loss of generality, suppose that $c\in (\iota_{\alpha\gamma}(a)- \iota_{\beta\gamma}(b))$. Then by (c),
\[c\in\iota_{\alpha\gamma}(a)\Rightarrow (c,\gamma)E(a,\alpha)\Rightarrow (c,\gamma)E(b,\beta)\Rightarrow c\in\iota_{\beta\gamma}(b),\] contradiction. Hence $(a,\alpha)\sim (b,\beta)$.
\eproof

\begin{lemma} \label{lemma_sim_E_rank_P} \begin{enumerate}[(a)]\item If $(a,\alpha)\in P_n$ and $(b,\beta)\sim(a,\alpha)$, then $(b,\beta)\in P_n$.
\item If $(a,\alpha)\in P_n$ and $(b,\beta)E(a,\alpha)$, then $(b,\beta)\in P_{n-1}$.
\end{enumerate}\end{lemma}

\bproof
(a) Suppose that $(a,\alpha)\in P_n$ and $(b,\beta)\sim(a,\alpha)$. If $\gamma\geq\alpha,\beta$, then $\iota_{\alpha\gamma}(a)=\iota_{\beta\gamma}(b)$.
But by the rank-preserving property of the $\iota_{\alpha\beta}$, we have\[(a,\alpha)\in P_n\Rightarrow a\in V_n(X^\alpha)\Rightarrow \iota_{\alpha\gamma}(a)\in V(X^\gamma)\Rightarrow \iota_{\beta\gamma}(b)\in V(X^\gamma)\Rightarrow b\in V_n(X^\beta)\Rightarrow (b,\beta)\in P_n.\]
(b) First observe that if $(a,\alpha)\in P_0$, and $(b,\beta)E(a,\alpha)$, then $b\in \iota_{0\beta}(a)\in V_0(X^\beta)=X^\beta$, i.e. $\iota_{0\beta}(a)$ is a member of the base set $X^\beta$ which has an element $b\in V(X^\beta)$ --- contradicting the definition of {\em base set}. Hence if $(a,\alpha)\in P_0$ there can be no $(b,\beta)$ such that $(b,\beta)E(a,\alpha)$.  Thus if $(b,\beta)E(a,\alpha)$, then $n\geq 1$. Then if $\gamma\geq\alpha,\beta$, it follows by the rank-preserving properties of the $\iota_{\alpha\beta}$ that \[\iota_{\beta\gamma}(b)\in\iota_{\alpha\gamma}(a)\in V_n(X^\gamma)\Rightarrow\iota_{\beta\gamma}(b)\in V_{n-1}(X^\gamma)\Rightarrow b\in V_{n-1}(X^\beta)\Rightarrow (b,\beta)\in P_{n-1}.\]

\eproof

Now define a sequence of sets $(W_n)_{n\in\mathbb N}$ by induction, as follows: First, if $(a,\alpha)\in P_0$, let $[a,\alpha]:=\{(b,\beta): (b,\beta)\sim(a,\alpha)\}$. By Lemma \ref{lemma_sim_E_rank_P}, $[a,\alpha]\subseteq P_0$. Then let
\[X^\lambda:=W_0:=\{[a,\alpha]:(a,\alpha)\in P_0\}.\]
If necessary, modify $X^\lambda$ so that it is a base set.

Now assume that $[b,\beta]$ and $W_k$ have already been defined for all $k<n$ and $(b,\beta)\in P_k$. Now define $[a,\alpha]$ for $(a,\alpha)\in P_n$, and define $W_n$ by\[[a,\alpha]:=\{[b,\beta]: (b,\beta)E(a,\alpha)\},\qquad W_n:=\{[a,\alpha]:(a,\alpha)\in P_n\}.\] (Observe that if $(b,\beta)E(a,\alpha)$, then $(b,\beta)\in P_{n-1}$ by Lemma \ref{lemma_sim_E_rank_P}, so that $[b,\beta]$ has already been defined.)
Then define \[W:=\bigcup_{n\in\mathbb N}W_n.\]
The elements of $W_0=X^\lambda$ act as atoms. The elements of $W-W_0$ act as sets.
Observe that $W$ is transitive over sets: If $x\in [a,\alpha]$, where $[a,\alpha]\in W-W_0$ is a set, then $x=[b,\beta]$ for some $(b,\beta)E(a,\alpha)$ and hence $x\in W$ also.  Note also that $W_n\subseteq V_n(X^\lambda)$, for all $n\in\mathbb N$.

\begin{lemma}\label{properties_[a,alpha]}\begin{enumerate}[(a)]\item If $\alpha\leq \beta<\lambda$ and $a\in V(X^\alpha)$, then $[a,\alpha]=[\iota_{\alpha\beta}(a),\beta]$.
\item If $\alpha,\beta<\lambda$, then $[b,\beta]\in [X^\alpha,\alpha]$ if and only if $b\in X^\beta$.
\end{enumerate}
\end{lemma}
\bproof
(a) If $(a,\alpha)\in P_0$, then \begin{align*}[\iota_{\alpha\beta}(a),\beta]&=\{(c,\gamma):(c,\gamma)\sim(\iota_{\alpha\beta}(a),\beta)\} \\&=\{(c,\gamma):\exists\delta\geq\beta,\gamma\;\Big(\iota_{\gamma\delta}(c)=\iota_{\beta\delta}(\iota_{\alpha\beta}(a))=\iota_{\alpha\delta}(a)\Big)\} 
\\&=\{(c,\gamma): (c,\gamma)\sim(a,\alpha)\}\\&=[a,\alpha]\end{align*} 
Similarly, if $(a,\alpha)\in P-P_0$, then
\begin{align*}[\iota_{\alpha\beta}(a),\beta]&=\{[c,\gamma]:(c,\gamma)E(\iota_{\alpha\beta}(a),\beta)\}\\&=\{[c,\gamma]:\exists\delta\geq\beta,\gamma\;\Big(\iota_{\gamma\delta}(c)\in\iota_{\beta\delta}(\iota_{\alpha\beta}(a))=\iota_{\alpha\delta}(a)\Big)\} ,
\\&=\{[c,\gamma]: (c,\gamma)E(a,\alpha)\},\\&=[a,\alpha].\end{align*}

(b) We have \begin{align*}[b,\beta]\in [X^\alpha,\alpha]&\Leftrightarrow (b,\beta)E(X^\alpha,\alpha),\\
&\Leftrightarrow \iota_{\beta\gamma}(b)\in\iota_{\alpha\gamma}(X^\alpha)\quad\text{for $\gamma\geq \alpha,\beta$},\\
&\Leftrightarrow \iota_{\beta\gamma}(b)\in X^\gamma,\\
&\Leftrightarrow  b\in X^\beta,\quad\text{by the rank-preserving property of $\iota_{\beta\gamma}$.}
\end{align*}
\eproof

For $\alpha\leq\lambda$, define maps $V(X^\alpha)\overset{\iota_{\alpha\lambda}}{\to}V(X^\lambda)$ as follows:
\[\iota_{\alpha\lambda}(a)=[a,\alpha] \quad\text{for $\alpha<\lambda$},\qquad \iota_{\lambda\lambda}=\text{id}_{V(X^\lambda)}.\]
We claim that each $\iota_{\alpha\lambda}$ is a bounded elementary embedding with the desired properties:
\begin{enumerate}[(i)]\item $\iota_{\lambda\lambda}=\text{id}_{V(X^\lambda)}$.
\item If $\alpha\leq\beta\leq\lambda$, then $\iota_{\beta\lambda}\circ\iota_{\alpha\beta}=\iota_{\alpha\lambda}$.
\item If $\alpha\leq\lambda$, then $\iota_{\alpha\lambda}(X^\alpha)=X^\lambda$.
\end{enumerate}
(i) holds by definition. Observe that if $\alpha\leq\beta\leq\lambda$,  and $(a,\alpha)\in P$, then 
\[\iota_{\beta\lambda}\circ\iota_{\alpha\beta}(a)=[\iota_{\alpha\beta}(a),\beta]=[a,\alpha]=\iota_{\alpha\lambda}(a), \] using Lemma \ref{properties_[a,alpha]}(a).
This proves (ii).

To prove (iii), note that \begin{align*}\iota_{\alpha\lambda}(X^\alpha)&=[X^\alpha,\alpha],\\&=\{[b,\beta]:[b,\beta]\in [X^\alpha,\alpha]\},\\&=\{[b,\beta]:b\in X^\beta=V_0(X^\beta)\},\\&=\{[b,\beta]:(b,\beta)\in P_0\},\\&=X^\lambda,
\end{align*} using Lemma \ref{properties_[a,alpha]}(b).

It remains to show that each $\iota_{\alpha\lambda}$ is a bounded elementary embedding. This follows by induction on the complexity of formulas. For atomic formulas, we have 
\[V(X^\alpha)\vDash a\in b\quad\Leftrightarrow\quad [a,\alpha]\in[b,\alpha] \quad\Leftrightarrow\quad V(X^\lambda)\vDash \iota_{\alpha\lambda}(a)\in\iota_{\alpha\lambda}(b).\]
Since $\iota_{\alpha\beta}$ is one-to-one, also $a=b$ if and only if $\iota_{\alpha\lambda}(a)=\iota_{\alpha\lambda}(b)$.

The propositional connectives $\land,\lnot$ are dealt with very easily.

 Suppose now that $\varphi(x_1,\dots,x_n)$ is of the form $\exists y\in x_1\;\psi(y,x_1,\dots, x_n)$.
If $V(X^\alpha)\vDash \exists y\in a_1\;\psi(y,a_1,\dots, a_n)$, then there is $b\in a_1$ such that $V(X^\alpha)\vDash\psi(b,a_1,\dots, a_n)$. By induction hypothesis, 
$V(X^\lambda)\vDash \psi(\iota_{\alpha\lambda}(b),\iota_{\alpha\lambda}(a_1),\dots, \iota_{\alpha\lambda}(a_n))$, where $\iota_{\alpha\lambda}(b)\in\iota_{\alpha\lambda}(a_1)$, and hence $V(X^\lambda)\vDash \exists y\in\iota_{\alpha\lambda}(a_1)\;\psi(y, \iota_{\alpha\lambda}(a_1),\dots, \iota_{\alpha\lambda}(a_n))$.

Conversely, suppose that $V(X^\lambda)\vDash \exists y\in\iota_{\alpha\lambda}(a_1)\;\psi(y, \iota_{\alpha\lambda}(a_1),\dots, \iota_{\alpha\lambda}(a_n))$.
Then there is $[b,\beta]\in [a_1,\alpha]$ such that $V(X^\lambda)\vDash \psi([b,\beta], [a_1,\alpha],\dots, [a_n,\alpha])$. Let $\gamma\geq\alpha,\beta$.
Then $[b,\beta]=[\iota_{\beta\gamma}(b),\gamma]=\iota_{\gamma\lambda}(\iota_{\beta\gamma}(b))$ and $[a_i,\alpha]=[\iota_{\alpha\gamma}(a_i),\gamma]=\iota_{\gamma\lambda}(\iota_{\alpha\gamma}(a_i))$, with $\iota_{\beta\gamma}(b)\in \iota_{\alpha\gamma}(a_1)$. Thus 
\[V(X^\lambda)\vDash \psi\Big(\iota_{\gamma\lambda}(\iota_{\beta\gamma}(b)), \iota_{\gamma\lambda}(\iota_{\alpha\gamma}(a_1)),\dots,\iota_{\gamma\lambda}(\iota_{\alpha\gamma}(a_n))\Big), \] where $\iota_{\gamma\lambda}(\iota_{\beta\gamma}(b))\in \iota_{\gamma\lambda}(\iota_{\alpha\gamma}(a_1))$. By induction hypothesis, we obtain
\[V(X^\gamma)\vDash \psi(\iota_{\beta\gamma}(b), \iota_{\alpha\gamma}(a_1),\dots,\iota_{\alpha\gamma}(a_n)), \] where $\iota_{\beta\gamma}(b)\in \iota_{\alpha\gamma}(a_1)$, from which we obtain
\[V(X^\gamma)\vDash \exists y\in  \iota_{\alpha\gamma}(a_1)\;\psi(y, \iota_{\alpha\gamma}(a_1),\dots,\iota_{\alpha\gamma}(a_n)).\]
But as $V(X^\alpha)\overset{\iota_{\alpha\gamma}}{\to}V(X^\gamma)$ is a bounded elementary embedding, it follows also that 
\[V(X^\alpha)\vDash \exists y\in a_1\psi(y,a_1,\dots,a_n).\] This completes the induction, and the proof that each map $V(X^\alpha)\overset{\iota_{\alpha\lambda}}{\to}V(X^\lambda)$ is a bounded elementary embedding.

\subsubsection{Construction of Polysaturated Extensions via Ultralimits}
Recall Theorem \ref{thm_enlargements_exist}, which states that for every  superstructure $V(X)$ there is a set $I$ and an ultrafilter over $I$ such that the induced  bounded elementary embedding $*:V(X)\to V(^*X)$ is an enlargement, where ${}^*X=X^I/\mathcal U$. The idea behind the proof of the existence of a polysaturated extension is to iterate this construction.

Let $V(X)$ be a superstructure over a base set $X$, and let $\kappa=|V(X)|$ be its cardinality. For $\alpha\leq\beta\leq\kappa^+$, we construct superstructures $V(X^\alpha)$ and
bounded elementary embeddings $*_{\alpha\beta}:V(X^\alpha)\to V(X^\beta)$ such that 
\begin{enumerate}[(i)]\item If $\alpha\leq\kappa^+$, then $*_{\alpha\alpha}=\text{id}_{V(X^\alpha)}$.
\item If $\alpha\leq\beta\leq \gamma\leq\kappa^+$, then $*_{\beta\gamma}\circ*_{\alpha\beta}=*_{\alpha\gamma}$.
\item If $\alpha\leq\beta\leq\kappa^+$, then $*_{\alpha\beta}(X^\alpha)=X^\beta$.
\end{enumerate}
We proceed by transfinite induction:

We define $V(X^0)=V(X)$, and $*_{00}=\text{id}_{V(X)}$.

Suppose now that superstructures $V(X^\alpha)$ and
bounded elementary embeddings $*_{\alpha\beta}:V(X^\alpha)\to V(X^\beta)$ have already been constructed for $\alpha\leq\beta<\lambda$, such that (i)-(iii) are satisfied. We now consider two cases:
\newline\underline{Case 1:} $\lambda$ is a limit ordinal. In that case construct $V(X^\lambda)$ and $*_{\alpha\lambda}$ as a limit of the $V(X^\alpha)$ and $ *_{\alpha\beta}$ for $\alpha,\beta<\lambda$. We have just seen that such a construction yields bounded elementary embeddings which preserve properties (i)-(iii).
\newline\underline{Case 2:} $\lambda=\gamma+1$ is a successor ordinal. In that case, let $V(X^\gamma)\overset{*_{\gamma\lambda}}{\to}V(X^\lambda)$ be an enlargement, as provided by Theorem \ref{thm_enlargements_exist}. Observe that $X^\lambda={}^{*_{\gamma\lambda}}X^\gamma$ by construction.
For $\alpha<\gamma$, define $*_{\alpha\lambda}=*_{\gamma\lambda}\circ*_{\alpha\gamma}$, and define $*_{\lambda\lambda}=\text{id}_{V(X^\lambda)}$. It is then straightforward to show  that that properties (i)-(iii) are satisfied.

Let's briefly recall the construction of the final step in the transfinite induction: $V(X^{\kappa^+})$: This is a limit step, so we have\[P_n:=\{(a,\alpha):\alpha<\kappa^+, a\in V_n(X^\alpha)\},\qquad P=\bigcup_nP_n.\]
The binary relations $\sim, E$ are given by
\[(a,\alpha)\sim(b,\beta)\Leftrightarrow\exists\gamma\geq\alpha,\beta\;(*_{\alpha\gamma}(a)=*_{\beta\gamma}(b)),\qquad (a,\alpha)E(b,\beta)\Leftrightarrow\exists\gamma\geq\alpha,\beta\;(*_{\alpha\gamma}(a)\in*_{\beta\gamma}(b)).\]
For $(a,\alpha)\in P_0$, we define $[a,\alpha]:=\{(b,\beta):(b,\beta)\sim (a,\alpha)\}$ and for  $(a,\alpha)\in P-P_0$ we put $[a,\alpha]:=\{[b,\beta]:(b,\beta)E(a,\alpha)$. Then we define
\[X^{\kappa^+}:=W_0:=\{[a,\alpha]:(a,\alpha)\in P_0\},\qquad W_n:=\{[a,\alpha]:(a,\alpha)\in P_n\},\qquad W:=\bigcup_nW_n.\] It then turns out that each $W_n\subseteq V_n(X^{\kappa^+})$, and that $W$ is a submodel of $V(X^{\kappa^+})$ which is transitive over sets. In addition,
\[[a,\alpha]=[*_{\alpha\beta}(a),\beta]\quad \text{ for }\alpha\leq\beta <\kappa^+.\]Finally we define\[*_{\alpha\kappa^+}:V(X^\alpha)\to V(X^{\kappa^+}):a\mapsto[a,\alpha].\]

Let $\mathcal A\subseteq W_n$ be a family of cardinality $\kappa$ with the f.i.p. We will show that $\bigcap\mathcal A\neq\varnothing$. Note that each $A\in\mathcal A$ is of the form $[a,\alpha]$ for some $a$ in some $V_n(X^\alpha)$.  Suppose that
\[\mathcal A=\{A_\xi:\xi<\kappa\}=\{[a_\xi,\alpha_\xi]:\xi<\kappa\} \quad\text{is an enumeration of  }\mathcal A.\]
Let $\beta:=\sup_{\xi<\kappa}\alpha_\xi$. Since $\kappa^+$ is a regular cardinal and each $\alpha_\xi<\kappa^+$, we have that $\beta<\kappa^+$. 
Define $a_\xi'=*_{\alpha_\xi\beta}(a_\xi)\in V_n(X^\beta)$, so that $A_\xi=[a_\xi, \alpha_\xi]=[a'_\xi,\beta]$, and let 
$\mathcal A':=\{a'_\xi:\xi<\kappa\}$
Observe that $*_{\beta\kappa^+}(a'_\xi)=[a'_\xi,\beta]=A_\xi$. As $\mathcal A$ has the f.i.p. and $*_{\beta\kappa^+}$ is a bounded elementary embedding, it follows easily that $\mathcal A'\subseteq V_n(X^\beta)$ has the f.i.p. as well. As $*_{\beta(\beta+1)}:V(X^\beta)\to V(X^{\beta+1})$ is an enlargement, it follows that\[\bigcap *_{\beta(\beta+1)}[\mathcal A']\neq\varnothing.\] Now if $x\in \bigcap *_{\alpha(\alpha+1}[\mathcal A']$, then $x\in *_{\beta(\beta+1)}a_\xi'$ for all $\xi<\kappa$. It follows that $*_{(\beta+1)\kappa^+}(x)\in *_{(\beta+1)\kappa^+}\circ *_{\beta(\beta+1)}(a_\xi')=*_{\beta\kappa^+}(a'_\xi)=A_\xi$ for all $\xi<\kappa$, and hence that .$*_{(\beta+1)\kappa^+}(x)\in\bigcap \mathcal A$. 
Thus $\bigcap\mathcal A\neq\varnothing$ whenever $\mathcal A\subseteq W_n$ is a family of $\leq\kappa$-many sets with the f.i.p.

Observe that $X^{\kappa^+}= *_{\alpha\kappa^+}(X_\alpha)$ for all $\alpha<\kappa^+$. We now show that each $*_{\alpha\kappa^+}:V(X^\alpha)\to V(X^{\kappa^+})$ is $\kappa^+$-saturated.

 By Theorem \ref{thm_fip_element}, we need only show that every family $\mathcal A\subseteq {}^*V_n(X^\alpha)$ of cardinality $\kappa$ with the f.i.p. has non-empty intersection. Let $\mathcal A$ be such a family. Observe that if $A\in\mathcal A$, then $A\in{}^*V_n(X^\alpha)=*_{\alpha\kappa^+}(V_n(X^\alpha))=[V_n(X^\alpha),\alpha]\in W_{n+1}$, and so $A\in W_n$, by transitivity of $W$ over sets and Lemma \ref{lemma_sim_E_rank_P}(b). It follows that $\mathcal A\subseteq W_{n}$. By what we have just seen, $\bigcap\mathcal A\neq\varnothing$.

Since $\kappa=|V(X)|$, it follows that if we define $*:=*_{0\kappa^+}$ and ${}^*X:=X^{\kappa^+}$, then the map $*:V(X)\to V({}^*X)$ is polysaturated.


\appendix
\section{A Refresher on Basic First-Order Logic and Model Theory} \label{appendix_logic}

\subsection{First-Order Languages, Models and Satisfaction}
A first-order language  $\mathcal L=(\mathcal R,\mathcal F)$ consists of a collection of {\em relation symbols} (predicate symbols) $\mathcal R$ and {\em function symbols} $\mathcal F$. If $\mathcal R=\{R_1,\dots, R_n\}$ and $\mathcal F=\{F_1,\dots,F_m\}$ are finite sets, we may write $\mathcal L =(R_1,\dots, R_n, F_1,\dots, F_m)$.

A  first-order structure (or model) for $\mathcal L$ is a set equipped with relations and functions that interpret these symbols. We will define what this means shortly. With every $R\in\mathcal R$ is associated an {\em arity} $n\in\mathbb N$, which indicates that $R$ is to be interpreted as an $n$-ary relation. Similarly, with every $F\in\mathcal F$ is associated an {\em arity} $n\in\mathbb N$, which indicates that $F$ is to be interpreted as an $n$-ary function. The function symbols of arity 0 are to be interpreted as {\em constants}. 

\begin{definition} \rm ($\mathcal L$-structure) Let $\mathcal L=(\mathcal R,\mathcal F)$ be a first-order language. An $\mathcal L$-structure is a tuple $\mathfrak A=(A,\mathcal L^\mathfrak A)$ where $\mathcal L^\mathfrak A=(\mathcal R^\mathfrak A,\mathcal F^\mathfrak A)$ consists of relations and functions on the set $A$. Specifically, for each $n$-ary relation symbol $R\in\mathcal R$ there corresponds an $n$-ary relation $R^\mathfrak A\in\mathcal R^\mathfrak A$ on $A$, and to each $n$-ary function symbol $F\in\mathcal F$ there corresponds an $n$-ary function $F^\mathfrak A\in\mathcal F^\mathfrak A$ on $A$ such that\[\mathcal R^\mathfrak A=\{R^\mathfrak A: R\in\mathcal R\},\qquad \mathcal F^\mathfrak A=\{F^\mathfrak A: F\in\mathcal F\}.\] 
In particular if $c\in\mathcal F$ is a nullary function, then $c^\mathfrak A$ is a constant element of $A$.\\The set $A$ is called the {\em universe} of $\mathfrak A$, and $\mathfrak A$ is said to be a model of $\mathcal L$.\\
If $\mathcal R=\{R_1,\dots, R_n\}$ and $\mathcal F=\{F_1,\dots,F_m\}$ are finite sets, we may write $\mathfrak A=(A, R_1^\mathfrak A, \dots R_n^\mathfrak A, F_1^\mathfrak A,\dots, F_m^\mathfrak A)$.  When the interpretation is clear, we may dispense with the $\mathfrak A$-superscripts entirely, and simply write $(A,R_1,\dots, R_n, F_1,\dots, F_m)$.
\endbox
\end{definition}

For example, if $X$ is a base set, then the superstructure $\mathbb U:=(V(X),\in)$ is a $\mathcal L_\in$-structure.

Apart from the relation-- and function symbols which define it, the first-order languages that we consider also come with various other symbols including:
\begin{itemize}\item Countably many {\em variable symbols} $x_n$ ($n\in\mathbb N$) --- But we will often use $x,y,z, \dots$ instead.
\item The {\em equality symbol} $=$, which is always to be interpreted as equality.
\item {\em Logical connectives} $\lnot$ (not) and $\land$ (and).
\item The {\em universal quantifier} $\forall$ (for all).
\item Punctuation symbols such as parentheses and commas.
\end{itemize}
 
 \begin{definition}\rm (Terms and Formulas)  Consider a first-order language $\mathcal L=(\mathcal R,\mathcal L)$.
 \begin{enumerate}[(a)]\item The  {\em terms} of $\mathcal L$ are defined inductively:
 \begin{enumerate}[(i)]\item Every variable and every constant symbol is a term.
 \item If $F\in\mathcal F$ is a $n$-ary function symbol and $t_1,\dots, t_n$ are terms, then the string $F(t_1,\dots, t_n)$ is a term.
 \item A string is a term if and only if it can be obtained via a finite number of applications of (i), (ii).
 \end{enumerate}

 \item The {\em atomic formulas} of $\mathcal L$ are the expressions of the  following type:
 \begin{enumerate}[(i)]\item $s=t$, where $s,t$ are terms.
 \item $R(t_1,\dots, t_n)$, where $R$ is an $n$-ary relation symbol and $t_1,\dots, t_n$ are terms.
 \end{enumerate}
 \item
 The  {\em formulas} of a language $\mathcal L=(\mathcal R,\mathcal L)$ are defined inductively:
  \begin{enumerate}[(i)]\item Every atomic formula is a formula.
  \item If $\varphi, \psi$ are formulas and $x$ is a variable, then $\lnot\varphi, (\varphi\land\psi)$ and $(\forall x)\;\varphi$ are formulas.
   \item A string is a formula if and only if it can be obtained via a finite number of applications of (i), (ii).
  \end{enumerate}
   \end{enumerate}\endbox
  \end{definition}

We will also introduce a few other symbols, to simplify notation, namely $\lor$ (or), $\to$ (then, implies), $\leftrightarrow$ (if and only if) and $\exists$ (there exists). Suppose that $\varphi,\psi$ are formulas and that $x$ is a variable.
\begin{itemize}\item $(\varphi\lor\psi)$ is an abbreviation for $\lnot(\lnot\varphi\land\lnot\psi)$.
\item $(\varphi\to\psi)$ is an abbreviation for $(\lnot\varphi\lor\psi)$.
\item $(\varphi \leftrightarrow\psi)$ is an abbreviation for  $((\varphi\to\psi)\land(\psi\to\varphi))$.
\item $(\exists x)\;\varphi$ is an abbreviation for $\lnot(\forall x)\;\lnot\varphi$.
\end{itemize}
In an effort to make formulas more readable, we may omit parentheses, or replace parentheses with brackets, etc.

A formula $\psi$ is a {\em subformula} of a formula $\varphi$ if $\psi$ is a consecutive string of symbols within the formula $\varphi$.

If $\varphi$ is a formula, then a variable $x$ is  said to be within the {\em scope} of a quantifier $\forall x$ (or $\exists x$) occurring in $\varphi$ if there is a subformula of the form $\forall x\;\psi$ (or $\exists x\;\psi)$ such that $x$ occurs in $\psi$.
A variable $x$ may occur a number of times within a formula $\varphi$. An occurrence of a variable $x$ in formula is said to be {\em bound} if it occurs within the scope of a quantifier; otherwise, the occurrence is said to be {\em free}.

A formula is said to be a {\em sentence} if it has no free variables, i.e. if every occurrence of a variable is bound.

\begin{remarks}\rm For nonstandard universes, the appropriate language is $\mathcal L_\in$, consisting of just one binary relation symbol $\in$. In addition, we typically work with a modification of the first order language, where the quantifiers are bounded, i.e. of the form $\forall y\in x$ and $\exists y\in x$.
\endbox\end{remarks}

\begin{definition}\label{defn_interpretation_satisfaction} \rm (Interpretation and Satisfaction) Let $\mathfrak A=(A,\mathcal L^\mathfrak A)$ be a model for a first-order language $\mathcal L$.
For this definition, we will  write $t(x_1,\dots, x_n)$ if  $t$ is a term {\em all} of whose variables are among $x_1,\dots, x_n$ --- they need not all occur, however. Similarly, if {\em all } the variables , we write $\varphi(x_1,\dots, x_n)$ if $\varphi$ is a formula {\em all} of whose variables are among $(x_1,\dots, x_n)$ --- again, they need not all occur.
\begin{enumerate}[(a)]
\item {\bf Terms:} The value $t[a_1,\dots, a_n]\in A$ of a term $t(x_1,\dots, x_n)$ at $a_1,\dots, a_n\in A$ is defined inductively as follows:
\begin{enumerate}[(i)] \item If $t\equiv x_i$ is a variable, then $t[a_1,\dots, a_n]:=a_i$.\newline
If $t\equiv c$ is a constant (i.e. a nullary function symbol), then $t[a_1,\dots, a_n]:=c^\mathfrak A$.
\item If $t\equiv F(t_1,\dots,t_m)$, where $F$ is an $m$-ary function symbol and $t_1,\dots, t_m$ are terms, then
\[t[a_1,\dots, a_n]:=F^\mathfrak A(t_1[a_1,\dots, a_n],\dots, t_m[a_1,\dots, a_n]).\]
\end{enumerate}
\item {\bf Formulas:} For a formula $\varphi(x_1,\dots, x_n)$, the {\em satisfaction relation} $\mathfrak A\vDash \varphi[a_1,\dots, a_n]$ is defined inductively, as follows.
\begin{enumerate}[(i)] \item If $\varphi$ is the atomic formula $t_1=t_2$, then $\mathfrak A\vDash \varphi[a_1,\dots, a_n]$ if and only if $t_1[a_1,\dots, a_n]=t_2[a_1,\dots, a_n]$.\\Similarly, if $\varphi$ is the atomic formula $R(t_1,\dots, t_m)$, then $\mathfrak A\vDash \varphi[a_1,\dots, a_n]$ if and only if the relation $R^\mathfrak A(t_1[a_1,\dots, a_n], \dots, t_m[a_1,\dots, a_n])$ holds in $\mathfrak A$.
\item If $\varphi$ is the formula $\lnot\psi$, then $\mathfrak A\vDash \varphi[a_1,\dots, a_n]$ if and only if not $\mathfrak A\vDash \psi[a_1,\dots, a_n]$.
\item If $\varphi$ is the formula $\psi\land \chi$, then $\mathfrak A\vDash \varphi[a_1,\dots, a_n]$ if and only if both $\mathfrak A\vDash \psi[a_1,\dots, a_n]$ and $\mathfrak A\vDash \chi[a_1,\dots, a_n]$ .
\item If $\varphi$ is the formula $\forall x_i\;\psi$ (where $x_i\in \{x_1,\dots, x_n\}$), then \[\mathfrak A\vDash \varphi[a_1,\dots, a_n]\quad\text{if and only if for every $a\in A$,}\quad \mathfrak A\vDash \psi[a_1,\dots, a_{i-1}, a,a_{i+1}, \dots, a_n].\]
\end{enumerate}
\end{enumerate}
\endbox\end{definition}

It is easy to verify that:
\begin{itemize}\item $\mathfrak A\vDash (\varphi\lor\psi)[a_1,\dots, a_n]$ if and only if $\mathfrak A\vDash \varphi[a_1,\dots, a_n]$ or $\mathfrak A\vDash \psi[a_1,\dots, a_n]$.
\item $\mathfrak A\vDash (\varphi\to\psi)[a_1,\dots, a_n]$ if and only if whenever $\mathfrak A\vDash \varphi[a_1,\dots, a_n]$, then also $\mathfrak A\vDash \psi[a_1,\dots, a_n]$.
\item $\mathfrak A\vDash (\exists x_i\;\varphi)[a_1,\dots, a_n]$ if and only if there is $a \in A$ such that $\mathfrak A\vDash \psi[a_1,\dots, a_{i-1}, a, a_{i+1},\dots,a_n]$.
\end{itemize}

The next lemma show that whether or not a $\mathfrak A\vDash \varphi[a_1,\dots, a_n]$ depends only on those $a_i$ which correspond to variables $x_i$ that have a free occurrence in $\varphi$.

\begin{lemma}
\begin{enumerate}[(a)]\item If $t$ is a term with variables among $x_1,\dots, x_n$, then $t[a_1,\dots, a_n]$ depends only on the values $a_i$ corresponding to variables $x_i$ which actually occur in $t$.  \\More precisely, suppose that the variables occurring in a term $t$ are among $x_1,\dots, x_n$. Suppose further that $a_1,\dots, a_p$ and $b_1,\dots, b_q$ are elements of $A$, where $p,q\geq n$, and that $a_i=b_i$ whenever $x_i$ actually occurs in $t$. Then $t[a_1,\dots, a_p]=t[b_1,\dots, b_q]$.
\item  If $\varphi$ is a formula with variables among $x_1,\dots, x_n$, then whether or not $\mathfrak A\vDash\varphi[a_1,\dots, a_n]$ depends only on those $a_i$ for which $x_i$ has a free occurrence in $\varphi$. More precisely, suppose that the variables occurring in a formula $\varphi$ are among $x_1,\dots, x_n$, where some occur freely, and the others bound. Suppose further that $a_1,\dots, a_p$ and $b_1,\dots, b_q$ are elements of $A$, where $p,q\geq n$, and $a_i=b_i$ whenever $i\leq n$ and $x_i$ has a free occurrence in $\varphi(x_1,\dots,x_n)$.

Then
\[\mathfrak A\vDash \varphi[a_1,\dots, a_p] \Longleftrightarrow \mathfrak A\vDash \varphi[b_1,\dots, b_q].\] 
\end{enumerate}\end{lemma}
\bproof  These facts are is easily proved via induction on the length of $t, \varphi$.\\
(a) If $t$ is the variable $x_i$, then $t[a_1,\dots, a_p]=t[b_1,\dots, b_q]$ whenever $a_i=b_i$. If $t$ is a constant symbol, the result is obvious. If $t$ is the term $F(t_1,\dots, t_m)$, then
\\$t[a_1,\dots, a_p]= F^\mathfrak A(t_1[a_1,\dots, a_p],\dots, t_m[a_1,\dots, a_p])$. But the length of $t_1,\dots, t_n$ is clearly less than the length of $t$, so by induction we have that $t_i[a_1,\dots,a_p]=t_i[b_1,\dots, b_q]$ for $i\leq m$. Clearly, therefore $t[a_1,\dots, a_p]=t[b_1,\dots, b_q]$.\\
(b) If $\varphi$ is an atomic formula, then all variables that occur in $\varphi$ are free, and hence the result follows by (a). If $\varphi$ is of the form $\lnot\psi$, then the free variables of $\varphi$ are the same as the free variables of $\psi$. Now
$\mathfrak A\vDash \varphi[a_1,\dots, a_p]$ if and only if not $\mathfrak A\vDash \psi[a_1,\dots, a_p]$. But as the length of the formula $\psi$ is shorter than that of $\varphi$, we have 
$\mathfrak A\vDash \psi[a_1,\dots, a_n]\Leftrightarrow\mathfrak A\vDash \psi[b_1,\dots, b_q]$, and hence $\mathfrak A\vDash \varphi[a_1,\dots, a_p]\Leftrightarrow \mathfrak A\vDash \varphi[b_1,\dots, b_q]$. The case where $\varphi$ is of the form $\psi\land\chi$ is dealt with in a similar fashion. Finally, if $\varphi$ is of the form $\forall x_i\;\psi$, then the free variables of $\psi$ are just the free variables of $\varphi$, plus (possibly) the variable $x_i$. Now 
\begin{align*}&\phantom{aaaaa} \mathfrak A\vDash \varphi[a_1,\dots, a_p],\\&\Leftrightarrow\quad \text{for every $a\in A$, }\mathfrak A\vDash \psi[a_1,\dots, a_{i-1},a,  a_{i+1},\dots, a_p],\\
&\Leftrightarrow\quad \text{for every $a\in A$, }\mathfrak A\vDash \psi[b_1,\dots, b_{i-1},a,  b_{i+1},\dots, b_q], \quad \text{(induction hypothesis)}\\
&\Leftrightarrow\quad \mathfrak A\vDash \varphi[b_1,\dots, b_q].
\end{align*}\eproof

By the above lemma,  we may henceforth write $\varphi(x_1,\dots, x_n)$ to indicate a formula whose {\em free} variables are among $x_1,\dots, x_n$.
\begin{corollary} Truth-values of sentences are fixed: If $\varphi$ is a sentence, then either $\mathfrak A\vDash \varphi[a_1,\dots, a_n]$ for {\em all} sequences $a_1,\dots, a_n\in  A$, or for none of them. \endbox\end{corollary}

Suppose that $\Sigma$ is a a set of $\mathcal L$-sentences, and that $\mathfrak A$ is an $\mathcal L$-model. We say that $\mathfrak A$ is a model of $\Sigma$, or that $\mathfrak A$ satisfies $\Sigma$ --- and write $\mathfrak A\vDash\Sigma$ --- if and only if $\mathfrak A\vDash \varphi$ for every $\varphi\in\Sigma$.

\subsection{Elementary Embeddings and Elementary Equivalence}
\begin{definition}\rm (Submodel) Suppose that $\mathfrak A=(A,\mathcal L^\mathfrak A)$ and $\mathfrak B=(B,\mathcal L^\mathfrak B)$ are models of a first-order language $\mathcal L$. We say that $\mathfrak A$ is a {\em submodel} of $\mathfrak B$ --- and write $\mathfrak A\subseteq\mathfrak B$ --- if and only if 
\begin{enumerate}[(a)]\item $A\subseteq B$.
\item If $R$ as an $n$-ary relation symbol of $\mathcal L$, then $R^\mathfrak A=R^\mathfrak B\restriction A$, i.e. for any $a_1,\dots, a_n\in A$, we have that $R^\mathfrak A(a_1,\dots, a_n)$ holds in $\mathfrak A$ if and only if $R^\mathfrak B(a_1,\dots, a_n)$ holds in $\mathfrak B$.
\item Similarly, if $F$ as an $n$-ary function symbol of $\mathcal L$, then $F^\mathfrak A=F^\mathfrak B\restriction A$.\\In particular, if $c$ is a constant symbol, then $c^\mathfrak A=c^\mathfrak B$.
\end{enumerate}
\endbox\end{definition}

\begin{definition}\label{defn_theory_diagram}\rm (Theory, Elementary Diagram) Suppose that $\mathfrak A=(A,\mathcal L^\mathfrak A)$ is a model of a first-order language $\mathcal L$.
\begin{enumerate}[(a)]\item The {\em theory} of $\mathfrak A$ is the set
${Th}(\mathfrak A)$ of all $\mathcal L$-sentences that are satisfied by $\mathfrak A$.
\item Let $X\subseteq A$. By $\mathcal L_X$, we mean the language $\mathcal L$ augmented with additional constant symbols $\{c_a: a\in X\}$.  The model $\mathfrak A_X:=(\mathfrak A, a)_{a\in X}$ denotes the expansion of $\mathfrak A$ to a model of $\mathcal L_X$, where, for $a\in X$,  the new constant $c_a$ is interpreted to be the element $a\in X$.
\item The {\em  elementary diagram} $\Gamma_\mathfrak A$ of $\mathfrak  A$ is the theory $\text{Th}(\mathfrak A_A)$ of the expansion $\mathfrak A_A:=(\mathfrak A, a)_{a\in A}$, i.e. it is the set of all sentences of $\mathcal L_A$ which hold in $\mathfrak A_A$. 
\end{enumerate}\endbox\end{definition}

Note that if $\varphi(x_1,\dots, x_n)$ is an $\mathcal L$--formula and $\mathfrak A=(A,\mathcal L^\mathfrak A)$ is an $\mathcal L$-model, then for any $a_1,\dots, a_n\in A$, we obtain an $\mathcal L_A$-sentence  $\varphi(c_{a_1},\dots, c_{a_n})$ with the property that\[\mathfrak A_A\vDash \varphi(c_{a_1},\dots, c_{a_n}) \quad\text{if and only if }\quad \mathfrak A \vDash \varphi[a_1,\dots, a_n].\]
To simplify notation, we will write $\varphi(a_1,\dots, a_n)$ instead of $\varphi(c_{a_1},\dots, c_{a_n}) $. Then the elementary diagram of $\mathfrak A$ is\[\Gamma_\mathfrak A=\{\varphi(a_1,\dots, a_n): \mathfrak A\vDash \varphi[a_1,\dots, a_n]\}.\]

\begin{definition}\rm (Elementary Equivalence, Elementary Embedding) Suppose that $\mathfrak A,\mathfrak B$ are two models of a first-order language $\mathcal L$. 
\begin{enumerate}[(a)]\item 
We say that $\mathfrak A, \mathfrak B$ are elementarily equivalent --- and write $\mathfrak A\equiv\mathfrak B$ --- if and only if $\mathfrak A,\mathfrak B$ satisfy the same $\mathcal L$-sentences, i.e. $\text{Th}(\mathfrak A)=\text{Th}(\mathfrak B)$.
\item We say that $\mathfrak A$ is an elementary submodel of $\mathfrak B$ --- and write $\mathfrak A\preceq \mathfrak B$ --- if and only if $\mathfrak A\subseteq \mathfrak B$ and for all $\mathcal L$-formulas $\varphi(x_1,\dots, x_n)$ and all $a_1,\dots, a_n\in A$,\[\mathfrak A\vDash \varphi[a_1,\dots, a_n]\quad\text{if and only if}\quad\mathfrak B\vDash \varphi[a_1,\dots, a_n].\]
\item We say that a map $f:A\to B$ is an {\em elementary embedding} --- and write $f:\mathfrak A\overset{\equiv}{\hookrightarrow}\mathfrak B$ --- if and only if for all $\mathcal L$-formulas $\varphi(x_1,\dots, x_n)$ and all $a_1,\dots, a_n\in A$,\[\mathfrak A\vDash \varphi[a_1,\dots, a_n]\quad\text{if and only if}\quad\mathfrak B\vDash \varphi[f(a_1),\dots, f(a_n)].\]
\end{enumerate}\endbox
\end{definition}

Observe the following trivial fact: \begin{lemma} If there is an elementary embedding $\mathfrak A\overset{\equiv}{\hookrightarrow}\mathfrak B$, then $\mathfrak A
\equiv\mathfrak B$.\endbox\end{lemma}

The following lemma follows by chasing through the above definitions:
\begin{lemma}\label{lemma_diagram_embedding} Suppose that $\mathfrak A,\mathfrak B$ are two models of a first-order language $\mathcal L$.
\begin{enumerate}[(a)]\item If $\mathfrak A\subseteq\mathfrak B$, 
then $\mathfrak A\preceq\mathfrak B$ if and only if $(\mathfrak B,a)_{a\in A}\vDash \Gamma_\mathfrak A$.\newline
\item Similarly, there is an elementary embedding $\mathfrak A\overset{\equiv}{\hookrightarrow}\mathfrak B$ if and only if there is an expansion $(\mathfrak B, b_a)_{a\in A}$ of $\mathfrak B$ such that $(\mathfrak B, b_a)_{a\in A}\vDash \Gamma_\mathfrak A$ (and then $a\mapsto b_a$ supplies the required elementary embedding).
\end{enumerate}\endbox
\end{lemma}

\subsection{Ultrafilters}

\begin{definition}\label{defn_fip_ultrafilter} (Ultrafilter)
Let $I$ be a set. 
\begin{enumerate}[(a)]\item A family $\mathcal  A\subseteq\mathcal P(I)$ is said to satisfy the {\em finite intersection property} (f.i.p.) if and only whenever $\mathcal B:=\{A_1,\dots, A_n\}$ is any finite subcollection of $\mathcal  A$, then $\bigcap\mathcal B$ is non-empty, i.e. $A_1\cap\dots\cap A_n\neq\varnothing$.

\item A {\em filter} over $I$ is a non-empty family $\mathcal F\subseteq\mathcal P(I)$ with the following properties:
\begin{enumerate}[(i)]\item $\varnothing\not\in\mathcal F$.
\item $F,G\in\mathcal F$ implies $F\cap G\in\mathcal F$.
\item $F\in\mathcal F$ and $G\supseteq F$ implies $G\in\mathcal F$.
\end{enumerate}
\item A filter $\mathcal U$ over $I$ is said to be an {\em ultrafilter} if and only if for every $A\subseteq I$ we have either $A\in\mathcal U$ or its complement $A^c\in\mathcal U$.
\item A filter $\mathcal F$ is said to be {\em countably incomplete} if and only if there is a countable subfamily $\{F_n:n\in\mathcal F\}$ of members of $\mathcal F$ such that $\bigcap_nF_n\not\in\mathcal F$.
\end{enumerate} 
\end{definition}

In Section \ref{section_existence_nonstandard_frameworks} it will be shown that any ultrafilter will induce a transfer map $*:V(X)\to V(Y)$ between superstructures. In order for that transfer map to give rise to a nonstandard framework --- i.e. so that ${}^\sigma C\subsetneq{}^*C$ for some countable $C\in V(X)$ --- we will require that the ultrafilter is countably incomplete.

Here follow some basic facts:
\begin{theorem} Let $I$ be a set.
\begin{enumerate}[(a)]\item Every  filter on $I$ has the f.i.p.
\item if $\mathcal A$ is a family of subsets of a set $I$ with the f.i.p., then \[\mathcal F:=\{F\subseteq I: \text{there are }A_1,\dots, A_n\in\mathcal A\text{ such that } A_1\cap\dots\cap A_n\subseteq F\}\] is a filter containing $\mathcal A$.
\item If $\mathcal A$ is a family of subsets of a set $I$ with the f.i.p. and $B\subseteq I$, then either $\mathcal A\cup\{B\}$ has the f.i.p., or else $\mathcal  A\cup \{B^c\}$ has the f.i.p. 
\item A filter over $I$ is an ultrafilter if and only if it is a maximal filter, i.e. if and only if it is not contained in any strictly larger filter. 
\item If $\mathcal A$ is a family of subsets of a set $I$ with the f.i.p., then there is an ultrafilter $\mathcal U$ such that $\mathcal A\subseteq\mathcal U$.
\item $\mathcal U$ is an ultrafilter over $I$ if and only if whenever $A_1,\dots, A_n\subseteq I$ are such that $A_1\cup\dots\cup A_n\in\mathcal U$, then there is $i\leq n$ such that $A_i\in\mathcal U$.
\item  $\mathcal U$ is a countably incomplete ultrafilter if and only if there is a partition $I$ into a sequence $I_n$ $(n\in\mathbb N)$ of disjoint non-empty sets with the property that $\bigcap_nI_n=I$, and $I_n\not\in\mathcal U$ for any $n$. \end{enumerate}
\end{theorem}

\bproof
(a) is obvious.\\
(b) is straightforward.\\
(c)  Suppose that $\mathcal  A\cup\{B\}$ does not have the f.i.p., then there exists $A_1,\dots, A_n\in\mathfrak A$ such that $A_1\cap\dots\cap A_n\cap B=\varnothing$, so that $A_1\cap\dots\cap A_n\subseteq B^c$. Then if $A'_1,\dots, A'_m\in \mathcal  A$, we have $A'_1\cap\dots\cap A'_m\cap B^c\supseteq A'_1\cap\dots\cap A'_m\cap A_1\cap\dots\cap A_n\neq\varnothing$, as $\mathcal  A$ has the f.i.p. and $A'_i, A_j\in\mathcal A$. Thus $\mathcal A\cup\{B^c\}$ has the f.i.p.\\
(d) Suppose that $\mathcal U$ is an ultrafilter over $I$. If $\mathcal F$ is a strictly larger filter and $F\in\mathcal F-\mathcal U$, then $F^c\in\mathcal U$, and hence $F\cap F^c=\varnothing\in\mathcal F$ --- contradicting the definition of filter. Conversely, an easy application of Zorn's Lemma shows that any filter can be extended to a maximal filter.  If $\mathcal F$ is a maximal filter and $B\subseteq I$, then either $\mathcal F\cup\{B\}$ or $\mathcal F\cup \{B^c\}$ has the f.i.p. Hence there is a filter $\mathcal G$ such that $\mathcal G\supseteq\mathcal F\cup\{B\}$ or $\mathcal G\supseteq\mathcal F\cup\{B^c\}$. But as $\mathcal F$ is maximal, we must have $\mathcal G=\mathcal F$. Thus either $B\in\mathcal F$ or $B^c\in\mathcal F$, proving that $\mathcal F$ is an ultrafilter. Thus the ultrafilters are precisely the maximal filters.\\
(e) By (b), there is a filter $\mathcal F$ such that $\mathcal A\subseteq\mathcal F$. By an easy application of Zorn's Lemma, there is a maximal filter $\mathcal U$ over $I$ such that $\mathcal F\subseteq\mathcal U$. By (d), $\mathcal U$ is an ultrafilter.\\
(f) Suppose that $\mathcal U$ is an ultrafilter with $A_1\cup\dots\cup A_n\in\mathcal U$. If $A_i\not\in\mathcal U$ for any $i\leq n$, then also $A_1^c\cap\dots\cap A_n^c=(A_1\cup\dots\cup A_n)^c\in\mathcal U$, and hence $\varnothing\in\mathcal U$ --- contradicting the definition of filter. Conversely, if $\mathcal F$ is a filter over $I$ with the stated property and $B\subseteq I$, then $I=B\cup B^c\in\mathcal F$, and hence either $B\in\mathcal F$ or $B^c\in\mathcal F$. Hence $\mathcal F$ is an ultrafilter.\\
(g)  Suppose that $\mathcal U$ is countably incomplete. Then there are $U_n\in\mathcal U$ (for $n\in \mathbb N$) such that $\bigcap_nU_n\not\in\mathcal U$. Since $\mathcal U$ is closed under finite intersections, we may assume that $U_1\supseteq  U_2\supseteq U_3\supseteq\dots$, and that $U_1=I$. By removing duplicates, we may assume that all the $U_n$ are distinct.
Now define
\[I_0:=\bigcap_n U_n,\quad \text{and}\quad I_{n}:=U_n-U_{n+1}\text{ for }n>0.\] Then $I=\bigcup_n I_n$, and clearly the sets $I_n$ are non-empty and partition $I$. By assumption, $I_0\not \in \mathcal U$. Furthermore, since $U_{n+1}^c\supseteq U_n-U_{n+1} = I_n$, we cannot have $I_n\in\mathcal U$. 
\\ Conversely, if $\{I_n:n\in\mathbb N\}$ is a partition of $I$ into disjoint non-empty sets, then $\bigcap_{m\neq n}I_m^c=I_n\neq\varnothing$, from which it follows that $\{I_n^c:n\in\mathbb N\}$ has the f.i.p. It follows that there is an ultrafilter $\mathcal U$ such that $I_n^c\in\mathcal U$ for each $n\in\mathbb N$.  As $\bigcap_n I_n^c=\varnothing\not\in\mathcal U$, the ultrafilter $\mathcal U$ is countably incomplete.
\eproof

\subsection{Ultraproducts and Ultrapowers}

Suppose that $I$ is a set, and that $\mathcal F$ is a a filter over $I$. Suppose further that $\mathfrak A_i=(A_i,\mathcal L^{\mathfrak A_i})$ $(i\in I)$ are non-empty models of a first-order language $\mathcal L$.  Recall that the product $\prod_{i\in I}A_i$ is the set of all choice functions $f:I\to\bigcup_{i\in I}A_i$, i.e. all those functions with the property that $f(i)\in A_i$ for all $i\in I$. Define a binary relation 
$\sim_\mathcal F$ on the set $\prod_{i\in I}A_i$ by\[f\sim_\mathcal F g\quad \text{if and only if }\quad \{i\in I: f(i)=g(i)\}\in\mathcal F.\] It is straightforward to verify that $\sim_\mathcal F$ is an equivalence relation. For example, if $f\sim_\mathcal F g$ and $g_\mathcal F H$, then\[\{i\in I: f(i)=h(i)\}\supseteq\{i\in I: f(i)=g(i)\}\cap\{i\in I: g(i)=h(i)\}\in\mathcal F,\] from which it follows that $f\sim_\mathcal F h$, i.e. that $\sim_\mathcal F$ is transitive.

Denote the equivalence relation corresponding to $f\in\prod_IA_i$ by $f/\mathcal F$, and let\[\prod_IA_i/\mathcal F:=\{f/\mathcal F: f\in \prod_IA_i\}\] denote the corresponding quotient set.

We now show how to equip $\prod_IA_i/\mathcal F$ with relations and functions which turn it into a $\mathcal L$-structure $\mathfrak B=\prod_I\mathfrak A_i/\mathcal F$.
For $n$-ary relation-- and function symbols $R, F$ of $\mathcal L$, let $R^{\mathfrak A_i}, F^{\mathfrak A_i}$ denote their corresponding interpretations in the model $\mathfrak A_i$. Define the relation $R^\mathfrak B$ and function $F^\mathfrak B$ on $B:=\prod_IA_i/\mathcal F$ by
\[R^\mathfrak B(f_1/\mathcal F,\dots, f_n/\mathcal F) \quad\text{if and only if}\quad \{i\in I: R^{\mathfrak A_i}(f_1(i),\dots, f_n(i))\}\in\mathcal F,\]
and
\[F^\mathfrak B(f_1/\mathcal F,\dots, f_n/\mathcal F) :=g/\mathcal F,\quad\text{where }g(i):=F^{\mathfrak A_i}(f_1(i),\dots, f_n(i)).\]
It is easy to verify that these notions are well-defined. For example, if $f_i/\mathcal F=g_i/\mathcal F$ for $i\leq n$  and $\{i\in I: R^i(f_1(i),\dots, f_n(i))\}\in\mathcal F$, then also
\[\{i:R^i(g_1(i),\dots, g_n(i))\}\supseteq \{i: R^i(f_1(i),\dots, f_n(i))\}\cap\{i: f_1(i)=g_1(i)\}\cap\dots\cap\{i\in I: f_n(i)=g_n(i)\}\in\mathcal F.\]

The $\mathcal L$-structure with base set $\prod_IA_i$ and corresponding relations and functions $R^\mathfrak B, F^\mathfrak B$ (for $R,F\in\mathcal L$) is called the {\em reduced product } of the models $\mathfrak A_i$, and denoted by $\prod_I\mathfrak A_i/\mathcal F$. If the $\mathfrak A_i$ are all identical, then we have a {\em reduced power}, and denote it by $\mathfrak A^I/\mathcal F$. If $\mathcal F$ is an ultrafilter over $I$, then a reduced product is called an {\em ultraproduct}, and a reduced power is called an {\em ultrapower}.

\begin{lemma} Suppose that  $\mathfrak A_i$  $(i\in I)$ are models of a first-order language $\mathcal L$, and that $\mathcal F$ is a filter over $I$.  If $t(x_1,\dots, x_n)$ is an $\mathcal L$-term, then the interpretation $t^\mathcal F$ of $t$ in $\prod_I\mathfrak A_i/\mathcal F$ is given by \[t^\mathfrak B(f_1/\mathcal F,\dots,f_n/\mathcal F)=g/\mathcal F,\quad\text{where}\quad g(i):=t^{\mathfrak A_i}(f_1(i),\dots, f_n(i)).\] Similarly, if $R$ is an $m$-ary relation symbol of $\mathcal L$, and $t_1,\dots, t_m$ are $\mathcal L$-terms, then \\$R^\mathfrak B\Big(t_1^\mathfrak B(f_1/\mathcal F,\dots, f_n/\mathcal F),\dots, t_m^\mathfrak B(f_1/\mathcal F,\dots, f_n/\mathcal F)\Big)$ holds in $\prod_I\mathfrak A_i/\mathcal F$ if and only if \[\Big\{i\in I: R^{\mathfrak A_i}\Big(t_1^{\mathfrak A_i}(f_1(i),\dots, f_n(i)),\dots ,t_m^{\mathfrak A_i}(f_1(i),\dots, f_n(i))\Big)\Big\}\in\mathcal F.\]
\end{lemma}
\bproof Induction: If $t\equiv F(x_1,\dots, x_n)$ for some $n$-ary function symbol $F\in\mathcal L$, the result follows by definition of $F^\mathfrak B$. If $t\equiv F(t_1(x_1,\dots, x_n),\dots, t_m(x_1,\dots, x_n))$ for some $m$-ary function symbol $F$ and terms $t_1,\dots, t_n$,  then by the inductive definition of $t^\mathfrak B$ (cf. Definition \ref{defn_interpretation_satisfaction}), we have  \[t^\mathfrak B(f_1/\mathcal F,\dots, f_n/\mathcal F):=F^\mathfrak B\Big(t_1^\mathfrak B(f_1/\mathcal F,\dots, f_n/\mathcal F),\dots, t_m^\mathfrak B(f_1/\mathcal F,\dots, f_n/\mathcal F)\Big).\]
By induction hypothesis, we have, for $k\leq m$,\[t_k^\mathfrak B(f_1/\mathcal F,\dots, f_n/\mathcal F)= g_k/\mathcal F,\quad\text{where}\quad g_k(i):=t_k^{\mathfrak A_i}(f_1(i),\dots, f_n(i)).\]
Thus  \[t^\mathfrak B(f_1/\mathcal F,\dots, f_n/\mathcal F)= F^\mathfrak B(g_1/\mathcal F,\dots, g_m/\mathcal F) =g/\mathcal F\quad\text{where}\quad g(i):= F^{\mathfrak A_i}(g_1(i),\dots, g_m(i)).\] But then by the inductive definition of $t^{\mathfrak A_i}$ we have
\[g(i)= F^{\mathfrak A_i}\Big(t_1^{\mathfrak A_i}(f_1(i),\dots, f_n(i)),\dots, t_m^{\mathfrak A_i}(f_1(i),\dots, f_n(i))\Big)=t^{\mathfrak A_i}(f_1(i),\dots, f_n(i)),\]
from which the result follows.

Now consider the case $R(t_1,\dots, t_m)$. For $k\leq m$, let $t_k^\mathfrak B(f_1/\mathcal F,\dots,f_n/\mathcal F)=g_k/\mathcal F$, where $g_k(i):=t^{\mathfrak A_i}(f_1(i),\dots, f_n(i))$.
Then we  have $R^\mathfrak B\Big(t_1^\mathfrak B(f_1/\mathcal F,\dots, f_n/\mathcal F),\dots, t_m^\mathfrak B(f_1/\mathcal F,\dots, f_n\mathcal F)\Big)$ if and only if $\mathcal R^\mathfrak B(g_1/\mathcal F,\dots, g_m/\mathcal F)$ if and only if $\{i\in I: R^{\mathfrak A_i}(g_1(i),\dots, g_n(i))\}\in\mathcal F$ if and only if $\Big\{i\in I: R^{\mathfrak A_i}\Big(t_1^\mathfrak{A_i}(f_1(i),\dots, f_n(i)),\dots ,t_m^\mathfrak{A_i}(f_1(i),\dots, f_n(i))\Big)\Big\}\in\mathcal F$.
\eproof

\begin{theorem}{\rm (\L os)} Suppose that  $\mathfrak A_i$  $(i\in I)$ are models of a first-order language $\mathcal L$, and that $\mathcal U$ is an ultrafilter over $I$.  If $\varphi(x_1,\dots, x_n)$ is an $\mathcal L$-formula whose free variables are among $x_1,\dots, x_n$, then
\[\prod_I\mathfrak A_i\Big/\mathcal U\vDash \varphi[f_1/\mathcal U,\dots, f_n/\mathcal U] \quad \Longleftrightarrow\quad \{i\in I: \mathfrak A_i\vDash\varphi[f_1(i),\dots, f_n(i)]\}\in\mathcal U.\]
\end{theorem}
\bproof Induction on the complexity of $\varphi$. If $\varphi$ is an atomic sentence, then the result follows from the previous lemma. If $\varphi$ is of the form $\psi\land\chi$, then by the recursive definition of the satisfaction relation (cf. Definition \ref{defn_interpretation_satisfaction}) and the induction hypothesis, we have
\begin{align*}\phantom{aaaaaaaaaaaa}&\prod_I\mathfrak A_i/\mathcal U\vDash\varphi[f_1/\mathcal U,\dots,f_n/\mathcal U] \\&\Longleftrightarrow\quad \prod_I\mathfrak A_i/\mathcal U\vDash \psi[f_1/\mathcal U,\dots,f_n/\mathcal U]\quad \text{and}\quad \prod_I\mathfrak A_i/\mathcal U\vDash \chi[f_1/\mathcal U,\dots,f_n/\mathcal U],\\&\Longleftrightarrow \quad\{i:\mathfrak A_i\vDash \psi[f_1(i),\dots, f_n(i)]\}\in\mathcal U\quad\text{and}\quad \{i:\mathfrak A_i\vDash \chi[f_1(i),\dots, f_n(i)]\}\in\mathcal U,\\&\Longleftrightarrow \quad \{i:\mathfrak A_i\vDash \varphi[f_1(i),\dots, f_n(i)]\}\in\mathcal U,
\end{align*}
using the fact that $\{i:\mathfrak A_i\vDash \varphi\}=\{i:\mathfrak A_i\vDash \psi\}\cap\{i:\mathfrak A_i\vDash \chi\}$.

If $\varphi$ is of the form $\lnot\psi$, then, using the fact that $\mathcal U$ is an ultrafilter, 
\begin{align*}\phantom{aaaaaaaaaaaa}&\prod_I\mathfrak A_i/\mathcal U\vDash\varphi[f_1/\mathcal U,\dots,f_n/\mathcal U],
 \\&\Longleftrightarrow\quad \prod_I\mathfrak A_i/\mathcal U\not\vDash \psi[f_1/\mathcal U,\dots,f_n/\mathcal U],
 \\&\Longleftrightarrow\quad\{i:\mathfrak A_i\vDash \psi[f_1(i),\dots, f_n(i)]\}\not\in\mathcal U,
  \\&\Longleftrightarrow\quad\{i:\mathfrak A_i\vDash \lnot\psi[f_1(i),\dots, f_n(i)]\}\in\mathcal U,
  \\&\Longleftrightarrow \quad \{i:\mathfrak A_i\vDash \varphi[f_1(i),\dots, f_n(i)]\}\in\mathcal U.\end{align*}

Finally, if $\varphi$ is of the form $\forall y\;\psi(y, x_1,\dots, x_n)$, then
\begin{align*}\phantom{aaaaaaaaaaaa}&\prod_I\mathfrak A_i/\mathcal U\vDash\varphi[f_1/\mathcal U,\dots,f_n/\mathcal U],
 \\&\Longleftrightarrow\quad\prod_I\mathfrak A_i/\mathcal U\vDash\psi[g/\mathcal U, f_1/\mathcal U,\dots,f_n/\mathcal U]\quad \text{for all }g/\mathcal U\in\prod_IA_i/\mathcal U,
  \\&\Longleftrightarrow\quad\{i\in I:\mathfrak A_i\vDash \psi(g(i), f_1(i),\dots, f_n(i))\}\in\mathcal U\quad \text{for all }g\in\prod_IA_i,
    \\&\Longleftrightarrow\quad\{i\in I:\mathfrak A_i\vDash \forall y\;\psi(y, f_1(i),\dots, f_n(i))\}\in\mathcal U\quad \text{for all }g\in\prod_IA_i,\\
    &\qquad\qquad\qquad\text{(since $g(i)$ can be chosen to be any member of $A_i$ whatsoever)}
      \\&\Longleftrightarrow \quad \{i:\mathfrak A_i\vDash \varphi[f_1(i),\dots, f_n(i)]\}\in\mathcal U.
 \end{align*}\eproof
 
 \begin{corollary} Suppose that $\mathfrak A=(A,\mathcal L^\mathfrak A)$ is a model of a first-order language $\mathcal L$. Let $\mathcal U$ be an ultrafilter over the set $I$. For $a\in A$,  let $c_a\in A^I$ be the constant map with value $a$.Then the map $h:\mathfrak A\to \mathfrak A^I/\mathcal U:a\mapsto c_a/\mathcal U$ is an elementary embedding.\newline Hence $\mathfrak A^I/\mathcal U\equiv \mathfrak A$.
 \end{corollary}
 \bproof If $\varphi(_1,\dots, x_n)$ is an $\mathcal L$-formula, then by \L os' Theorem,
 \[\mathfrak A^I/\mathcal U\vDash \varphi[c_{a_1}/\mathcal U, \dots, c_{a_n}/\mathcal U] \quad\text{if and only if}\quad \{i\in I:\mathfrak A\vDash \varphi[a_1,\dots, a_n]\}\in\mathcal U.\] But $\{i\in I:\mathfrak A\vDash \varphi[a_1,\dots, a_n]\}$ is either $I$ (if $\mathfrak A\vDash \varphi[a_1,\dots, a_n]$), or $\varnothing$ (if $\mathfrak A\not\vDash \varphi[a_1,\dots, a_n]$).\eproof
 
 \begin{theorem}\label{thm_compactness} Let $\Sigma$ be a set of sentences of a first order language. If every finite subset of $\Sigma$ has a model, then $\Sigma$ has a model.
 \end{theorem}
 \bproof Let $I:=\mathcal P^{<\omega}(\Sigma)-\{\varnothing\}$ be the set of all non-empty finite subsets of $\Sigma$. For each $i\in I$, let $\mathfrak A_i$ be an $\mathcal L$-model so that $\mathfrak A_i\vDash i$. For $i\in I$, define $I_i\subseteq I$ by\[I_i:=\{j\in I: i\subseteq j\}.\] Observe that 
 $I_i\cap I_j=I_{i\cup j}$, so that the family $\{I_i: i\in I\}$ has the f.i.p. Let $\mathcal U$ be an ultrafilter over $I$ such that $\{I_i: i\in I\}\subseteq\mathcal U$. We claim that $\prod_I\mathfrak A_i/\mathcal U\vDash\Sigma$.
 
 For suppose that $\varphi\in\Sigma$, and let $i_0:=\{\varphi\}\in I$. Observe that $I_{i_0}=\{j\in I: \varphi\in j\}$. Now $\varphi\in j$ implies $\mathfrak A_j\vDash \varphi$, since $\mathfrak A_j\vDash j$. Thus
 \[\{j\in I: \mathfrak A_j\vDash \varphi\}\supseteq \{j\in I: \varphi\in j\}=I_{i_0}\in\mathcal U.\] Hence $\{j\in I: \mathfrak A_j\vDash \varphi\}\in\mathcal U$, so that by \L os' Theorem $\prod_I\mathfrak A_i/\mathcal U\vDash \varphi$. As this is true for any $\varphi\in \Sigma$, it follows that $\prod_I\mathfrak A_i/\mathcal U\vDash \Sigma$.
 \eproof

\section{Existence of Good Ultrafilters}\label{appendix_good_ultrafilters}
\fancyhead[LO]{Existence of Good Ultrafilters}

Recall the following definition: 
\begin{definition}\rm Let $\kappa$ be a cardinal and let $\mathcal E\subseteq\mathcal P(\kappa)$.
\begin{enumerate}[(a)]
\item An {\em order-reversal} is a map $p:\mathcal P^{<\omega}(\kappa)\to\mathcal E$ such that  whenever $s,t \in \mathcal P^{<\omega}(\kappa)$ and $s\subseteq t$, then $p(s)\supseteq p(t)$.
\item An {\em anti-additive} map is a map $p:\mathcal P^{<\omega}(\kappa)\to\mathcal E$ such that if $s,t \in \mathcal P^{<\omega}(\kappa)$, then $p(s\cup t)=p(s)\cap p(t)$. Clearly every anti-additive map is an order-reversal.
\item 
 An ultrafilter $\mathcal U$ over a set $I$ is {\em $\kappa$-good} if for every cardinal $\alpha<\kappa$ and every  order-reversal $p: \mathcal P^{<\omega}(\alpha)\to\mathcal U$ there is an anti-additive map $q: \mathcal P^{<\omega}(\alpha)\to\mathcal U$ such that $q\leq p$, i.e. such that $q(s)\subseteq p(s)$ for all $s\in\mathcal P^{<\omega}(\kappa)$.
 \end{enumerate}\endbox
\end{definition}

The following lemma slightly simplifies the verification of the  goodness condition for successor cardinals:

\begin{lemma} Let $\kappa$ be a cardinal, and let $\mathcal U$ be an ultrafilter over a set $I$. Then $\mathcal U$ is $\kappa^+$--good if and only if  for every order-reversal $p:\mathcal P^{<\omega}(\kappa)\to \mathcal U$ there is an anti-additive $q:\mathcal P^{<\omega}(\alpha)\to \mathcal U$ such that $q\leq p$. \end{lemma}\
\bproof The  $(\Rightarrow)$-direction is obvious.\\
 For the $(\Leftarrow)$-direction, suppose that $\alpha\leq \kappa$ and that $p:\mathcal P^{<\omega}(\alpha)\to \mathcal U$ is an order-reversal. Note that if $s$ is a finite subset of $\kappa$, then $s\cap\alpha$ is a finite subset of $\alpha$. Thus we can define an order-reversing map \[\bar{p}:\mathcal P^{<\omega}(\kappa)\to \mathcal U:a\mapsto p(s\cap\alpha).\]
By assumption there is an anti-additive $\bar{q}:\mathcal P^{<\omega}(\kappa)\to \mathcal U$ such that $\bar{q}\leq\bar{p}$.Define $q$ to be the restriction  $q:=\bar{q}\restriction\mathcal P^{<\omega}(\alpha)$. Then $q:\mathcal P^{<\omega}(\alpha)\to\mathcal U$ is an clearly anti-additive map. In addition, if $s\in\mathcal P^{<\omega}(\alpha)$, then  \[q(s)=\bar{q}(s)\subseteq\bar{p}(s)=p(s\cap\alpha)=p(s),\] so $q\leq p$.\eproof

For an infinite cardinal $\kappa$, let $\mathcal P^*(\kappa)$ be the set of all subsets of $\kappa$ with cardinality $\kappa$:\[\mathcal P^*(\kappa)=\{X\subseteq\kappa: |X|=\kappa\}.\]
If $h:X\to\mathcal Y$ is a function whose range is a family of sets, then we say that $h$ is a {\em disjoint function}  if the sets $h(x), x\in X$ are disjoint.
\begin{lemma} \label{lemma_dominate_disjoint} Suppose that $\kappa$ is an infinite cardinal, and that $f:\kappa\to\mathcal P^*(\kappa)$. Then there is a disjoint function $h:\kappa\to\mathcal P^*(\kappa)$ such that $h\leq f$.
\end{lemma}

\bproof
Since $|\kappa\time\kappa|=\kappa$, there is a bijection $\iota:\kappa\to\kappa\times\kappa$. Define $\mathcal G$ to be the set of all one-to-one functions $g$ with the following properties:
\begin{enumerate}[(i)]\item $\dom(g)\in\kappa$.
\item If $\xi\in \dom(g)$, then $g(\xi)\in f(\pi_0\circ\iota(\xi))$, where $\pi_0:\kappa\times\kappa\to \kappa:(\zeta,\eta)\mapsto \zeta$.\\ Thus if $\xi\in \iota^{-1}[\{\zeta\}\times\kappa]$, then $g(\xi)\in f(\zeta)$.
\end{enumerate}
Let $G\subseteq\mathcal G$ be a maximal chain, ordered by inclusion, and define $g^*:=\bigcup G$. Then clearly $g^*$ is a one-to-one function satisfying (ii), with $\dom(g^*)\leq\kappa$. We shall show that $\dom(g^*)=\kappa$. For suppose that $\dom(g^*)=\alpha<\kappa$. Then $\iota(\alpha)=(\zeta,\eta)$ for  some $\zeta,\eta<\kappa$, and $|\ran(g^*)|=|\alpha|<\kappa$. As $|f(\zeta)|=\kappa$, we may choose $\gamma\in f(\zeta)-\ran(g^*)$, and define $g'\in\mathcal G$ by $g'=g^*\cup\{(\alpha, \gamma)\}$. Then $g'$ is a one-to-one function satisfying (ii) with $\dom(g')=\alpha+1$, and hence $G\cup\{g'\}$ is a chain in $\mathcal G$ which extends $G$ --- contradicting the maximality of $G$.

Now define $h:\kappa\to\mathcal P^*(\kappa)$ by\[h(\zeta):=\{g^*(\xi): \xi\in\iota^{-1}[\{\zeta\}\times\kappa]\}=\{g^*(\xi):\pi_0\circ \iota(\xi)=\zeta\}.\] Then clearly $|h(\zeta)|=\kappa$ as $|\iota^{-1}[\{\zeta\}\times\kappa]|=\kappa$. In addition, if $x\in h(\zeta)$, then $x=g^*(\xi)\in f(\pi_0\circ\iota(\xi))=f(\zeta)$ (by (ii)), and hence $h(\zeta)\subseteq f(\zeta)$.Moreover, the sets $h(\zeta), \zeta<\kappa$ are disjoint, as $g^*$ is one-to-one and the sets $\iota^{-1}[\{\zeta\}\times\kappa]$ (for $\zeta<\kappa$) are disjoint.
\eproof

The following corollary is merely a restatement of the previous lemma: 

\begin{corollary}\label{cor_dominate_disjoint}Suppose that $\kappa$ is an infinite cardinal, and that $\{X_\alpha:\alpha<\kappa\}$ is a family of sets, each of cardinality $\kappa$. Then there exist sets $Y_\alpha,\alpha<\kappa$ such that
\begin{enumerate}[(i)]\item Each $Y_\alpha\subseteq X_\alpha$.
\item Each $|Y_\alpha|=\kappa$,
\item The sets $Y_\alpha, \alpha<\kappa$ are mutually disjoint.
\end{enumerate}
\end{corollary}
\bproof Since $|\bigcup_{\alpha<\kappa}X_\alpha|=\kappa$, we may assume without loss of generality that each $X_\alpha\in\mathcal P^*(\kappa)$. Define $f:\kappa\to\mathcal P^*(\kappa):\alpha\mapsto X_\alpha$, and apply  Lemma \ref{lemma_dominate_disjoint} to obtain $h\leq f$. Now define $Y_\alpha:=h(\alpha)$.
\eproof

\begin{definition}\rm
\begin{enumerate}[(a)]\item 
A  {\em partition} $P$ of a set $X$ is a collection of disjoint subsets of $X$ whose union is $X$. These subsets are called the {\em cells} of the partition. 
\item A partition is $P$ of $X$ is said to be {\em large} if $|P|= X$.

\item If $\mathcal E\subseteq\mathcal P(X)$ and $\Pi$ a family of partitions of $X$, then the pair $(\Pi,\mathcal E)$ is {\em consistent} if whenever $E\in\mathcal E$ and $C_1,\dots, C_n$ are cells chosen from {\em distinct partitions} $P_1,\dots, P_n\in\Pi$, we have $E\cap\bigcap_{i\leq n}C_i\neq \varnothing$.
\item
A family $\mathcal E\subseteq\mathcal P(X)$ is called a {\em $\pi$-system} if it is closed under finite intersections.
\end{enumerate}
\endbox\end{definition}

\begin{lemma}\label{lemma_large_consistent} Suppose that $\kappa$ is an infinite cardinal. Let $\mathcal E\subseteq\mathcal P^*(\kappa)$ be a family such that $|\mathcal E|\leq\kappa$. Then there exists a family $\Pi$ of large partitions of $\kappa$ such that $|\Pi|=2^\kappa$, and such that $(\Pi,\mathcal E)$ is consistent.
\end{lemma}
\bproof
Suppose that $\{E_\alpha:\alpha<\kappa\}$ is an enumeration of $\mathcal E$. By Corollary \ref{cor_dominate_disjoint}, there are mutually disjoint sets $I_\alpha,\alpha<\kappa$ such that $I_\alpha\subseteq E_\alpha$ and such that $|I_\alpha|=\kappa$.

Let $\Gamma:=\{(s,r): s\in\mathcal P^{<\omega}(\kappa), r:\mathcal P(s)\to\kappa\}$, and observe that $|\Gamma|=\kappa$. Since $|I_\alpha|=\kappa$, we can enumerate $\Gamma$ along each $I_\alpha$, i.e. there is an enumeration $\{(s_\xi,r_\xi):\xi<\kappa\}$ (with repetitions!) of $\Gamma$ such that 
\[\Gamma=\{(s_\xi,r_\xi):\xi\in I_\alpha\} \qquad \text{ for each }\alpha<\kappa.\]

For each $J\in\mathcal P(\kappa)$, define
\[f_J:\kappa\to\kappa:\xi\mapsto\left\{\aligned r_\xi(s_\xi\cap J)\quad&\text{if }\xi\in \bigcup_{\alpha<\kappa}I_\alpha,\\
0\qquad&\text{otherwise}.\endaligned\right..\] 

We now show that for  each sequence $\alpha,\gamma_1,\dots,\gamma_n\in\kappa$,  and each sequence $J_1,\dots ,J_n$ of {\em distinct} subsets of $ \kappa$  \[\text{there exists $\xi\in I_\alpha$ such that}\quad f_{J_i}(\xi)=\gamma_i,\qquad i=1,\dots, n.\tag{$\star$} \]
For since $J_1,\dots,J_n$ are distinct subsets of $\kappa$, the symmetric differences $J_i\Delta J_j$ are non-empty, and thus we may choose an element $x_{ij}\in J_i\Delta J_j$ for each pair $i,j$ with $1\leq 1 <j\leq n$. Define $s:=\{x_{ij}; 1\leq i<j\leq n\}$, and observe that the sets $s\cap J_i$ are distinct (for $i=1,\dots, n$). Now let $r:\mathcal P(s)\to\kappa$ be  any map such that $r(s\cap J_i)=\gamma_i$. Then $(s,r)\in\Gamma$, so there is $\xi\in I_\alpha$ so that $(s_\xi,r_\xi)=(s,r)$ Then $f_{J_i}(\xi)=\gamma_i$  for all $i=1,\dots, n$, as required.

Observe that each  $f_J:\kappa\to\kappa$ is surjective, as can be seen from $(\star)$ with $n=1$. Thus if we define $P_J$ to be the partition of $\kappa$ given  by \[P_J:=\{f_{J}^{-1}\{\gamma\}:\gamma<\kappa\},\]
then $|P_J|=\kappa$, i.e $P_J$ is large. 

Let $\Pi:=\{P_J:J\subseteq\kappa\}$. To show that $(\Pi,\mathcal E)$ is consistent, suppose that  $E\in\mathcal E$ and that $C_1,\dots, C_n$ belong to $P_{J_1},\dots, P_{J_n}$, where the $J_i$ are distinct subsets of $\kappa$.  Then there exists $\alpha<\kappa$ such that $E=E_\alpha$, and there exist $\gamma_1,\dots,\gamma_n<\kappa$ so that $C_i=f_{J_i}^{-1}\{\gamma_i\}$ for $ i=1,\dots, n$. But then by $(\star)$ there exists $\xi\in I_\alpha$ such that $f_{J_i}(\xi)=\gamma_i$, from which it follows that  \[\xi\in I_\alpha\cap\bigcap_{i\leq n} C_i.\] Since $I_\alpha\subseteq E_\alpha$, it follows that $E_\alpha\cap \bigcap_{i\leq n}C_i\neq \varnothing$. Hence $(\Pi,\mathcal E)$ is consistent.

In particular, if $J_1, J_2$ are distinct subsets of $\kappa$, then $C_1\cap C_2\neq \varnothing$ for any $C_1\in P_{J_1}, C_2\in P_{J_2}$. It follows that $P_{J_1}, P_{J_2}$ have no cells in common, and thus that $P_{J_1}\neq P_{J_2}$, i.e. all the $P_J$ are distinct. Hence $|\Pi|=|\mathcal P(\kappa)|=2^\kappa$.
\eproof


\begin{lemma}\label{lemma_successor_step_1}
Suppose that $\Pi$ is a family of partitions of a cardinal $\kappa$ and that $\mathcal E\subseteq\mathcal P(\kappa)$ is  a $\pi$-system. Suppose further that $(\Pi,\mathcal E)$ is consistent. Let $J\subseteq \kappa$. Then either $(\Pi, \mathcal E\cap\{J\})$ is consistent, or else there is a cofinite $\Pi'\subseteq\Pi$ so that $(\Pi',\mathcal E\cap\{\kappa-J\})$ is consistent. \end{lemma}

\bproof If $(\Pi, \mathcal E\cap\{J\})$ is not consistent, there is a set $E\in\mathcal E$,  and  distinct partitions $P_1,\dots, P_n\in \Pi$ and corresponding cells $C_i\in P_i$ such that
\[J\cap E\cap\bigcap_{i\leq n}C_i=\varnothing.\] It follows that $E\cap\bigcap_{i\leq n}C_i\subseteq \kappa-J$.  Let $\Pi':=\Pi-\{P_1,\dots, P_n\}$. To prove that $(\Pi',\mathcal E\cup\{\kappa-J\})$ is consistent, take  $E'\in\mathcal E$, distinct partitions $P_1',\dots, P_m'\in\Pi'$, and cells $C_1',\dots, C_m'$ in $P_1',\dots, P_m'$. We need only show $(\kappa-J)\cap E'\cap\bigcap_{j=1}C_j'\neq \varnothing$.  Now as $(\Pi,\mathcal E)$ is consistent and $\mathcal E$ is closed under finite intersections, we have\[E\cap E'\cap\bigcap_{i\leq n}C_i\cap\bigcap_{j\leq m}C_j'\neq \varnothing.\] As $E\cap\bigcap_{i\leq n}C_i\subseteq \kappa-J$, we immediately
see that also $(\kappa-J)\cap E'\cap\bigcap_{j\leq m}C_j',\neq \varnothing$. 
\eproof

\begin{lemma}\label{lemma_successor_step_2}
Suppose that $\Pi$ is a family of large partitions of a cardinal $\kappa$ and that $\mathcal E\subseteq\mathcal P(\kappa)$ is a $\pi$-system. Suppose further that $(\Pi,\mathcal E)$ is consistent. Let $p:\mathcal P^{<\omega}(\kappa)\to\mathcal E$ be an order-reversal, and let $P\in\Pi$. Then there is a $\pi$-system $\mathcal E' \supseteq\mathcal E$ and an anti-additive  map $q:\mathcal P^{<\omega}(\kappa)\to\mathcal E'$ such that $q\leq p$ and $(\Pi-\{P\},\mathcal E')$ is consistent.
\newline In addition, we can take $\mathcal E'$ to be the $\pi$-system generated by $\mathcal E\cup\ran(q)$.
\end{lemma}

\bproof
Let $\{C_\alpha: \alpha<\kappa\}$ be an enumeration of the cells of $P$, without repetition, and let $\{t_\alpha:\alpha<\kappa\}$ be an enumeration of $\mathcal P^{<\omega}(\kappa)$. For each $\alpha<\kappa$ define a map $q_\alpha:\mathcal P^{<\omega}(\kappa)\to\mathcal P(\kappa)$ by\[q_\alpha(s):=\left\{\aligned p(t_\alpha)\cap  C_\alpha\quad&\text{if }s\subseteq t_\alpha,\\
\varnothing\qquad&\text{else.}\endaligned\right.\]
Since $p$ is an order-reversal, $q_\alpha(s)\subseteq p(t_\alpha)\subseteq p(s)$. Furthermore, since $(\Pi,\mathcal E)$ is consistent, $p(t_\alpha)\in\mathcal E$ and $P\in \Pi$, we always have $p(t_\alpha)\cap C_\alpha\neq \varnothing$. Hence if $s\subseteq t_\alpha$, then $q_\alpha(s)\neq \varnothing$.

Next observe that each $q_\alpha$ is  anti-additive: Indeed, if both $s,s'\subseteq t_\alpha$, then clearly $q_\alpha(s\cup s')=p(t_\alpha)\cap C_\alpha = q_\alpha(s)\cap q_\alpha(s')$. On the other hand, if one of $s,s'$ is  not contained in  $t_\alpha$, then $q_\alpha(s\cap s')=\varnothing = q_\alpha(s)\cap q_\alpha(s')$. 

Define a function $q$ on $\mathcal P^{<\omega}(\kappa)$ by 
\[q(s):=\bigcup_{\alpha < \kappa} q_\alpha(s).\]
Then since each $q_\alpha(s)\subseteq p(s)$ we have $q\leq p$. Further note that if $s,s'\in \mathcal P^{<\omega}(\kappa)$ (not necessarily distinct), and $\alpha\neq \beta$, then $q_\alpha(s)\cap q_\beta(s')=\varnothing$, since $C_\alpha\cap C_\beta=\varnothing$. In particular, $q(s)$ is a disjoint union of subsets of the cells $C_\alpha$, and hence $q$ is also anti-additive:
\[\aligned q(s)\cap q(s')&=\bigcup_{\alpha<\kappa}q_\alpha(s) \cap \bigcup_{\beta<\kappa}q_\beta(s'),\\
&=\bigcup_{\alpha<\kappa}\bigcup_{\beta<\kappa} \big(q_\alpha(s)\cap q_\beta(s')\big),\\
&= \bigcup_{\alpha<\kappa}\big(q_\alpha(s)\cap q_\alpha(s')\big),\\
&=\bigcup_{\alpha<\kappa}q_\alpha(s\cup s'),\\
&=q(s\cup s').
\endaligned\]
Hence $\ran(q)$ is a $\pi$-system. Define $\mathcal E'$ to be the $\pi$-system generated by $\mathcal E\cup\ran(q)$:
\[\mathcal E':=\{E\cap q(s):E\in\mathcal E, s\in \mathcal P^{<\omega}(\kappa)\},\] and note that $q:\mathcal P^{<\omega}(\kappa)\to\mathcal E'$.  We claim that $(\Pi-\{P\}, \mathcal E')$ is consistent. For suppose that, for $i=1,\dots,n$ we have $D_i\in P_i\in\Pi-\{P\}$, and that $E\in\mathcal E, s\in \mathcal P^{<\omega}(\kappa)$. Then $s=t_\alpha$ for some $\alpha<\kappa$, and hence $q_\alpha(s)=p(t_\alpha)\cap C_\alpha=p(s)\cap C_\alpha$. Now since $p:\mathcal P^{<\omega}(\kappa)\to\mathcal E$ and $\mathcal E$ is a $\pi$-system, we have $p(s)\cap E\in\mathcal E$. Since $(\Pi,\mathcal E)$ is consistent, \[
p(s)\cap E\cap C_\alpha\cap\bigcap_{i\leq n}D_i\neq \varnothing,\] and hence \[ \big(E\cap q(s)\big)\cap \bigcap_{i\leq n}D_i\supseteq\big(E\cap q_\alpha(s)\big)\cap \bigcap_{i\leq n}D_i=p(s)\cap E\cap C_\alpha\cap\bigcap_{i\leq n}D_i\varnothing.\]
\eproof

We are now ready to prove the existence of good ultrafilters under certain conditions:

\begin{theorem}\label{thm_good_ultrafilters_exist} Suppose that $\kappa$ is an infinite cardinal. Then there exists a $\kappa^+$-good countably incomplete ultrafilter over $\kappa$.
\end{theorem}

\bproof
Start with a sequence $I_n\downarrow\varnothing$ of subsets of $\kappa$, each of cardinality $\kappa$. (This exists since $|\omega\times\kappa|=\kappa$.) Let $\mathcal F_0$ be the filter over $\kappa$ generated by the sets $I_n$. Then any filter $\mathcal F$ which extends $\mathcal F_0$ will be countably incomplete. 

By Lemma \ref{lemma_large_consistent} there is a family $\Pi_0$ of large partitions of $\kappa$ such that $|\Pi_0|=2^\kappa$, and such that $(\Pi_0,\mathcal F_0)$ is consistent. We now use transfinite induction to define, for  ordinals $\xi<2^\kappa$, a  sequence $\Pi_\xi$ of partitions of $\kappa$, and a sequence $\mathcal F_\xi$ of filters over $\kappa$ such that \begin{enumerate}[(i)]
\item If $\eta\leq\xi<2^\kappa$, then $\Pi_\eta\supseteq \Pi_\xi$ and $\mathcal F_\eta\subseteq \mathcal F_\xi$.
\item $|\Pi_\xi|=2^\kappa$.
\item $|\Pi_{\xi+1}-\Pi_{\xi}|<\omega$.
\item $\Pi_\lambda=\bigcap_{\xi<\lambda}\Pi_\xi$ if $\lambda<2^\kappa$ is a limit ordinal.
\item Each pair $(\Pi_\xi, \mathcal F_\xi)$ is consistent.
\end{enumerate}

Suppose that $\{p_\xi: \xi<2^\kappa\}$ is an enumeration of all order-reversing maps $\mathcal P^{<\omega}(\kappa)\to \mathcal P(\kappa)$. Similarly, let $\{J_\xi:\xi<2^\kappa\}$ enumerate $\mathcal P(\kappa)$.
Assume now that we have defined $(\Pi_\eta,\mathcal F_\eta)$ for all $\eta<\xi$, so that (i)-(iv) are satisfied.

There are three cases in the construction of the $\kappa^+$-good ultrafilter $\mathcal U$ by transfinite induction. Cases 1 and 2 deal with odd and even successor ordinals. Case 1 ensures that the filter we construct is an ultrafilter, whereas Case 2 ensures that it is good. Case 3, which is the simplest, deals with the step at limit ordinals.

\underline{Case 1:} $\xi$ is a successor ordinal of the form $\lambda +2n+1$, where $\lambda$ is a limit ordinal and $n<\omega$.\\
Let $J_\eta$ be the first subset of $\kappa$ in the enumeration of $\mathcal P(\kappa)$ which is not in $\mathcal F_{\xi-1}$. By Lemma \ref{lemma_successor_step_1}, we can find a partition $\Pi_\xi\subseteq\Pi_{\xi-1}$ such that $|\Pi_{\xi-1}-\Pi_\xi|<\omega$ (so that $|\Pi_\xi|=2^\kappa$ also), such that either $(\Pi_\xi, \mathcal F_{\xi-1}\cup\{J_\eta\})$ is consistent, or $(\Pi_\xi, \mathcal F_{\xi-1}\cup\{\kappa-J_\eta\})$ is consistent. In the former case, define $\mathcal F_\xi$ to be the filter generated by $\mathcal F_{\xi-1}\cup\{J_\eta\}$, and in the latter, the filter generated by $\mathcal F_{\xi-1}\cup\{\kappa-J_\eta\}$. Thus 
$(\Pi_\xi,\mathcal F_\xi)$ is consistent.

\underline{Case 2:} $\xi$ is a successor ordinal of the form $\lambda +2n+2$, where $\lambda$ is a limit ordinal.\\
In that case let $p_\eta$ be the first function $\mathcal P^{<\omega}(\kappa)\to \mathcal F_{\xi-1}$ in the enumeration of  order-reversing maps $\mathcal P^{<\omega}(\kappa)\to \mathcal P(\kappa)$ that has not already been dealt with. Pick $P\in\Pi_{\xi-1}$. By Lemma \ref{lemma_successor_step_2}, there is an anti-additive $q:\mathcal P^{<\omega}(\kappa)\to \mathcal F_\xi$ such that $q\leq p_\eta$, and such that $(\Pi_\xi,\mathcal F_\xi)$ is consistent, where $\Pi_\xi:=\Pi_{\xi-1}-\{P\}$ and $\mathcal F_\xi$ is the filter generated by $\mathcal F_{\xi-1}\cup\ran(q)$.

\underline{Case 3:} $\xi$ is a limit ordinal.\\
In that case, define $\Pi_\xi:=\bigcap_{\eta<\xi} \Pi_\eta$ and $\mathcal F_\xi:=\bigcup_{\eta<\xi}\mathcal F_\eta$. Since at each stage $\eta<\xi$ we remove only finitely many partitions, at most $|\xi\cdot\omega|<2^\kappa$--many partitions have been removed from $\Pi_0$ to form $\Pi_\xi$, and hence $|\Pi_\xi|=2^\kappa$. It is easy to see that $(\Pi_\xi,\mathcal F_\xi)$ is consistent: For suppose that $F\in\mathcal F_\xi$ and $C_1,\dots, C_m$ are cells from distinct partitions $P_1,\dots, P_m\in \Pi_\xi$. Then there is $\eta<\xi$ such that $F\in \mathcal F_\eta$, and moreover $C_1,\dots, C_m$ are cells from distinct partitions in $\Pi_\eta$, as $\Pi_\xi\subseteq\Pi_\eta$. Since $(\Pi_\eta,\mathcal F_\eta)$ is consistent, $F\cap\bigcap_{i\leq m}C_i\neq \varnothing$.

This completes the transfinite induction. 

Now define $\mathcal U:=\bigcup_{\xi<2^\kappa}\mathcal F_\xi$.  By Case 1, for every $J\subseteq \kappa$, either $J\in\mathcal U$ or $\kappa-J\in\mathcal U$. Hence $\mathcal U$ is an ultrafilter. Since $\mathcal U\supseteq\mathcal F_0$ it is a countably incomplete ultrafilter.

Suppose now that $p:\mathcal P^{<\omega}(\kappa)\to\mathcal U$ is an order-reversal.  Observe\footnote{Recall K\"onig's Theorem: If $\alpha_i<\beta_i$ are cardinals, then $\sum_i\alpha_i<\prod_i\beta_i$. Now if $\alpha_i<2^\kappa$, and $\beta_i=2^\kappa$ for $i<\kappa$, then $\sum_{i<\kappa}\alpha_i<\prod_{i<\kappa}2^\kappa=2^\kappa$, i.e. $\sum_{i<\kappa}\alpha_i<2^\kappa$. Hence $\text{cf}(2^\kappa)>\kappa$.} that since $\text{cf}(2^\kappa)>\kappa$ and $|\mathcal P^{<\omega}(\kappa)|=\kappa$, there is a least $\eta<\kappa$ so that $\ran(p)\subseteq\mathcal F_\eta$. Thus $p$ will be dealt with at some stage $\xi\geq \eta$ by Case 2, which guarantees the existence of an anti-additive $q\leq p$ which maps into $\mathcal F_\xi\subseteq\mathcal U$. It follows that $\mathcal U$ is $\kappa^+$-good.
\eproof

\bibliographystyle{alpha}
\bibliography{Nonstandard_Universes_bib.org.tug}
\end{document}